אוניברסיטת בן-גוריון בנגב
Ben-Gurion University of the Negev

FACULTY OF ENGINEERING SCIENCES

# HIGH-EFFICIENCY SELF-ADJUSTING SWITCHED CAPACITOR DC-DC CONVERTER WITH BINARY RESOLUTION

by
Alexander Kushnerov

A thesis submitted to
the Department of Electrical Engineering and Computer Science
in partial fulfillment of the requirements
for the degree of
Master of Science

Supervised by: Professor Sam Ben-Yaakov
and Professor Eugene Paperno

Beer-Sheva, August 2009

**High-Efficiency Self-Adjusting Switched Capacitor DC-DC Converter
with Binary Resolution**

by

Alexander Kushnerov

Dipl. Ing. (Omsk State Technical University) 2002

A thesis submitted in partial fulfillment
of the requirements for the degree of

Master

in

Electrical Engineering and Computer Science

in the

GRADUATE DIVISION

of the

BEN-GURION UNIVERSITY OF THE NEGEV

Jury:

Professor Sam Ben-Yaakov
Professor Eugene Paperno
Professor Raul Rabinovici
Doctor Doron Shmilovitz

Defense date:
March 04, 2010

# ABSTRACT

Switched-Capacitor Converters (SCC) suffer from a fundamental power loss deficiency which make their use in some applications prohibitive. The power loss is due to the inherent energy dissipation when SCC operate between or outside their output target voltages. This drawback was alleviated in this work by developing two new classes of SCC providing binary and arbitrary resolution of closely spaced target voltages. Special attention is paid to SCC topologies of binary resolution. Namely, SCC systems that can be configured to have a no-load output to input voltage ratio that is equal to any binary fraction for a given number of bits.

To this end, we define a new number system and develop rules to translate these numbers into SCC hardware that follows the algebraic behavior. According to this approach, the flying capacitors are automatically kept charged to binary weighted voltages and consequently the resolution of the target voltages follows a binary number representation and can be made higher by increasing the number of capacitors (bits). The ability to increase the number of target voltages reduces the spacing between them and, consequently, increases the efficiency when the input varies over a large voltage range.

The thesis presents the underlining theory of the binary SCC and its extension to the general radix case. Although the major application is in step-down SCC, a simple method to utilize these SCC for step-up conversion is also described, as well as a method to reduce the output voltage ripple. In addition, the generic and unified model is strictly applied to derive the SCC equivalent resistor, which is a measure of the power loss. The theoretical predictions are verified by simulation and experimental results.

# ACKNOWLEDGEMENTS

First and foremost, I would like to thank my advisor Professor Sam Ben-Yaakov for his patient guidance. He gave me a chance to start my education here and provided a very good environment for my research and study. Besides his outstanding expertise in both theoretical and practical matters, his amicable disposition and accessibility have provided for constructive and fruitful work. Professor Ben-Yaakov has taught me how to structure my ideas more rigorously and I believe over these past years that I have absorbed some of his creative approach to research.

This thesis could not have been developed without the pioneering work and continual support of Meir Shashoua, co-founder of "K. S. Waves Ltd.", Tel-Aviv, Israel. I am grateful to Mr. Shashoua for his original ideas and asking questions that have periodically resurfaced in my mind and conducted me to sharper thinking. The material assistance from "K. S. Waves Ltd." is deeply appreciated and deserves special acknowledgement.

I am greatly indebted to the graduate student Michael Evzelman for his friendship and optimism. He is always ready to discuss anything that could be connected with motorcycles, electronics or programming. My gratitude is also to Alexei Smirnov, the author of the NL5 circuit simulator, for his careful reading of the manuscript and help with English corrections.

At the end of this study I wish to thank my mother – her love, understanding, support and sacrifices are always an encouragement to me.

# TABLE OF CONTENTS

# TABLE OF CONTENTS (cont'd)

# LIST OF FIGURES

# LIST OF FIGURES (cont'd)

# 1. INTRODUCTION

## 1.1 Background review

The purpose of a DC-DC converter is to provide a predetermined and constant output voltage to a load from a poorly specified or fluctuating input voltage source. Linear regulators and switching converters are two common types of DC-DC converters. In a linear regulator the output current comes directly from the power supply, therefore the efficiency is approximately defined as the ratio of the output voltage to supply voltage. It is obvious that a worse efficiency will be obtained when the supply voltage is much larger than the output voltage. Switching converters are more efficient than linear regulators due to intercepted energy transfer. This is done by periodically switching energy storing components to deliver a portion of energy from the power supply to the output. Switching DC-DC converters (except for resonant converters) can be divided into two large groups: inductive and capacitive.

The inductive converters using one or several inductors have been a power supply solution in all kinds of applications for many years due to the wide variety of possibilities in current and voltage requirements. Generally, the inductors in such a converter are bulky, not realizable on-chip and are the cause of two difficult problems. One problem are high voltage spikes that must be damped or recuperated otherwise, the switches which are not rated for such constraints can blow, while the rest of circuit can be damaged. The other problem with inductive converters is a pulsating input current, which can produce an electromagnetic interference (EMI) from other equipment and conductor lines. This interference may penetrate into susceptible devices and lead to unreliable operation. So, the pulsating input current requires a special filter and sometimes shielding. All these factors increase the board space and inductive converter cost.

The capacitive converters based on switched capacitors are widespread in applications requiring small power and no isolation between input and output. They feature relatively low noise, minimal radiated EMI, and in most cases are fabricated as integrated circuits which have made capacitive converters popular for use in power management for mobile devices. An additional goal of such converters is the option for unloaded operation with no need for dummy loads or complex control. However, capacitive converters suffer from inherent power loss during charging and discharging of a capacitor connected in parallel with the voltage source or another capacitor. Theory predicts that this power loss is proportional to the squared voltage difference taking place before the corresponding circuit has been configured. As a result, capacitive converters exhibit a rather high efficiency if the capacitors pre-charged to certain voltages are paralleled with components maintaining similar voltages.

The most known type of capacitive DC-DC converter is called a charge pump; for historical reasons it is often considered to a step-up converter built from capacitors and diodes, which are used as switches. Nowadays, when charge pumps are built around transistor switches, their circuitry does not differ in principle from the step-down switched capacitor DC-DC converters. The cornerstone of both circuits is a reconfigurable array of switches and capacitors generally called "flying capacitors". These capacitors are charged from the input voltage and then discharged to the load thus providing charge transfer and a constant output voltage.

It is a well-known phenomenon that when a capacitive converter operates at the target output to input voltage ratios, the efficiency is high and may exceed 90%. This is due to the fact that, at these voltage ratios, the capacitors do not see appreciable voltage variations. When the same capacitive converter operates between or outside the target voltage ratios, the efficiency drops dramatically. Obviously, in practice one would expect the conversion ratio to change and hence there is no way to escape the losses. However, there are several "lossless" techniques to provide regulation of the output voltage. In most cases, these techniques change the rate at which the charge is transferred to the output and this leads to an increased output voltage ripple. In general, capacitive converters feature a set of discrete target voltage ratios that can be contrasted with the continuous transfer function of inductive converters.

The down side of capacitive converters is the larger number of switches and respective drivers complicating the converter circuitry. Another problem of capacitive converters is a high inrush current during start-up that must be limited by soft-start circuitry.

## 1.2. Motivation and relevance

*"Discontent is the first necessity of progress."*
Thomas A. Edison

Switched Capacitor Converters (SCC) suffer from a fundamental power loss which is a severe limitation because of the common requirement to regulate output voltage. The power loss is due to the inherent energy dissipation when a capacitor is charged or discharged by a voltage source or another capacitor [1-8]. Hence, SCC exhibit rather high efficiency only when operating at the target voltages at which the voltage differences that charge and discharge the capacitors are small. Earlier studies attempted to overcome the power loss by proposing SCC with an increased number of target voltages [9-14]. However, the common disadvantage of these SCC is that the target voltages are spread apart.

It is thus evident that there is a need and it will be highly advantageous to design a SCC that has a large number of target voltages that are spaced at high resolution over the range of interest and thereby improve the efficiency. Another desired feature is a simple way to increase the resolution only by changing the control scheme. In addition, it would be desirable to obtain a smooth transition from one target voltage ratio to another. It is yet another demand to regulate the output voltage while maintaining high efficiency. It would also be desirable to provide low output voltage ripple over a wide range of target voltage ratios.

This work presents the theory that underlines the operation of the multi-target SCC and allows one to design new SCC satisfying the above requirements. The theory is based on the redundancy of the positional number systems [38-42], which is used to develop two new SCC classes providing binary and arbitrary resolution of target voltages. In these new SCC classes, the flying capacitors are automatically kept charged to radix-$r$-weighted voltages, while the gap between neighboring target voltages is defined by the resolution. Both the radix and the resolution can be made higher by increasing the number of flying capacitors.

# 2. BASICS OF SWITCHED CAPACITOR CIRCUITS

## 2.1 Transient and limitation of current spike

*"The immediate effect is likely to be what it's always been - a spike in violence."*
Donald Rumsfeld

The principle of a gradual change of energy in any physical system, and specifically in an electrical circuit, means that the energy stored in electric or magnetic fields cannot change instantaneously [1-3]. For the sake of simplicity, however, the assumption is made in transient analysis that the switching occurs quite instantaneously [4-8].

Let $t = 0$ be the instant of time when switching starts, and two additional instants: just prior and just after switching be $t = 0_-$ and $t = 0_+$ respectively. In mathematical language, the value of the function $f(0_-)$ is the "limit from the left", as $t$ approaches zero from the left, while $f(0_+)$ is the "limit from the right", as $t$ approaches zero from the right. According to the above principle, the voltage (charge) of a capacitor just after switching is equal to the voltage (charge) just prior to switching:

$$v_C(0_+) = v_C(0_-) \tag{2.1.1}$$

$$q(0_+) = q(0_-) \tag{2.1.2}$$

Defining an ideal switch as a zero-resistance device that gets opened or closed in zero time, we consider the charging circuit shown in Fig. 2.1.1(a), where the voltage source $V_S$, the switch $Sw$ and the capacitor $C_1$ are ideal. When $Sw$ is turned-on, the capacitor voltage $v_1$ changes abruptly from zero to $V_S$. In other words, the charging of $C_1$ is accompanied by an infinitely high current pulse during an infinitesimal time.

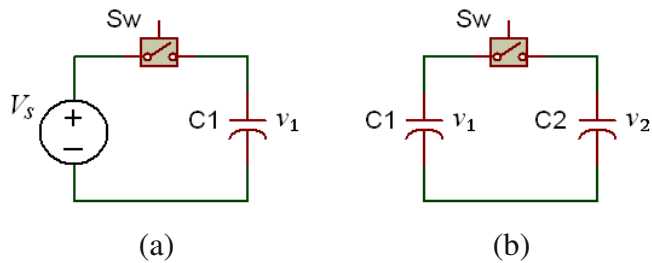

(a)         (b)

Figure 2.1.1: Switched circuits including the ideal capacitors.

Using the above designations, we can write $v_1(0_-) = 0$, $v_1(0_+) = V_S$ and $v_1(0_+) \neq v_1(0_-)$ which contradicts (2.1.1). In transient analysis, the last expression is called an incorrect initial condition for the chosen mathematical model of an ideal switched circuit.

Consider now the switched circuit of Fig. 2.1.1(b), where the ideal switch $Sw$ serves to discharge the ideal capacitor $C_1$ pre-charged to the voltage $v_1(0_-) = V_S$ into another empty ideal capacitor $C_2$. According to the law of charge conservation, the total charge on two capacitors $C_1$ and $C_2$ connected in parallel is the sum of the initial charges $q_1(0_-) = C_1 V_S$ and $q_2(0_-) = 0$, while the final voltages are:

$$v_1(0_+) = v_2(0_+) = \frac{C_1}{C_1 + C_2} V_S \tag{2.1.3}$$

So, the contradiction of $v_C(0_+) \neq v_C(0_-)$ is observed again and, as in the previous case, the discharging of $C_1$ will be accompanied by an infinitely high current pulse during an infinitesimal time. This contradiction can be refuted since any circuit with a real capacitor has in practice some resistance and inductance connected in series. The series inductance is generally and is neglected in present analysis.

The infinite current spike is prevented in the switched circuits shown in Fig. 2.1.2, so that the initial conditions are correct. Note that such circuits may be composed by taking into consideration just the resistances of the connecting wires.

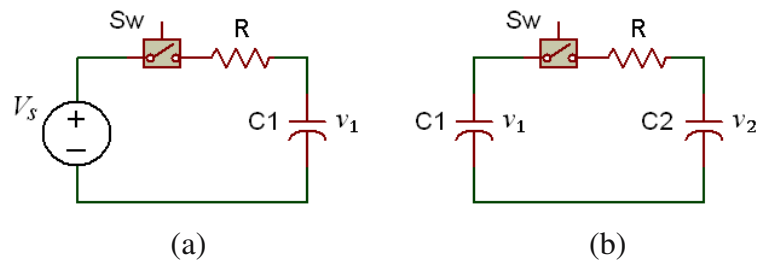

(a)                  (b)

Figure 2.1.2: Switched circuits including the serial resistor.

The charging circuit with the resistor is shown in Fig. 2.1.2(a), it is described by a first order differential equation because it comprises only one capacitor $C_1$. So, the aim is to calculate the complete response of the first order circuit to the voltage step $V_S$. According to the Kirchhoff Current Law $i_R(t) = i_{C1}(t)$, this is:

$$\frac{V_S - v_1(t)}{R} = C_1 \frac{dv_1(t)}{dt} \tag{2.1.4}$$

Rearranging the above equation:

$$\frac{dv_1(t)}{v_1(t) - V_S} = -\frac{1}{RC_1} dt \qquad (2.1.5)$$

Now it can be simply integrated:

$$\int \frac{dv_1(t)}{v_1(t) - V_S} = -\frac{1}{RC_1} \int dt + D \qquad (2.1.6)$$

The integration of both sides yields:

$$\ln[v_1(t) - V_S] = -\frac{t}{RC_1} + D \qquad (2.1.7)$$

Since the time constant $\tau = RC_1$,

$$v_1(t) - V_S = e^{D - t/\tau} \qquad (2.1.8)$$

An initial voltage across $C_1$ will be

$$v_1(0) = V_0 = e^D + V_S \qquad (2.1.9)$$

So, the complete response is:

$$v_1(t) = V_0 e^{-t/\tau} + V_S (1 - e^{-t/\tau}) \qquad (2.1.10)$$

Note that the first term in (2.1.10) is the natural response, while the second term is the forced response. Both terms and the complete response were calculated in MathCAD and are presented in Fig. 2.1.3 together with the following current, which is limited by $I_0 = (V_S - V_0)/R$.

The time constant may be easily found from Fig. 2.1.3 by drawing a tangent line to the response curve at $t = 0$. The intercept point of the tangent and the asymptotic limit projected to the time axis yields the time constant. The units of the time constant are seconds [$\tau$] = $\Omega \cdot$F, therefore it is considered as an interval during which the voltage drops (grows) relatively to its initial value. At the end of the $\tau$ interval, the voltage is $e^{-1} \approx 0.368$ of its initial value, while at the end of $5\tau$ the voltage ratio is less than 0.01. Because of this fact, it is usual to presume that the duration of the transient response is about $5\tau$. Note that, precisely speaking, the transient response declines to zero in infinite time, since $e^{-t} \to 0$, when $t \to \infty$.

$$V_S := 1 \quad V_0 := \sqrt{0.5} \cdot V_S \quad R := 1 \quad C_1 := 10^{-6} \quad \tau := R \cdot C_1 \quad n := 5 \cdot \tau \quad j := \frac{n}{10^4} \quad t := 0, j..n$$

$$up(t) := V_S \cdot \left(1 - \exp\left(\frac{-t}{\tau}\right)\right) \quad dn(t) := V_0 \cdot \exp\left(\frac{-t}{\tau}\right) \quad v_1(t) := up(t) + dn(t) \quad i(t) := \frac{V_S - V_0}{R} \cdot \exp\left(\frac{-t}{\tau}\right)$$

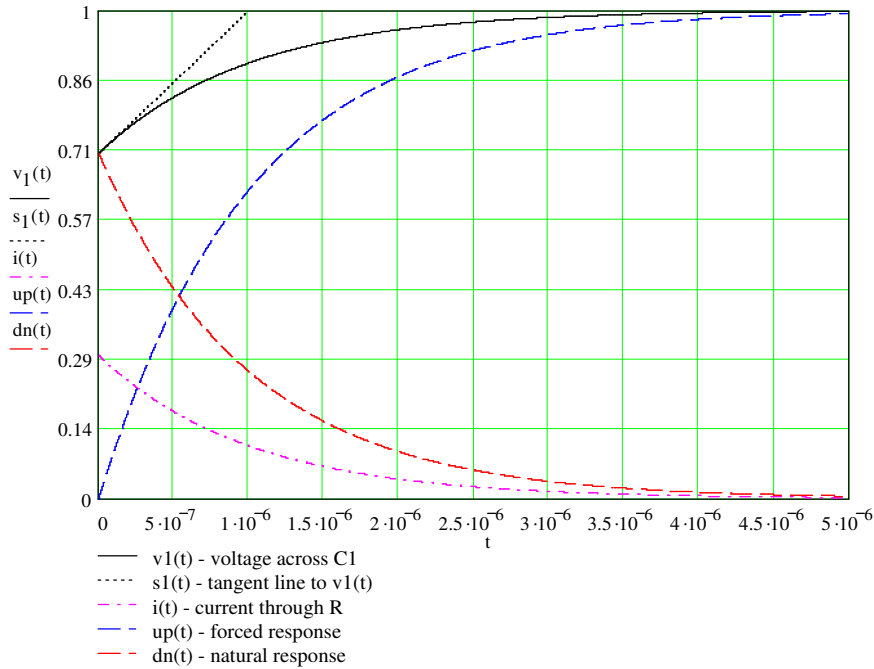

v1(t) - voltage across C1
s1(t) - tangent line to v1(t)
i(t) - current through R
up(t) - forced response
dn(t) - natural response

Figure 2.1.3: Complete response of the charging circuit and its components.

When $t \gg 5\tau$, the voltage $v_1(t)$ is considered to be equal to the voltage $v_1(0_+) = V_S$ as in the ideal switched circuit shown in Fig. 2.1.1(a). The amount of charge transferred by an exponentially decaying current is equal to the product of its initial value and the time constant.

$$q = \int_0^\infty i(t)dt = I_0 \int_0^\infty e^{-t/\tau} dt = I_0(-t)e^{-t/\tau} \Big|_0^\infty = I_0 \tau \tag{2.1.11}$$

This result justifies using an impulse function $\delta$ to represent the very large, approaching infinity, magnitude of the current pulse, applied for a very short (approaching zero) time interval, whereas their product stays finite, as shown in Fig. 2.1.4.

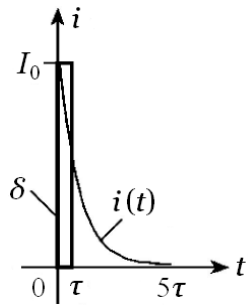

Figure 2.1.4: A large and fast decaying $i(t)$ and an equivalent impulse.

The other switched circuit with the resistor $R$ is shown in Fig. 2.1.2(b) and serves to discharge the capacitor $C_1$ (pre-charged to the voltage $V_S$) and simultaneously to charge the empty capacitor $C_2$. To find the voltages across these capacitors, consider a bi-directional current flow and compose two first-order differential equations:

$$\begin{cases} v_1(t) = v_2(t) - RC_1 \dfrac{dv_1(t)}{dt} \\ v_1(t) = v_2(t) + RC_2 \dfrac{dv_2(t)}{dt} \end{cases} \quad \text{or} \quad \begin{cases} v_2(t) = v_1(t) + RC_1 \dfrac{dv_1(t)}{dt} \\ v_2(t) = v_1(t) - RC_2 \dfrac{dv_2(t)}{dt} \end{cases} \quad (2.1.12)$$

Take the Laplace Transform of both systems:

$$\begin{cases} v_1(s) = v_2(s) - RC_1[s \cdot v_1(s) - V_S] \\ v_1(s) = v_2(s) \cdot [1 + s \cdot RC_2] \end{cases} \quad \begin{cases} v_2(s) = v_1(s) + RC_1[s \cdot v_1(s) - V_S] \\ v_2(s) = v_1(s) - RC_2 \cdot s \cdot v_2(t) \end{cases} \quad (2.1.13)$$

The solutions in the Laplace domain are:

$$v_1(s) = \frac{C_1 V_S (1 + s \cdot RC_2)}{s^2 R \cdot C_1 C_2 + s(C_1 + C_2)} \qquad v_2(s) = \frac{C_1 V_S}{s^2 R \cdot C_1 C_2 + s(C_1 + C_2)} \quad (2.1.14)$$

Taking the Inverse Laplace Transform of both equations (2.1.14), we obtain the voltages in the time domain. To simplify the expressions we introduce the time constant $\tau = R \dfrac{C_1 C_2}{C_1 + C_2}$ since the current flows through serially connected capacitors.

$$v_1(t) = \frac{C_1 V_S}{C_1 + C_2} + \frac{C_2 V_S}{C_1 + C_2} e^{-t/\tau} \qquad v_2(t) = \frac{C_1 V_S}{C_1 + C_2}[1 - e^{-t/\tau}] \quad (2.1.15)$$

$$\text{while } i(t) = \frac{v_1(t) - v_2(t)}{R} = \frac{V_S}{R} e^{-t/\tau} \quad (2.1.16)$$

The boundary values are:

$$i(0) = I_0 = \frac{V_S}{R} \qquad V_1 = \lim_{t \to \infty} v_1(t) = \frac{C_1 V_S}{C_1 + C_2} \qquad V_2 = \lim_{t \to \infty} v_2(t) = \frac{C_1 V_S}{C_1 + C_2} \quad (2.1.17)$$

It is evident from (2.1.16) that the asymptotic limits are the same voltages $v_1(0_+) = v_2(0_+)$ as derived in (2.1.3) by using the charge conservation law. At the instant of switching, the current is limited by $i(0) = V_S / R$ and reaches 0.01 of this value at $5\tau$.

As follows from (2.1.14), the transient rates for $C_1$ and $C_2$ are different and defined by the time constants $RC_1$ and $RC_2$ respectively. Since in this particular case, $C_1$ is pre-charged to $V_S$, its discharging can be considered as the natural response, while the charging of empty $C_2$ matches the definition of a forced response. Both the voltages of (2.1.15) and the current through $R$ given by (2.1.16) were calculated in MathCAD and depicted in Fig. 2.1.5.

$$R := 1 \quad C_1 := 10^{-6} \quad C_2 := 2 \cdot C_1 \quad V_s := 1 \quad \tau := \frac{C_1 \cdot C_2}{C_1 + C_2} \cdot R \quad n := 5 \cdot \tau \quad i := \frac{n}{10^4} \quad t := 0, i..n$$

$$v_1(t) := \frac{C_1 \cdot V_s + C_2 \cdot V_s \cdot \exp\left(\frac{-t}{\tau}\right)}{C_1 + C_2} \qquad v_2(t) := \frac{C_1 \cdot V_s \cdot \left(1 - \exp\left(\frac{-t}{\tau}\right)\right)}{C_1 + C_2} \qquad i(t) := \frac{v_1(t) - v_2(t)}{R}$$

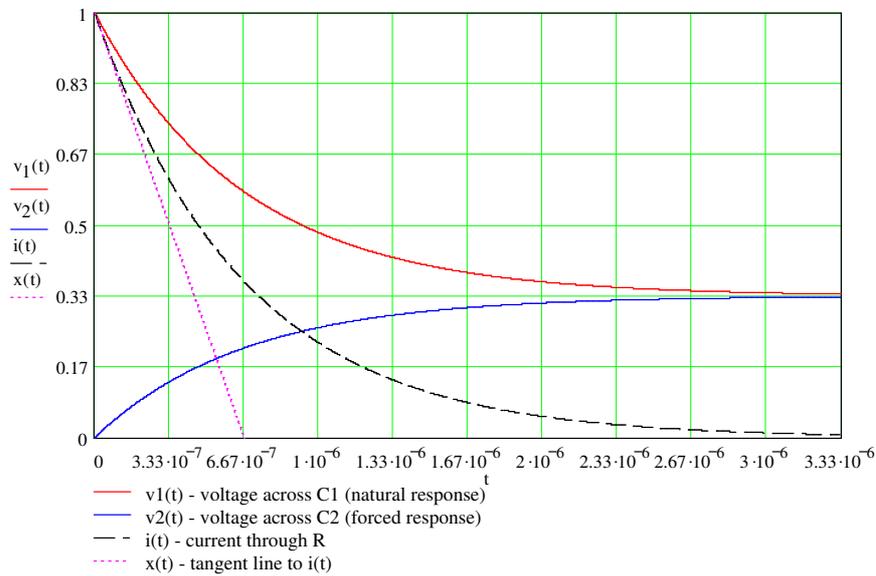

v1(t) - voltage across C1 (natural response)
v2(t) - voltage across C2 (forced response)
i(t) - current through R
x(t) - tangent line to i(t)

Figure 2.1.5: Complete response of the discharging circuit and its components.

## 2.2 Inherent energy loss at voltage difference

*"In mathematics you don't understand things. You just get used to them."*
John von Neumann

The transient in the switched circuits considered in the previous section is accompanied by either an infinitely high pulse or an exponentially decaying current. Energy is lost in both cases; however, in each case the nature of energy loss is different. In the first case of the ideal switched circuit, it is common to presume that the energy loss is radiation caused by the infinitely high current pulse. The other case is more close to practice because the current is limited by the series resistor, which is heated and dissipates energy. As known, the energy stored in the capacitor is:

$$E = \frac{CV^2}{2} = \frac{Q^2}{2C} = \frac{V}{2Q} \qquad (2.2.1)$$

Consider again the charging circuit in Fig. 2.1.1 (a), where the ideal capacitor $C_1$ is pre-charged to the voltage $V_0$ and holds an initial energy $E_0 = C_1 V_0^2/2$. After $C_1$ is charged instantaneously to $V_S$, the final energy $E_1 = C_1 V_S^2/2$ and the voltage difference $\Delta V = V_S - V_0$, its square defines the energy loss:

$$E_1 - E_0 = \Delta E = \frac{C_1 (\Delta V)^2}{2} \qquad (2.2.2)$$

The same energy is dissipated as heat when the capacitor $C_1$ is charging through the resistor $R$ as shown in Fig. 2.1.2(a). According to the Joule-Lenz law, the power is $P = I^2 R$ and its integrated value is the heating loss:

$$E_h = I_0^2 R \cdot \int_0^\infty \left(e^{-t/\tau}\right)^2 dt = \frac{(\Delta V)^2}{2R} \tau = \frac{C_1 (\Delta V)^2}{2} \qquad (2.2.3)$$

In the particular case when $V_0 = 0$ the final energy $E_1$ equals the dissipated (radiated) energy, therefore half of the energy delivered by the source is lost. This fact corresponds to the law of energy conservation and can be proved by taking the integral of the delivered power:

$$E_d = V_S \int_0^\infty C_1 \frac{dv_1}{dt} dt = C_1 V_S \int_0^{V_S} dv_1 = CV_S^2 \qquad (2.2.4)$$

The above considerations can be applied to the ideal discharging circuit in Fig. 2.1.1(b), where the energy $E = E_{C1} + E_{C2}$ because it is stored in both capacitors $C_1$ and $C_2$. The initial voltages across the capacitors are $v_1(0_-) = V_S$, and $v_2(0_-) = 0$, by substitution into (2.2.1) the initial energy $E_0 = C_1 V_S^2 / 2$. After the circuit has closed, the final energy is given by the voltages $v_1(0_+) = v_2(0_+) = \dfrac{C_1}{C_1 + C_2} V_S$ derived in (2.1.3), so that $E_1 = \dfrac{(C_1 V_S)^2}{2(C_1 + C_2)}$, while the energy loss:

$$E_0 - E_1 = \Delta E = \frac{C_1 C_2}{C_1 + C_2} \frac{V_S^2}{2} \tag{2.2.5}$$

As in the previous case this energy loss should be compared with the heating loss when the energy is dissipated by the resistor during current flow. The corresponding switched circuit is shown in Fig. 2.1.2 (b). Substituting (2.1.16) into the power integral we obtain:

$$E_h = I_0^2 R \cdot \int_0^\infty \left(e^{-t/\tau}\right)^2 dt = \frac{V_S^2}{2R} \tau = \frac{C_1 C_2}{C_1 + C_2} \frac{V_S^2}{2} \tag{2.2.6}$$

So, the dissipated energy is equal to the energy loss found in (2.2.4) and caused by radiation, in the case of $C_1 = C_2 = C$ the loss will be $\Delta E = E_h = \dfrac{C V_S^2}{4} = \dfrac{1}{2} E_0$. The more interesting situation is when $C_1$ and $C_2$ are pre-charged to different $V_1$ and $V_2$ respectively.

The initial energy in this case $E_0 = \dfrac{C_1 V_1^2}{2} + \dfrac{C_2 V_2^2}{2}$. To know the final energy we have to find new values of the final voltages using the charge conservation law:

$$v_1(0_+) = v_2(0_+) = \frac{C_1}{C_1 + C_2} V_1 + \frac{C_2}{C_1 + C_2} V_2 \tag{2.2.7}$$

Substitution of these values into (2.2.1) yields $E_1 = \dfrac{(C_1 V_1 + C_2 V_2)^2}{2(C_1 + C_2)}$. The energy loss will be again proportional to the squared voltage difference:

$$E_0 - E_1 = \Delta E = \frac{C_1 C_2}{C_1 + C_2} \frac{(\Delta V)^2}{2} \tag{2.2.8}$$

## 2.3 Target voltages and SCC equivalent circuit

As mentioned above, SCC feature a set of discrete target voltages that can be contrasted with the continuous transfer function of inductor-based converters. This set of target voltages is closely related to the SCC efficiency over the full range of input voltages [26].

The target voltage is the no-load output voltage and is equal to some multiple $n$ of the input voltage. In general, $n$ is a function of the number of flying capacitors and the way that they are connected to the input and output and among themselves. Such interconnections are called hereinafter "SCC topologies". So, the target voltage is independent of the values of the flying capacitors and determined only by SCC topology, while $n$ can be a positive or negative rational number [11], [14], [18]. At each target voltage, the SCC efficiency reaches a maximum value and drops when the desired output voltage lies between or outside the target voltages.

For example, the commercial SCC [12] operates at the fixed output voltage $V_o = 1.8V$ and has two peaks of efficiency shown in Fig. 2.3.1. This SCC can be switched between two conversion ratios $n = 1/2$ and $2/3$.

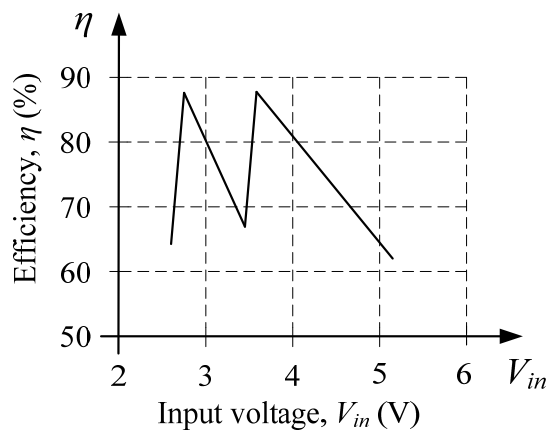

Figure 2.3.1: The output characteristics of a commercial SCC.

When the input voltage is lower than about 3.5V, the conversion ratio is set to $n = 2/3$ and for input voltage above 3.5V it is switched to $n = 1/2$. Consequently, high efficiency is observed when the input voltage is about 2.7V (1.8/(2/3)) and at 3.6V (1.8/(1/2)). When the SCC operates between and outside these two target voltages, the efficiency drops as the difference between the output voltage and 1.8V increases.

Any SCC can be modeled by an equivalent circuit that includes a voltage source $V_{TRG}$ and an internal resistance $R_{eq}$ as depicted schematically in Fig. 2.3.2 [15-26].

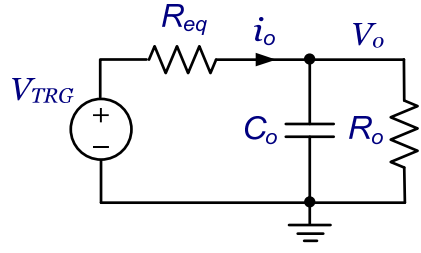

Figure 2.3.2: The SCC equivalent circuit.

In the model presentation of Fig. 2.3.2, the power losses are conveniently described as a function of the load current which simplifies the formulation of the input to output voltage ratio as well as the efficiency:

$$\frac{V_o}{V_{TRG}} = \frac{R_o}{R_o + R_{eq}} \qquad (2.3.1)$$

$$V_{TRG} = nV_{in} \qquad (2.3.2)$$

$$\eta = \frac{V_o}{V_{TRG}} = \frac{1}{n} \cdot \left( \frac{R_o}{R_o + R_{eq}} \right) \cdot \frac{V_o}{V_{in}} \qquad (2.3.3)$$

It is clear, that the highest efficiency will be achieved if $n$ is manipulated such that $V_{TRG}$ is made only slightly higher than the desired $V_o$, leaving a small voltage drop on $R_{eq}$. It is further clear that the best results can be obtained if the resolution by which $n$ is altered is high and when its values are evenly spaced. Previous attempts to improve the efficiency by changing $n$ on-the-fly gave SCC configurations with a limited number of target voltages, namely with a coarse resolution of $n$. As a result, the efficiency drops significantly when the required $n$ is in between the sparsely spread values of $n$.

SCC can be operated in open loop or closed loop configurations. In the open loop case, $n$ and $R_{eq}$ are fixed. In this case, the output voltage will not be regulated and will depend on $V_{in}$ and the load resistance $R_o$. In this situation, it is advantageous to reduce $R_{eq}$ as much as possible to keep the efficiency high. Regulation can be achieved by either changing n or $R_{eq}$ (or both) [11].

The no-load voltage can be changed by changing on-the-fly the SCC topology and hence altering $n$, while $R_{eq}$ can be changed by adding resistance to the circuit e.g. by placing a linearly controlled MOSFET in the charging/discharging paths. Other possibilities to vary $R_{eq}$ are frequency change, frequency dithering and duty cycle control [13], [14].

## 2.4 Demystifying the Equivalent Resistor Issue

> *"To improve is to change, to be perfect is to change often."*
> Winston Churchill

In this section we derive the equivalent resistor expression for a simple case of voltage follower SCC depicted in Fig. 2.4.1, where $R_1$ and $R_2$ represent the "on" resistances of $S_1$ and $S_2$ respectively, while *ESR* is the series loss component of the flying capacitor $C$. The analysis is based on the generic and unified average model [26], [54] and made under the assumption that the output capacitor $C_o$ is sufficiently large, so that the output voltage ripple is neglected.

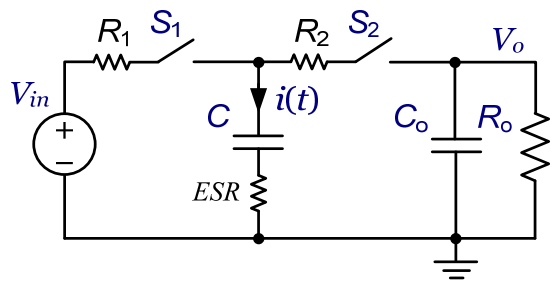

Figure 2.4.1: Voltage follower SCC.

Two clocks $\varphi_1$ and $\varphi_2$ shown in Fig. 2.4.2 alternately turn on/off the corresponding switches $S_1$ and $S_2$. The clocks are non-overlapping due to a dead time $p$, so that the total "on" duration $T_{on} = t_1 + t_2$ is smaller than the switching period $T_s$. During the interval $t_1$ the capacitor $C$ is charged by $V_{in}$ through $S_1$ and discharged to $V_o$ during the interval $t_2$ through $S_2$.

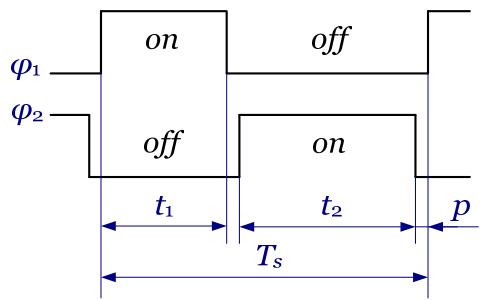

Figure 2.4.2: Two non-overlapping clocks $\varphi_1$ and $\varphi_2$.

Let $V_1$ and $V_2$ be the initial voltages across the capacitor $C$ at the instants just prior to its connection to the voltages $V_{in}$ and $V_o$ respectively. Since the initial voltages can be replaced by the voltage sources, the capacitor $C$ is charged by $\Delta V_1 = V_{in} - V_1$ during $t_1$ and discharged to $\Delta V_2 = V_2 - V_o$ during $t_2$.

It is convenient to consider a generic charge/discharge circuit presented in Fig. 2.4.3, where $V_C$ is the initial voltage across the capacitor $C$ and $R$ is the total loop resistance ("on" resistance of the switch $S$ plus the capacitor $ESR$).

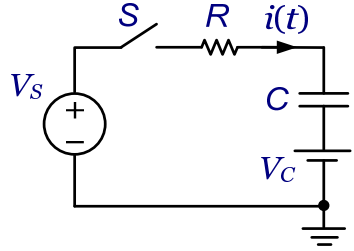

Figure 2.4.3: Generic charge/discharge circuit.

The switch $S$ remains turned on during $t_S$, so that both the energy dissipated by $R$ and the transferred charge can be found using $\Delta V = V_S - V_C = I_0 R$ and $i(t) = I_0 \cdot e^{-t/RC}$.

$$E_R = I_0^2 R \cdot \int_0^{t_S} \left(e^{-t/RC}\right)^2 dt = \frac{(\Delta V)^2 C}{2} \cdot (1 - e^{-2 t_S / RC}) \quad (2.4.1)$$

$$Q = \frac{\Delta V}{R} \cdot \int_0^{t_S} e^{-t/RC} dt = \Delta V \cdot C \cdot (1 - e^{-t_S/RC}) \quad (2.4.2)$$

Designating $\beta_1 = \dfrac{t_1}{(R_1 + ESR) \cdot C}$ and $\beta_2 = \dfrac{t_2}{(R_2 + ESR) \cdot C}$ we can relate the above results to the voltage follower SCC in Fig. 2.4.1. The energy losses for each interval $t_1$ and $t_2$ are:

$$E_1 = \frac{(\Delta V_1)^2 C}{2} \cdot (1 - e^{-2\beta_1}) \qquad E_2 = \frac{(\Delta V_2)^2 C}{2} \cdot (1 - e^{-2\beta_2}) \quad (2.4.3)$$

In the steady state, the charge transferred during $t_1$ and $t_2$ is the same:

$$Q = \Delta V_1 \cdot C \cdot (1 - e^{-\beta_1}) = \Delta V_2 \cdot C \cdot (1 - e^{-\beta_2}) \quad (2.4.4)$$

Since the average current $I_{av} = Q/T_s = f_s Q$, we can write

$$I_{av} = f_s \Delta V_1 C \cdot (1 - e^{-\beta_1}) = f_s \Delta V_2 C \cdot (1 - e^{-\beta_2}) \quad (2.4.5)$$

Rearranging the terms of (2.4.5) yields:

$$\Delta V_1 = \frac{I_{av}}{f_s C \cdot (1 - e^{-\beta_1})} \qquad \Delta V_2 = \frac{I_{av}}{f_s C \cdot (1 - e^{-\beta_2})} \quad (2.4.6)$$

These voltage differences are substituted into (2.4.3), so that:

$$E_1 = \frac{I_{av}^2(1-e^{-2\beta_1})}{2f_s^2 C \cdot (1-e^{-\beta_1})^2} \qquad E_2 = \frac{I_{av}^2(1-e^{-2\beta_2})}{2f_s^2 C \cdot (1-e^{-\beta_2})^2} \qquad (2.4.7)$$

Because the total energy loss $E_R = E_1 + E_2$,

$$E_R = \frac{I_{av}^2}{2f_s^2 C} \cdot \left[\frac{1-e^{-2\beta_1}}{(1-e^{-\beta_1})^2} + \frac{1-e^{-2\beta_2}}{(1-e^{-\beta_2})^2}\right] \qquad (2.4.8)$$

Or after simplification:

$$E_R = \frac{I_{av}^2}{2f_s^2 C} \cdot \left(\frac{1+e^{-\beta_1}}{1-e^{-\beta_1}} + \frac{1+e^{-\beta_2}}{1-e^{-\beta_2}}\right) \qquad (2.4.9)$$

The total average power loss $P_T = E_R/T_s = E_R f_s$, so that:

$$P_T = \frac{I_{av}^2}{2f_s C} \cdot \left(\frac{1+e^{-\beta_1}}{1-e^{-\beta_1}} + \frac{1+e^{-\beta_2}}{1-e^{-\beta_2}}\right) \qquad (2.4.10)$$

Comparing (2.4.10) with $P_T = I_{av}^2 \cdot R_{eq}$ we conclude that the equivalent resistor is:

$$R_{eq} = \frac{1}{2f_s C} \cdot \left(\frac{1+e^{-\beta_1}}{1-e^{-\beta_1}} + \frac{1+e^{-\beta_2}}{1-e^{-\beta_2}}\right) \qquad (2.4.11)$$

Employing the definition of $\coth\left(\frac{x}{2}\right) = \frac{1+e^{-x}}{1-e^{-x}}$, rewrite (2.4.11) as:

$$R_{eq} = \frac{1}{2f_s C} \cdot \left[\coth\left(\frac{\beta_1}{2}\right) + \coth\left(\frac{\beta_2}{2}\right)\right] \qquad (2.4.12)$$

For the particular case of $\beta_1 = \beta_2 = \beta$, the general expression (2.4.12) is reduced to:

$$R_{eq} = \frac{1}{f_s C} \cdot \coth\left(\frac{\beta}{2}\right) \qquad (2.4.13)$$

Assuming zero dead time and $\beta = T_s/2RC$, we can rewrite (2.4.13) as:

$$R_{eq} = 2R \cdot \beta \cdot \coth\left(\frac{\beta}{2}\right) \qquad (2.4.14)$$

Consider an extreme case of (2.4.14) when $\beta \to 0$:

$$\lim_{\beta \to 0} R_{eq} = 2R \cdot \lim_{\beta \to 0}\left[\beta \cdot \coth\left(\frac{\beta}{2}\right)\right] = 4R \qquad (2.4.15)$$

This seemingly surprising result has a simple explanation. In the circuit of Fig. 2.4.1 the momentary current during each switching phase is $2I_o$ (to make the average current $I_{av} = I_o$), so that the losses are $(2I_o)^2 R = I_o^2 4R$.

An additional extreme case for (2.4.13) is $\beta \to \infty$ that results in $R_{eq}$ reduced to the well known expression:

$$\lim_{\beta \to \infty} R_{eq} = \frac{1}{f_s C} \cdot \lim_{\beta \to \infty}\left[\coth\left(\frac{\beta}{2}\right)\right] = \frac{1}{f_s C} \qquad (2.4.16)$$

To demonstrate how both the above limits (2.4.15) and (2.4.16) are reached we built the graphs of the corresponding terms $\beta \cdot \coth(\beta/2)$ and $\coth(\beta/2)$ as depicted in Fig. 2.4.4.

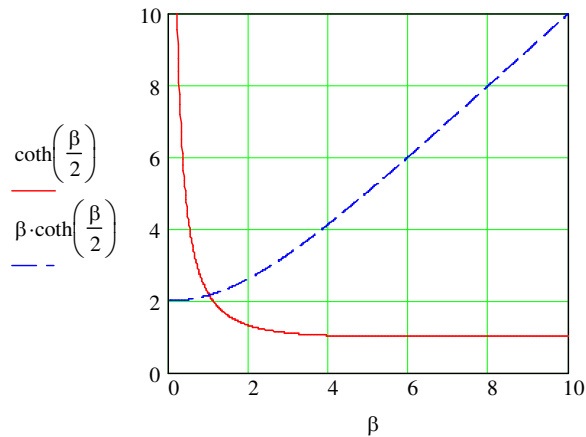

Figure 2.4.4: Functions $\coth(\beta/2)$ and $\beta \cdot \coth(\beta/2)$.

It is evident that for $\beta \approx 5$, the term $\coth(\beta/2) \approx 1$. This fact can be simply explained since the time constant $\tau = RC$, $\beta = t/\tau$ and the transient is quite finished after $t = 5\tau$. Thus, (2.4.16) corresponds to the case of full charging/discharging of the flying capacitors. On the other hand, when $\beta \to 0$, the term $\coth(\beta/2) \to \infty$, while $\beta \cdot \coth(\beta/2) \to 2$. Since (2.4.14) is written under the assumption of $\beta = T_s/2RC$, where $T_s \neq 0$, for $\beta \to 0$ we need $\tau \to \infty$.

In practice, R is relatively small and one can get $\tau \to \infty$ with sufficiently large flying capacitors. So, (2.4.15) corresponds to the case of partial charging/discharging while the current through R is constant.

# 3. PROPOSED CLASS OF SCC WITH BINARY RESOLUTION

## 3.1 Extended Binary (EXB) Representation

*"It is through science that we prove, but through intuition that we discover."*
Jules H. Poincare

As mentioned above, the total SCC efficiency over the full range of input voltages can be improved by increasing the number of target voltages. In order to design a step-down SCC with closely spaced multiple target voltages, we have developed an Extended Binary (EXB) representation. According to this approach, the flying capacitors are automatically kept charged to binary weighted voltages and, consequently, the resolution of the target voltages is binary. The resolution can be made higher by increasing the number of flying capacitors.

For the resolution $n$, consider a set of fractions $M_n$ in the range (0, 1) with odd numerators 1, 3, …, $2^n - 1$ and denominator $2^n$. Any fraction $M_n$ can be represented in the form:

$$M_n = A_0 + \sum_{j=1}^{n} +A_j 2^{-j} \qquad (3.1.1)$$

where $A_0$ can be either 0 or 1, and $A_j$ can take any of three values -1, 0, 1.

The expression (3.1.1) defines the Extended Binary (EXB) representation, which differs from its conventional binary counterpart since $A_j$ can be -1. Because of the three values -1, 0, 1 for $A_j$, the EXB representation is akin to binary signed-digit (BSD) representation of integer numbers, for example:

$$\begin{aligned}
5 &= 0 + 4 + 2 - 1 \rightarrow \{0\ 1\ 1\ \text{-1}\} \\
5 &= 8 - 4 + 0 + 1 \rightarrow \{1\ \text{-1}\ 0\ 1\} \\
5 &= 8 + 0 - 2 - 1 \rightarrow \{1\ 0\ \text{-1}\ \text{-1}\}
\end{aligned} \qquad (3.1.2)$$

As seen from (3.1.2), the BSD representation for a given integer is not unique and this property is used mostly for carry-save fast computer arithmetic. We have modified the BSD representation for fractions $M_n$ limited in the range (0, 1). As a result, the coefficient $A_0$ in the EXB representation (3.1.1) is not allowed to be "1".

Because of the redundancy that comes from the BSD representation, any fraction $M_n$ can be represented by a number of EXB codes, for example:

$$5/8 = 0 + 2^{-1} + 2^{-2} - 2^{-3} \rightarrow \{0\ 1\ 1\ \text{-}1\}$$
$$5/8 = 1 - 2^{-1} + 0 + 2^{-3} \rightarrow \{1\ \text{-}1\ 0\ 1\} \qquad (3.1.3)$$
$$5/8 = 1 + 0 - 2^{-2} - 2^{-3} \rightarrow \{1\ 0\ \text{-}1\ \text{-}1\}$$

In the next section we provide a simple procedure to spawn all the EXB codes for a given fraction $M_n$. This procedure will be followed by a number of corollaries, which are crucial to define and explicate the operation of the EXB based SCC.

## 3.2 Spawning the EXB codes and its corollaries

> *"Get your facts first then you can distort them as you please."*
> Mark Twain

In order to generate all the EXB codes corresponding to a given fraction $M_n$ within the range (0, 1), we use a procedure that involves adding and subtracting the coefficient $A_j = 1$ to the conventional binary code of $M_n$.

**Spawning the EXB codes.** This procedure is iterative and starts from any $A_j = 1$ in the conventional binary code of $M_n$. Adding "1" to this $A_j$ results in "0" and "1" from the left as the carry. To maintain the value of $M_n$ we subtract "1" from the obtained $A_j$, and spawn thereby a new EXB code. The procedure repeats for all $A_j = 1$ in the original code and for all $A_j = 1$ in each spawned EXB code.

In example (3.2.1), four alternative EXB codes are spawned from the conventional binary code of $M_3 = 3/8$. The EXB codes for other fractions $M_n$ with the resolution $n = 1\ldots3$ are summarized in Table 3.2.1.

$$
\begin{array}{cccc}
\quad\ \ \downarrow & \quad\ \ \downarrow & \ \downarrow\quad\ \  & \downarrow\quad\quad\  \\
2^0\ 2^{-1}\ 2^{-2}\ 2^{-3} & 2^0\ 2^{-1}\ 2^{-2}\ 2^{-3} & 2^0\ 2^{-1}\ 2^{-2}\ 2^{-3} & 2^0\ 2^{-1}\ 2^{-2}\ 2^{-3} \\
\ \ 0\ \ 0\ \ 1\ \ 1 & \ \ 0\ \ 0\ \ 1\ \ 1 & \ \ 0\ \ 1\ \ 0\ \text{-}1 & \ \ 0\ \ 1\ \text{-}1\ \ 1 \\
+\ 0\ \ 0\ \ 0\ \ \mathbf{1} & +\ 0\ \ 0\ \ \mathbf{1}\ \ 0 & +\ 0\ \ \mathbf{1}\ \ 0\ \ 0 & +\ 0\ \ \mathbf{1}\ \ 0\ \ 0 \\ \hline
\ \ 0\ \ 1\ \ 0\ \ 0 & \ \ 0\ \ 1\ \ 0\ \ 1 & \ \ 1\ \ 0\ \ 0\ \text{-}1 & \ \ 1\ \ 0\ \text{-}1\ \ 1 \\
+\ 0\ \ 0\ \ 0\ \text{-}\mathbf{1} & +\ 0\ \ 0\ \text{-}\mathbf{1}\ \ 0 & +\ 0\ \text{-}\mathbf{1}\ \ 0\ \ 0 & +\ 0\ \text{-}\mathbf{1}\ \ 0\ \ 0 \\ \hline
\ \ 0\ \ 1\ \ 0\ \text{-}1 & \ \ 0\ \ 1\ \text{-}1\ \ 1 & \ \ 1\ \text{-}1\ \ 0\ \text{-}1 & \ \ 1\ \text{-}1\ \text{-}1\ \ 1
\end{array}
\qquad (3.2.1)
$$

**Corollary 1:** For the resolution $n$, the minimum number of EXB codes is $n + 1$.

This is because each of the "1"s in the conventional binary code with resolution $n$ generates a new EXB code and a carry. Further iterations cause the carry to propagate, so that each "0" in the conventional binary code is turned to "1", which is also operated on to spawn a new code. So, the minimum number of codes is the original code plus $n$ that is, $n + 1$.

**Corollary 2:** Each $A_j = 1$ in either the conventional binary or spawned EXB code yields at least one $A_j = -1$ in the same position $j$ of another EXB code.

This is because the spawning procedure involves subtracting "1" from $A_j = 0$.

Both the above corollaries are very important and, as detailed in the following, provide the self-adjusting target voltage $M_n \cdot V_{in}$ at the output of the EXB based SCC, irrespectively of the values of the used capacitors.

Table 3.2.1: The EXB codes of $M_n$, $n = 1…3$.

| M₃ = 1/8 | | | | M₂ = 2/8 | | | | M₃ = 3/8 | | | | M₁ = 4/8 | | | |
|---|---|---|---|---|---|---|---|---|---|---|---|---|---|---|---|
| $A_0$ | $A_1$ | $A_2$ | $A_3$ | $A_0$ | $A_1$ | $A_2$ | $A_3$ | $A_0$ | $A_1$ | $A_2$ | $A_3$ | $A_0$ | $A_1$ | $A_2$ | $A_3$ |
| 1 | -1 | -1 | -1 | 1 | -1 | -1 | 0 | 1 | -1 | 0 | -1 | 1 | -1 | 0 | 0 |
| 0 | 1 | -1 | -1 | 0 | 1 | -1 | 0 | 0 | 1 | 0 | -1 | 0 | 1 | 0 | 0 |
| 0 | 0 | 1 | -1 | 0 | 0 | 1 | 0 | 1 | -1 | -1 | 1 | | | | |
| 0 | 0 | 0 | 1 | | | | | 0 | 1 | -1 | 1 | | | | |
| | | | | | | | | 0 | 0 | 1 | 1 | | | | |

Table 3.2.1: cont'd.

| M₃ = 5/8 | | | | M₂ = 6/8 | | | | M₃ = 7/8 | | | |
|---|---|---|---|---|---|---|---|---|---|---|---|
| $A_0$ | $A_1$ | $A_2$ | $A_3$ | $A_0$ | $A_1$ | $A_2$ | $A_3$ | $A_0$ | $A_1$ | $A_2$ | $A_3$ |
| 1 | 0 | -1 | -1 | 1 | -1 | 1 | 0 | 1 | 0 | 0 | -1 |
| 1 | -1 | 1 | -1 | 1 | 0 | -1 | 0 | 1 | 0 | -1 | 1 |
| 0 | 1 | 1 | -1 | 0 | 1 | 1 | 0 | 1 | -1 | 1 | 1 |
| 1 | -1 | 0 | 1 | | | | | 0 | 1 | 1 | 1 |
| 0 | 1 | 0 | 1 | | | | | | | | |

## 3.3 Combinatorial method to obtain EXB codes

*"The true delight is in the finding out rather than in the knowing."*
Isaac Asimov

Due to the spawning procedure described in the previous section, we have derived important properties of the EXB codes. However, from the viewpoint of performance, this procedure is slow because each EXB code of $M_n$ is obtained by the series, digit-by-digit adding and subtracting the coefficient $A_j$ = 1. The alternative combinatorial method proposed in this section is parallel, and therefore faster than the previous one.

According to the definition, the EXB representation of $M_n$ contains $n$ coefficients $A_j$, which can take any of three values: -1, 0, 1. We consider all combinations of these values arranged at $n$ positions as a matrix **M** of $3^n$ rows by $n$ columns. This matrix is obtained by the full factorial design, where each level is 3, and the number of levels is $n$. Note that **M** defines the representations of the numbers $\frac{1-3^n}{2},\ldots,\frac{3^n-1}{2}$ in the balanced ternary number system.

For the sake of an exact integer calculation, we multiply both the sides of the EXB formula (3.1.1) by $2^n$. As a result each EXB weight $2^{-j}$ is replaced by $2^{n-j}$ and we have $n$ powers of two, which compose a column-vector **K** of length $n$. Multiplying the matrix **M** by this vector yields a column-vector **F** of length $3^n$. To indicate $A_0$ = 0 and $A_0$ = 1 in the EXB codes we introduce a column-vector **B** of the same length. The positive elements of the vector **F** correspond to $A_0$ = 0 and transferred to the vector **B** as zeros, while the negative elements correspond $A_0$ = 1 and transferred as ones.

We complete the negative elements of the vector **F** to the positive by adding $2^n$ and obtain a column-vector **F'**. So, this vector will contain $3^n$ elements, which are the numerators $m = 1,\ldots, 2^n - 1$ of all $M_n$. The search for certain $m$ results in several row indexes, while the same rows of **B** and **M** compose the EXB codes for given $M_n$.

For the resolution $n = 3$, this combinatorial method is demonstrated step-by-step in the following, while the full factorial design matrix and the used vectors are shown in Fig. 3.3.1.

1) All combinations of -1, 0, 1 in 3 positions are given by the full factorial design matrix **M** of $3^3$ rows by 3 columns obtained with the MATLAB command fullfact([3 3 3])–2.

2) The column-vector **K** = [4; 2; 1], so that the product of **M** and **K** is the column-vector **F** of length $3^3$ comprising the numbers from $1 - 2^3$ through $2^3 - 1$.

3) The positive numbers of the vector **F** are transferred to the vector **B** as zeros, while the negative numbers are transferred as ones.

4) At the same time the negative numbers in **F** are completed to the positive by adding $2^3$, so that the completed vector **F'** contains the numbers from 0 through $2^3 - 1$.

5) Since **F'** contains $3^3$ elements, any of $m = 1, \ldots, 2^3 - 1$ appears in **F'** more than once, and search for certain $m$ results in several row indexes. The same rows of **B** and **M** compose the EXB codes for given $M_n$ ($n = 1\ldots3$). Such a gathering is demonstrated in Fig. 3.3.1, where $m = 3$ that corresponds to $M_3 = 3/8$.

Comparing the obtained codes with the codes presented in Table 3.2.1 we conclude that the combinatorial method yields the same result as the procedure spawning the EXB codes.

$$
\mathbf{B} = \begin{bmatrix} 1 \\ 1 \\ 1 \\ 1 \\ 1 \\ 1 \\ 1 \\ 1 \\ 1 \\ 1 \\ 1 \\ 1 \\ 1 \\ 0 \\ 0 \\ 0 \\ 0 \\ 0 \\ 0 \\ 0 \\ 0 \\ 0 \\ 0 \\ 0 \\ 0 \\ 0 \\ 0 \end{bmatrix} \quad \mathbf{M} = \begin{bmatrix} -1 & -1 & -1 \\ -1 & -1 & 0 \\ -1 & -1 & 1 \\ -1 & 0 & -1 \\ -1 & 0 & 0 \\ -1 & 0 & 1 \\ -1 & 1 & -1 \\ -1 & 1 & 0 \\ -1 & 1 & 1 \\ 0 & -1 & -1 \\ 0 & -1 & 0 \\ 0 & -1 & 1 \\ 0 & 0 & -1 \\ 0 & 0 & 0 \\ 0 & 0 & 1 \\ 0 & 1 & -1 \\ 0 & 1 & 0 \\ 0 & 1 & 1 \\ 1 & -1 & -1 \\ 1 & -1 & 0 \\ 1 & -1 & 1 \\ 1 & 0 & -1 \\ 1 & 0 & 0 \\ 1 & 0 & 1 \\ 1 & 1 & -1 \\ 1 & 1 & 0 \\ 1 & 1 & 1 \end{bmatrix} \quad \mathbf{K} = \begin{bmatrix} 4 \\ 2 \\ 1 \end{bmatrix} \quad \mathbf{F} = \begin{bmatrix} -7 \\ -6 \\ -5 \\ -5 \\ -4 \\ -3 \\ -3 \\ -2 \\ -1 \\ -3 \\ -2 \\ -1 \\ -1 \\ 0 \\ 1 \\ 1 \\ 2 \\ 3 \\ 1 \\ 2 \\ 3 \\ 3 \\ 4 \\ 5 \\ 5 \\ 6 \\ 7 \end{bmatrix} \quad \mathbf{F'} = \begin{bmatrix} 1 \\ 2 \\ 3 \\ 3 \\ 4 \\ 5 \\ 5 \\ 6 \\ 7 \\ 5 \\ 6 \\ 7 \\ 7 \\ 0 \\ 1 \\ 1 \\ 2 \\ 3 \\ 1 \\ 2 \\ 3 \\ 3 \\ 4 \\ 5 \\ 5 \\ 6 \\ 7 \end{bmatrix}
$$

Figure 3.3.1: Full factorial design EXB matrix of size 27×3 and used vectors.

## 3.4 Translating the EXB codes to SCC topologies

For certain EXB fraction $M_n$, we consider a step-down SCC system that includes a voltage source $V_{in}$, a set of *n* flying capacitors $C_j$ and an output capacitor $C_o$ connected in parallel with the load $R_o$. These components are connected in accordance with the EXB codes of $M_n$ in such a way that $C_o$ is continuously charged. In particular, the EXB coefficient $A_0$ is responsible for the connection of $V_{in}$, while the connection of each flying capacitor $C_j$ is determined by the EXB coefficient $A_j$. Irrespective of the connection of $V_{in}$ the flying capacitors $C_j$ are always connected serially. To configure the EXB based SCC topologies we use the following rules:

1) If $A_0 = 1$, then $V_{in}$ is connected.
2) If $A_0 = 0$, then $V_{in}$ is not connected.
3) If $A_j = -1$, then $C_j$ is charged.
4) If $A_j = 0$, then $C_j$ is not connected.
5) If $A_j = 1$, then $C_j$ is discharged.

As an example we translate all the EXB codes of $M_3 = 3/8$ presented in Table 3.4.1 to the corresponding SCC topologies. Since the resolution $n = 3$ we need three flying capacitors $C_1$, $C_2$ and $C_3$, the serial connection of which is determined by $A_1$, $A_2$ and $A_3$ respectively. Thus, each EXB code of $M_3 = 3/8$ leads to a specific SCC topology as depicted in Figure 3.4.1.

Table 3.4.1.

| $M_3 = 3/8$ | | | |
|---|---|---|---|
| $A_0$ | $A_1$ | $A_2$ | $A_3$ |
| 1 | -1 | -1 | 1 |
| 0 | 1 | -1 | 1 |
| 1 | -1 | 0 | -1 |
| 0 | 1 | 0 | -1 |
| 0 | 0 | 1 | 1 |

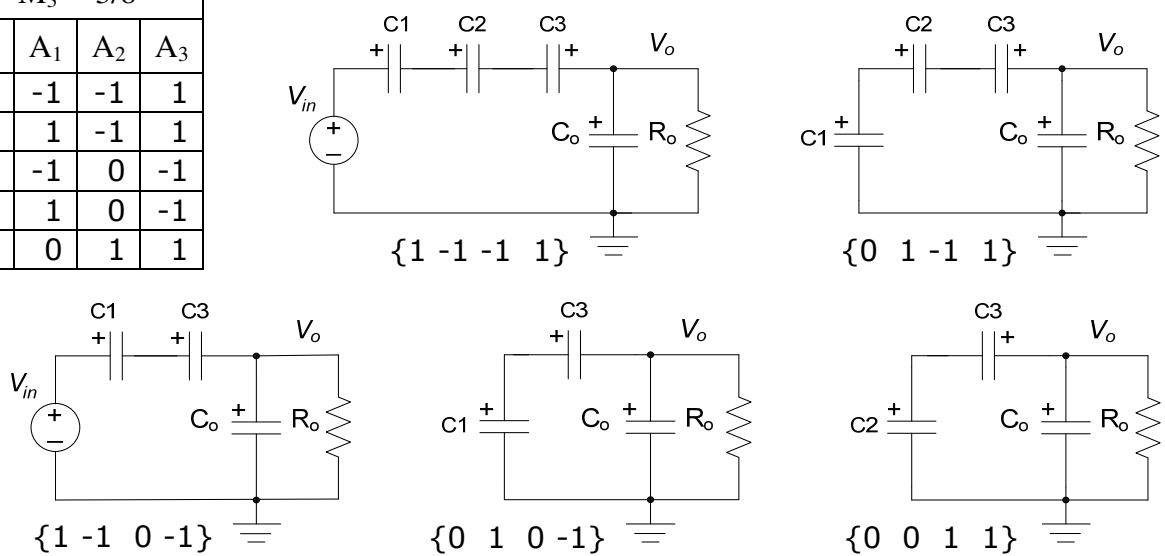

Figure 3.4.1: SCC topologies configured from the EXB codes of $M_3 = 3/8$.

We assume that in each SCC topology of Fig. 3.4.1, the flying capacitors $C_1$, $C_2$ and $C_3$ keep the voltages $V_1 = 2^{-1} \cdot V_{in}$, $V_2 = 2^{-2} \cdot V_{in}$ and $V_3 = 2^{-3} \cdot V_{in}$ respectively. Multiplying $V_{in}$ and these voltages by the corresponding coefficients $A_0$, $A_1$, $A_2$ and $A_3$ in the EXB codes of $M_3 = 3/8$, we find their algebraic sum, which is equal to the target voltage $V_o = 3/8 \cdot V_{in}$.

Generally, translating all the EXB codes of certain $M_n$ to the SCC topologies, we ought to obtain the target voltage $V_o = M_n \cdot V_{in}$, under the condition that each flying capacitor $C_j$ keeps the voltage $V_j = 2^{-j} \cdot V_{in}$. In the following we show that all the voltages in the EXB based SCC are self-adjusting to the above specified values and this property is due to Corollaries 1 and 2 of the procedure for spawning the EXB codes.

## 3.5 Self-adjusting voltages in the EXB based SCC

*"Make everything as simple as possible, but not simpler."*
Albert Einstein

In this section we consider the EXB based SCC under the assumption that in each SCC topology all the capacitors voltages remain constant but of unknown values. Applying Kirchhoff's Voltage Law (KVL) to $w$ different SCC topologies we compose a system of $w$ linear equations. If this system has a unique solution, we obtain the target and binary weighted voltages across the output and flying capacitors respectively.

The KVL states that the algebraic sum of all voltages around any closed path in a circuit is zero. Any SCC topology is a closed path circuit because the flying capacitors are charged and discharged thus providing charge transfer. The output voltage of the SCC is assumed to be constant and the KVL is applied to the voltages across the capacitors engaged in a SCC topology.

First we consider the simplest voltage halving SCC defined by $M_1 = 1/2$. The topologies of this SCC are depicted in Fig. 3.5.1, where $C_1$ and $C_o$ keep the voltages $V_1$ and $V_o$ respectively.

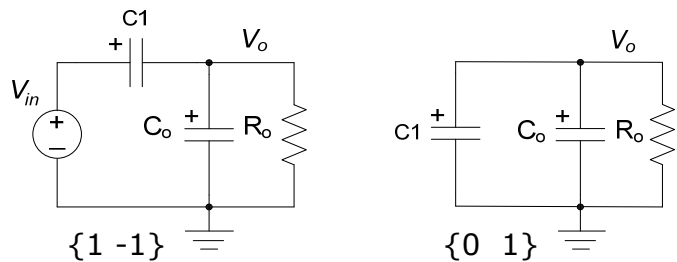

Figure 3.5.1: Topologies of the voltage halving SCC.

The system of linear equations for both topologies of Fig. 3.5.1 is:

$$\begin{cases} V_{in} - V_1 = V_o \\ V_1 = V_o \end{cases} \quad (3.5.1)$$

The solution of (3.5.1) is trivial: $\quad V_o = V_1 = \dfrac{1}{2} V_{in} \quad (3.5.2)$

Generally, a system of equations for the EXB based SCC may be composed directly from the corresponding EXB codes. As an example, we show that the EXB codes of $M_3 = 3/8$ lead not only to the SCC topologies of Fig. 3.5.2, but also to the system (3.5.3).

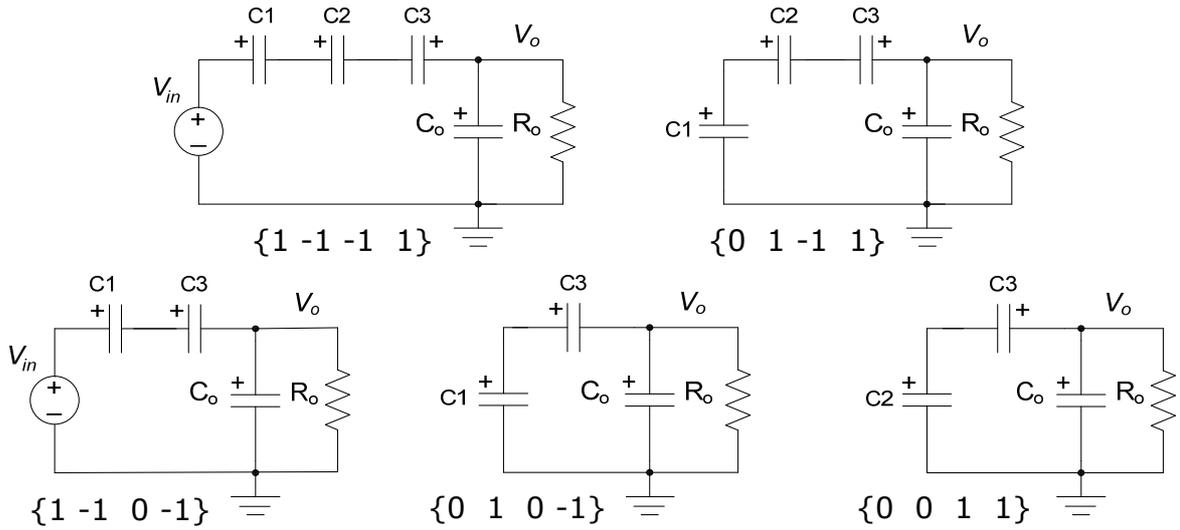

Figure 3.5.2: SCC topologies configured from the EXB codes of $M_3 = 3/8$.

$$\begin{cases} V_{in} - V_1 - V_2 + V_3 = V_o \\ 0 + V_1 - V_2 + V_3 = V_o \\ V_{in} - V_1 + 0 - V_3 = V_o \\ 0 + V_1 + 0 - V_3 = V_o \\ 0 + 0 + V_2 + V_3 = V_o \end{cases} \quad (3.5.3)$$

The number of equations in (3.5.3) is identical to the number $w$ of all the EXB codes of $M_3 = 3/8$ and equals to 5, while the number of unknowns is equal to 4 and defined as the resolution $n = 3$ plus one. Grouping the unknowns in (3.5.3) at the left hand side yields:

$$\begin{cases} -V_1 - V_2 + V_3 - V_o = -V_{in} \\ V_1 - V_2 + V_3 - V_o = 0 \\ -V_1 + 0 - V_3 - V_o = -V_{in} \\ V_1 + 0 - V_3 - V_o = 0 \\ 0 + V_2 + V_3 - V_o = 0 \end{cases} \quad (3.5.4)$$

The system of equations (3.5.4) contains two non-zero free terms as the negative value of $V_{in}$. Generally, the connection of $V_{in}$ is provided by Corollary 1 of the procedure spawning the EXB codes as follows. Consider the conventional binary code of $M_n$, where the coefficient $A_j$ takes either "1" or "0". Due to Corollary 1, the case $A_j = 0$ is turned to $A_j = 1$, which is used to generate the coefficient $A_0 = 1$ responsible for the connection of $V_{in}$.

So, the EXB based SCC is described by a system of linear equations which contains at least one non-zero free term. In linear algebra, such a system is called non-homogeneous.

Returning to (3.5.4) we normalize it to $V_{in}$:

$$\begin{cases} -1 \cdot x_1 - 1 \cdot x_2 + 1 \cdot x_3 - 1 \cdot x_4 = -1 \\ 1 \cdot x_1 - 1 \cdot x_2 + 1 \cdot x_3 - 1 \cdot x_4 = 0 \\ -1 \cdot x_1 + 0 \cdot x_2 - 1 \cdot x_3 - 1 \cdot x_4 = -1 \\ 1 \cdot x_1 + 0 \cdot x_2 - 1 \cdot x_3 - 1 \cdot x_4 = 0 \\ 0 \cdot x_1 + 1 \cdot x_2 + 1 \cdot x_3 - 1 \cdot x_4 = 0 \end{cases} \quad (3.5.5)$$

where

$$\begin{aligned} x_1 &= V_1/V_{in} \\ x_2 &= V_2/V_{in} \\ x_3 &= V_3/V_{in} \\ x_4 &= V_o/V_{in} \end{aligned} \quad (3.5.6)$$

The conventional brief notation for (3.5.5) is **Ax = b**, where

$$\mathbf{A} = \begin{bmatrix} -1 & -1 & 1 & -1 \\ 1 & -1 & 1 & -1 \\ -1 & 0 & -1 & -1 \\ 1 & 0 & -1 & -1 \\ 0 & 1 & 1 & -1 \end{bmatrix} \quad \mathbf{x} = \begin{bmatrix} x_1 \\ x_2 \\ x_3 \\ x_4 \end{bmatrix} \quad \mathbf{b} = \begin{bmatrix} -1 \\ 0 \\ -1 \\ 0 \\ 0 \end{bmatrix} \quad (3.5.7)$$

In order to investigate the solvability of (3.5.5) we supplement the coefficient matrix **A** with the vector **b** and form thereby the augmented matrix **A1**:

$$\mathbf{A1} = \begin{bmatrix} -1 & -1 & 1 & -1 & -1 \\ 1 & -1 & 1 & -1 & 0 \\ -1 & 0 & -1 & -1 & -1 \\ 1 & 0 & -1 & -1 & 0 \\ 0 & 1 & 1 & -1 & 0 \end{bmatrix} \quad (3.5.8)$$

According to the Kronecker-Capelli theorem [31], [36], a non-homogeneous system has at least one solution if and only if the rank of its coefficient matrix **A** is equal to the rank of its augmented matrix **A1**. This theorem has a corollary that specifies the number of solutions.

The solution is unique if and only if the rank the augmented matrix **A1** equals the number of unknowns. If the rank of **A** equals the rank of **A1**, but is less than the number of unknowns, the system has an infinite number of solutions. If, on the other hand, the rank of **A1** is greater than the rank of **A**, the system has no solutions.

Note that the rank of any matrix is equal to the row rank and equal to the column rank, and as a consequence, the maximum number of linearly independent rows of a matrix is equal to the maximum number of its linearly independent columns.

For the resolution $n$ the number of columns in the matrix **A** is $n + 1$, while the number of rows is provided by Corollary 1 to be $w \geq n + 1$. According to the above theorem we conclude that the rank of **A** as well as the rank of **A1** must be equal to $n + 1$, which is exactly the number of unknowns. So, we have to prove rigorously that the procedure spawning the EXB codes provides exactly $n + 1$ linearly independent rows (or columns) in both the matrices **A** and **A1**.

From a practical point of view the EXB based SCC with high resolution $n > 16$ will be very expensive for realization and therefore we suppose $1 \leq n \leq 16$. For each $n$ we calculated the rank of **A** and the rank of **A1** numerically in MATLAB, and ensured that the procedure spawning the EXB codes leads to a system from $n + 1$ linearly independent equations, while its unique solution is $M_n$. Giving a rigorous theoretical proof of this for any $n$ is planned for future work.

In the case of system (3.5.5), the rank of **A** is 4, which is equal to the rank of **A1** and equal to the number of unknowns. Thus, system (3.5.5) has a unique solution and comprises one redundant equation. Solving (3.5.5) in the form $\mathbf{x} = \mathbf{A}^{-1}\mathbf{b}$ we obtain:

$$\begin{aligned} x_1 &= 1/2 & V_1 &= 1/2 \cdot V_{in} \\ x_2 &= 1/4 & V_2 &= 1/4 \cdot V_{in} \\ x_3 &= 1/8 & \text{or} \quad V_3 &= 1/8 \cdot V_{in} \\ x_4 &= 3/8 & V_o &= 3/8 \cdot V_{in} \end{aligned} \quad (3.5.9)$$

Since solution (3.5.9) is unique, the voltages $V_j$ ($j = 1\ldots3$) across the flying capacitors $C_j$ are self-adjusting to the specified values of $2^{-j}\cdot V_{in}$, so that there is no need for any closed loop control scheme [48-50] to assure that these voltages are reached and hence the output will always self stabilize to the expected target voltage $V_o = 3/8 \cdot V_{in}$.

Generally, for a given EXB fraction $M_n$ with resolution $n$, there are two ways to provide the self-adjusting voltages $V_j = 2^{-j}\cdot V_{in}$ ($j = 1\ldots n$) and $V_o = M_n \cdot V_{in}$. One way is to configure only those SCC topologies, which correspond to $n + 1$ linearly independent equations. The other way is to introduce redundant SCC topologies in addition to the $n + 1$ mentioned above.

By configuring these $w \geq n + 1$ different SCC topologies periodically we provide continuous charge transfer and, consequently, current flow through the load. An intuitive explanation for this is that while the topologies change each flying capacitor goes through a sequence of charge and discharge. This is assured by Corollary 2 of the procedure spawning the EXB codes, which states that for each "1" in the EXB code there is at least one "-1" in the same position. This means that the theory predicts that in the EXB based SCC, all capacitors are going through the sequence of charge and discharge.

According to the number $w$ of different SCC topologies we introduce the same number $w$ of time intervals $t_k$, where $k = 1, \ldots, w$. During the interval $t_k$, the corresponding $k$-th topology does not change and repeats as depicted conceptually for $M_n = 3/8$ in Fig. 3.5.3 with a period $T = \sum_{k=1}^{w} t_k$ at the intervals $t_{k+p \cdot w}$, where $p = 0, 1, 2, \ldots$ is the number of period.

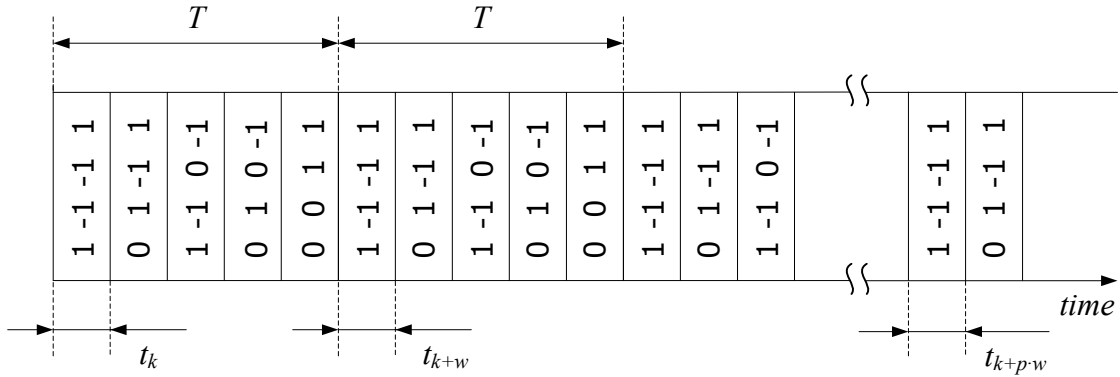

Figure 3.5.3: The perpetual EXB sequences of the SCC with $M_3 = 3/8$.

Irrespective of the order in which the topologies repeat, the voltages $V_j$ and $V_o$ eventually reach and stay at the specified values of $2^{-j} \cdot V_{in}$ and $M_n \cdot V_{in}$ even if the SCC starts with zero or arbitrary voltages across the capacitors or when it is subjected to a disturbance. Due to this property, the EXB based SCC can be considered to be hardware for solving a system of linear equations by an iterative method.

# 3.6 Method to reduce output voltage ripple

*"To every action there is always opposed an equal reaction."*
Isaac Newton

Output voltage ripple in the EXB based SCC can be reduced when the charging and discharging of each flying capacitor $C_j$ are interleaved. Since the connection polarity of $C_j$ is determined by the sign of coefficient $A_j$, we can form an alternating sequence of EXB codes, where the sign of the next $A_j$ is opposite to the previous one. Evidently, a significant ripple reduction can be achieved by a balanced sequence, where the alternating nonzero coefficients $A_j$ are spaced with a constant number of zeros. As shown in the following, the balanced sequence can be formed by replicating some of the EXB codes.

For the sake of an exact integer calculation, we multiply each EXB weight $2^{-j}$ by $2^n$ and compose a row-vector **Y** of length *n*. The elements of **Y** can represent the numbers from 1 through $2^n - 1$ in accordance with conventional binary codes. We consider these codes as a matrix **B** of two level full factorial design created in MATLAB with the command B = ff2n(n), where n is the number of columns. Multiplying **B** and **Y** element-by-element yields the weight matrix **U**. To introduce $A_j = \pm 1$ we use the same binary matrix **B** as exponent the for -1 and obtain thereby a sign matrix **P**.

Having two matrices **U** and **P** we find their product **M** of size $2^n \times 2^n$, each column of which contains the numbers from $1 - 2^n$ through $2^n - 1$. The negative elements of **M** are completed to the positive by adding $2^n$, so that each column of **M** will comprise the numbers from 1 through $2^n - 1$. For each number, we have 8 pairs of indices [*i, j*] corresponding to the rows of **Y** and **P**. Multiplying these rows element-by-element yields the balanced sequence.

For the resolution *n* = 3, the matrices **U** and **P** are given in Table 3.6.1 and Table 3.6.2 respectively, while the balanced sequences are presented in Table 3.6.3.

Table 3.6.1.

| $U_1$ | $U_2$ | $U_3$ |
|---|---|---|
| 0 | 0 | 0 |
| 0 | 0 | 1 |
| 0 | 2 | 0 |
| 0 | 2 | 1 |
| 4 | 0 | 0 |
| 4 | 0 | 1 |
| 4 | 2 | 0 |
| 4 | 2 | 1 |

Table 3.6.2.

| $P_1$ | $P_2$ | $P_3$ |
|---|---|---|
| 1 | 1 | 1 |
| 1 | 1 | -1 |
| 1 | -1 | 1 |
| 1 | -1 | -1 |
| -1 | 1 | 1 |
| -1 | 1 | -1 |
| -1 | -1 | 1 |
| -1 | -1 | -1 |

Table 3.6.3: Balanced EXB sequences for $M_n$, $n = 1\ldots3$.

| № | $M_3 = 1/8$ | | | | $M_2 = 2/8$ | | | | $M_3 = 3/8$ | | | | $M_1 = 4/8$ | | | |
|---|---|---|---|---|---|---|---|---|---|---|---|---|---|---|---|---|
| | $A_0$ | $A_1$ | $A_2$ | $A_3$ | $A_0$ | $A_1$ | $A_2$ | $A_3$ | $A_0$ | $A_1$ | $A_2$ | $A_3$ | $A_0$ | $A_1$ | $A_2$ | $A_3$ |
| 1 | 0 | 0 | 0 | 1 | 0 | 0 | 1 | 0 | 0 | 0 | 1 | 1 | 1 | -1 | 0 | 0 |
| 2 | 0 | 0 | 1 | -1 | 0 | 1 | -1 | 0 | 0 | 1 | 0 | -1 | 0 | 1 | 0 | 0 |
| 3 | 0 | 0 | 0 | 1 | 0 | 0 | 1 | 0 | 1 | -1 | -1 | 1 | 1 | -1 | 0 | 0 |
| 4 | 0 | 1 | -1 | -1 | 1 | -1 | -1 | 0 | 0 | 1 | 0 | -1 | 0 | 1 | 0 | 0 |
| 5 | 0 | 0 | 0 | 1 | 0 | 0 | 1 | 0 | 0 | 0 | 1 | 1 | 1 | -1 | 0 | 0 |
| 6 | 0 | 0 | 1 | -1 | 0 | 1 | -1 | 0 | 1 | -1 | 0 | -1 | 0 | 1 | 0 | 0 |
| 7 | 0 | 0 | 0 | 1 | 0 | 0 | 1 | 0 | 0 | 1 | -1 | 1 | 1 | -1 | 0 | 0 |
| 8 | 1 | -1 | -1 | -1 | 1 | -1 | -1 | 0 | 1 | -1 | 0 | -1 | 0 | 1 | 0 | 0 |

Table 3.6.3: cont'd.

| № | $M_3 = 5/8$ | | | | $M_2 = 6/8$ | | | | $M_3 = 7/8$ | | | |
|---|---|---|---|---|---|---|---|---|---|---|---|---|
| | $A_0$ | $A_1$ | $A_2$ | $A_3$ | $A_0$ | $A_1$ | $A_2$ | $A_3$ | $A_0$ | $A_1$ | $A_2$ | $A_3$ |
| 1 | 1 | 0 | -1 | -1 | 1 | 0 | -1 | 0 | 1 | 0 | 0 | -1 |
| 2 | 1 | -1 | 0 | 1 | 1 | -1 | 1 | 0 | 1 | 0 | -1 | 1 |
| 3 | 0 | 1 | 1 | -1 | 1 | 0 | -1 | 0 | 1 | 0 | 0 | -1 |
| 4 | 1 | -1 | 0 | 1 | 0 | 1 | 1 | 0 | 1 | -1 | 1 | 1 |
| 5 | 1 | 0 | -1 | -1 | 1 | 0 | -1 | 0 | 1 | 0 | 0 | -1 |
| 6 | 0 | 1 | 0 | 1 | 1 | -1 | 1 | 0 | 1 | 0 | -1 | 1 |
| 7 | 1 | -1 | 1 | -1 | 1 | 0 | -1 | 0 | 1 | 0 | 0 | -1 |
| 8 | 0 | 1 | 0 | 1 | 0 | 1 | 1 | 0 | 0 | 1 | 1 | 1 |

Comparing Table 3.6.3 with Table 3.2.1 which contains all the EXB codes of the same $M_n$, we conclude that each balanced sequence can be formed by replicating and further arranging these EXB codes. Since the sign of each $A_j$ ($j > 0$) in Table 3.6.3 alternates, the charging and discharging of each flying capacitor $C_j$ are interleaved and this reduces output voltage ripple. Moreover, due to a constant number of zeros spacing nonzero $A_j$ the ripple reduction will be significant.

## 3.7 The EXB based SCC in step-up mode

*"Algebra is generous; she often gives more than is asked of her."*
Jean le Rond d'Alembert

In this section, we demonstrate how the step-down EXB based SCC can be utilized for step-up conversion. Considering the fact that system (3.5.3) describing the step-down SCC with $M_3 = 3/8$ is solvable it should also be solvable if the indices of $V_o$ and $V_{in}$ are interchanged. Such a manipulation leads to the system (3.7.1) and from the hardware viewpoint means switching the input and output of the step-down SCC as depicted in Fig 3.7.1.

$$\begin{cases} V_o - V_1 + 0 - V_3 = V_{in} \\ 0 + V_1 + 0 - V_3 = V_{in} \\ V_o - V_1 - V_2 + V_3 = V_{in} \\ 0 + V_1 - V_2 + V_3 = V_{in} \\ 0 + 0 + V_2 + V_3 = V_{in} \end{cases} \quad (3.7.1)$$

The solution of (3.7.1) is:

$$\begin{aligned} V_1 &= 4/3 \cdot V_{in} \\ V_2 &= 2/3 \cdot V_{in} \\ V_3 &= 1/3 \cdot V_{in} \\ V_o &= 8/3 \cdot V_{in} \end{aligned} \quad (3.7.2)$$

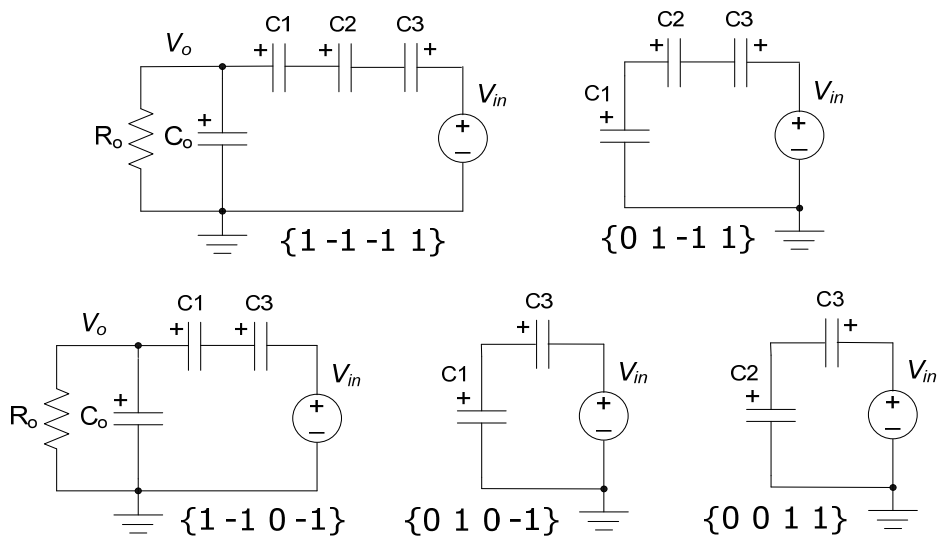

Figure 3.7.1: Topologies of the step-up SCC reciprocal to the case of $M_3 = 3/8$.

Generally, the conversion ratios of the step-up EXB based SCC with resolution $n$ are reciprocal to their step-down counterparts $M_n$ and defined by a set of fractions with numerator $2^n$ and denominators $1, \ldots, 2^n - 1$. These fractions have no resolution in the common sense and behave as $1/x$. For $n = 1\ldots5$, the step-up conversion ratios $1/M_n$ are depicted in Fig. 3.7.2.

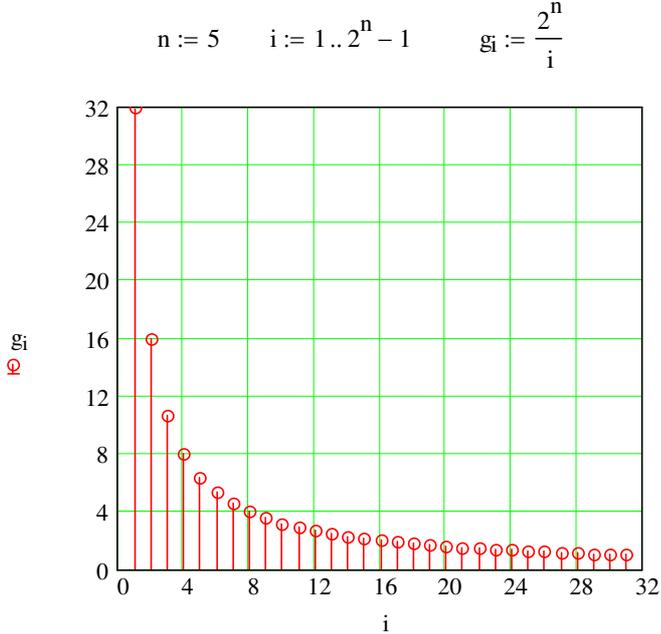

Figure 3.7.2: Step-up conversion ratios $1/M_n$, $n = 1\ldots5$.

Note that the highest conversion ratio in this SCC embodiment is equal to $2^n$. Although a number of step-up SCC with the conversion ratio $2^n$ have been proposed earlier [43–46], [52] there is no published report of a SCC with intermediate binary conversion ratios.

## 3.8 Some investigations into redundancy

*"Perfection is achieved, not when there is nothing more to add,*
*but when there is nothing left to take away."*

Antoine de Saint Exupery

As shown above, for the EXB fraction with resolution $n$, we can spawn $w \geq n + 1$ codes, which are translated to $w$ different SCC topologies. Applying the KVL to these topologies yields a system of $w$ linear equations with $n + 1$ unknowns. Since the self-adjusting voltages correspond to a unique solution of this system, the number of linearly independent equations in the system must be $n + 1$. Thus, we have $w - n - 1$ redundant equations, which can be eliminated as well as the corresponding SCC topologies.

For the system written in matrix form we can apply the conventional Gaussian or Gauss-Jordan elimination. However, these methods modify the coefficient matrix, while the aim is to point out its linearly dependent rows. It is also desirable to obtain the expression for each redundant row, which can be produced by more complicated methods [34], [35] that use a memory matrix and special pivoting technique. Any nonzero matrix **A** can be reduced by Gaussian elimination to an infinite number of row echelon forms by using different sequences of row operations. However, all row echelon forms of **A** correspond to exactly one matrix, which is called the reduced row echelon form and obtained by Gauss-Jordan elimination.

The proposed method to identify dependent rows of the coefficient matrix **A** is based on the fact that the row space of **A** is the column space of its transposed counterpart **A'**, while the rank of both matrices is the same.

For example: the EXB based SCC with conversion ratio $M_3 = 3/8$ is described by the system $\mathbf{Ax} = \mathbf{b}$, where:

$$\mathbf{A} = \begin{bmatrix} -1 & -1 & 1 & -1 \\ 1 & -1 & 1 & -1 \\ -1 & 0 & -1 & -1 \\ 1 & 0 & -1 & -1 \\ 0 & 1 & 1 & -1 \end{bmatrix} \quad \mathbf{x} = \begin{bmatrix} x_1 \\ x_2 \\ x_3 \\ x_4 \end{bmatrix} \quad \mathbf{b} = \begin{bmatrix} -1 \\ 0 \\ -1 \\ 0 \\ 0 \end{bmatrix} \quad (3.8.1)$$

Transposing **A** yields:

$$\mathbf{A'} = \begin{bmatrix} -1 & 1 & -1 & 1 & 0 \\ -1 & -1 & 0 & 0 & 1 \\ 1 & 1 & -1 & -1 & 1 \\ -1 & -1 & -1 & -1 & -1 \end{bmatrix} \quad (3.8.2)$$

Subjecting **A'** to Gauss-Jordan elimination we obtain its reduced row echelon form:

$$\mathbf{F} = \begin{bmatrix} 1 & 0 & 0 & -1 & 0 \\ 0 & 1 & 0 & 1 & 0 \\ 0 & 0 & 1 & 1 & 0 \\ 0 & 0 & 0 & 0 & 1 \end{bmatrix} \quad (3.8.3)$$

Let us consider in **F** the columns where "1" is leading; their indices correspond to the indices of independent rows of **A**. These rows can be identified simply when the matrix **F** is taken by absolute value and then summed over the columns. The resulting vector **s** will comprise elements greater than one (3.8.4), while their indices will be the indices of dependent rows of **A**.

$$\mathbf{s} = \begin{bmatrix} 1 & 1 & 1 & 3 & 1 \end{bmatrix} \quad (3.8.4)$$

The only element in (3.8.4) that is greater than one is 3 which has index 4, so the fourth row of **A** can be safely eliminated.

Since the flying capacitors in the EXB based SCC are always connected in series, the charge delivered in each SCC topology depends on the number of nonzero coefficients $A_j$ in the corresponding row of **A**. Sorting the rows of **A** by the number of zeros in descending order before executing the above elimination procedure allows one to reduce both the adjustment duration and the equivalent resistor (detailed in the following).

# 4. PROPOSED CLASS OF SCC WITH ARBITRARY RESOLUTION

## 4.1 Generic fractional (GFN) representation

In the above proposed EXB based SCC, the number of target voltages is dependent on the resolution $n$ and equal to $2^n - 1$. At each target voltage, the SCC efficiency reaches a maximum value. However, the efficiency drops when the desired output voltage lies between or outside the target voltages. This drawback can be alleviated by increasing the resolution $n$ and, consequently, the number of flying capacitors. Another approach is to introduce one or several additional target voltages between neighbor EXB target voltages at the same number of flying capacitors.

In order to realize this approach we have developed a Generic Fractional Number (GFN) representation, where the radix is not restricted by 2 as in the EXB representation, but can take an arbitrary integer value. In the GFN based SCC, the voltages across the flying capacitors are defined by the corresponding powers of the radix, while the resolution determines the gap between neighboring GFN target voltages. Both the radix and the resolution can be made higher by increasing the number of flying capacitors.

For the resolution $n$ and the radix $r$, consider a set of fractions $N_n(r)$ in the range (0, 1) with numerators 1, 2, ..., $r^n - 1$ and denominator $r^n$. Any fraction $N_n(r)$ can be represented in the next form:

$$N_n(r) = A_0 + \sum_{j=1}^{n} A_j r^{-j} \qquad (4.1.1)$$

where $r$ is an integer greater than one, $A_0$ can be either 0 or 1, and $A_j$ can take any of the values: -($r$-1), ..., -1, 0, 1, ..., ($r$-1). Expression (4.1.1) defines the Generic Fractional Numbers (GFN) representation, which differs from the conventional high-radix representations (e.g. decimal) since $A_j$ can be negative. Because $A_j$ can take any of the values $1-r$, ..., $r-1$, the GFN representation is akin to the Generalized Signed Digit (GSD) representation, for example:

$$\begin{aligned} 5 &= 1 \cdot 3^2 - 1 \cdot 3^1 - 1 \cdot 3^0 \rightarrow \{1\ \text{-}1\ \text{-}1\} \\ 5 &= 0 \cdot 3^2 + 2 \cdot 3^1 - 1 \cdot 3^0 \rightarrow \{0\ \ 2\ \text{-}1\} \\ 5 &= 1 \cdot 3^2 - 2 \cdot 3^1 + 2 \cdot 3^0 \rightarrow \{1\ \text{-}2\ \ 2\} \end{aligned} \qquad (4.1.2)$$

As seen from (4.1.2), the GSD representation for a given integer is not unique and this property is used mostly for carry-save fast computer arithmetic. We have modified the GSD representation for the fractions $N_n(r)$ limited in the range (0, 1). As a result, the coefficient $A_0$ is not allowed to be "-1". Because of the redundancy coming from the GSD representation, any $N_n(r)$ can be represented by a number of GFN codes, for example $N_2(3) = 5/9$:

$$\begin{aligned}
5/9 &= 1 - 1 \cdot 3^{-1} - 1 \cdot 3^{-2} \rightarrow \{1\ \text{-}1\ \text{-}1\} \\
5/9 &= 0 + 2 \cdot 3^{-1} - 1 \cdot 3^{-2} \rightarrow \{0\ \ 2\ \text{-}1\} \\
5/9 &= 1 - 2 \cdot 3^{-1} + 2 \cdot 3^{-2} \rightarrow \{1\ \text{-}2\ \ 2\}
\end{aligned} \qquad (4.1.3)$$

In the next section, we provide a simple procedure to spawn all the GFN codes for a given fraction $N_n(r)$. This procedure will be followed by a number of corollaries, which are crucial to define and explain the operation of the GFN based SCC.

## 4.2 Spawning the GFN codes and its corollaries

As mentioned above, the GFN code of $N_n(r)$ with non-negative coefficients $A_j$ is identical to the representation of $N_n(r)$ in the conventional number system with radix $r$. This code is called hereinafter the original code. In order to generate all the GFN codes corresponding to a given GFN fraction $N_n(r)$ within the range (0, 1), we use a procedure that involves adding and subtracting the coefficient $A_j = r - 1$ to the original code.

**Spawning the GFN codes.** This procedure is iterative and starts from any $A_j > 0$ in the original code of $N_n(r)$. Adding $r - 1$ to this $A_j$ results in $A_j < r - 1$ and "1" from the left as the carry. To maintain the value of $N_n(r)$, we subtract $r - 1$ from the obtained $A_j$, and spawn thereby a new GFN code. The procedure repeats for all $A_j > 0$ in the original code and for all $A_j > 0$ in each spawned GFN code.

In example (4.2.1), three alternative GFN codes are spawned from the original code of $N_2(3) = 4/9$. The GFN codes for other $N_n(3)$ where $n = 1, 2$ are summarized in Table 4.2.1.

$$
\begin{array}{ccc}
\downarrow & \downarrow & \downarrow \\
3^0\ 3^{-1}\ 3^{-2} & 3^0\ 3^{-1}\ 3^{-2} & 3^0\ 3^{-1}\ 3^{-2} \\
\phantom{+}0\ \ 1\ \ 1 & \phantom{+}0\ \ 1\ \ 1 & \phantom{+}0\ \ 2\ -2 \\
+\ 0\ \ 0\ \ \mathbf{2} & +\ 0\ \ \mathbf{2}\ \ 0 & +\ 0\ \ \mathbf{2}\ \ 0 \\
\hline
\phantom{+}0\ \ 2\ \ 0 & \phantom{+}1\ \ 0\ \ 1 & \phantom{+}1\ \ 1\ -2 \\
+\ 0\ \ 0\ \mathbf{-2} & +\ 0\ \mathbf{-2}\ \ 0 & +\ 0\ \mathbf{-2}\ \ 0 \\
\hline
\phantom{+}0\ \ 2\ -2 & \phantom{+}1\ -2\ \ 1 & \phantom{+}1\ -1\ -2
\end{array}
\qquad (4.2.1)
$$

**Corollary 1:** For resolution $n$, the minimum number of GFN codes is $n + 1$.

This is because each of the $A_j > 0$ in the original code with resolution $n$ generates a new GFN code and a carry. Further iterations cause the carry to propagate, so that each "0" in the original code is turned to "1", which is also operated on to spawn a new GFN code. So, the minimum number of GFN codes is the original code plus $n$ that is, $n + 1$.

**Corollary 2:** Each $A_j > 0$ in either the original or the GFN code yields at least one $A_j < 0$ in the same position $j$ of another GFN code. This is because the spawning procedure involves subtracting $r - 1$ from $A_j < r - 1$.

Both the above corollaries are very important and, as detailed in the following, provide the self-adjusting target voltage $N_n(r) \cdot V_{in}$ at the output of the GFN based SCC, irrespectively of the values of the capacitors used.

Table 4.2.1: The GFN codes of $N_n(3)$, $n = 1, 2$.

| $N_2(3) = 1/9$ | | | $N_2(3) = 2/9$ | | | $N_1(3) = 3/9$ | | | $N_2(3) = 4/9$ | | |
|---|---|---|---|---|---|---|---|---|---|---|---|
| $A_0$ | $A_1$ | $A_2$ | $A_0$ | $A_1$ | $A_2$ | $A_0$ | $A_1$ | $A_2$ | $A_0$ | $A_1$ | $A_2$ |
| 1 | -2 | -2 | 1 | -2 | -1 | 1 | -2 | 0 | 1 | -1 | -2 |
| 0 | 1 | -2 | 0 | 1 | -1 | 0 | 1 | 0 | 0 | 2 | -2 |
| 0 | 0 | 1 | 0 | 0 | 2 | | | | 1 | -2 | 1 |
| | | | | | | | | | 0 | 1 | 1 |

Table 4.2.1: cont'd.

| $N_2(3) = 5/9$ | | | $N_1(3) = 6/9$ | | | $N_2(3) = 7/9$ | | | $N_2(3) = 8/9$ | | |
|---|---|---|---|---|---|---|---|---|---|---|---|
| $A_0$ | $A_1$ | $A_2$ | $A_0$ | $A_1$ | $A_2$ | $A_0$ | $A_1$ | $A_2$ | $A_0$ | $A_1$ | $A_2$ |
| 1 | -1 | -1 | 1 | -1 | 0 | 1 | 0 | -2 | 1 | 0 | -1 |
| 0 | 2 | -1 | 0 | 2 | 0 | 1 | -1 | 1 | 1 | -1 | 2 |
| 1 | -2 | 2 | | | | 0 | 2 | 1 | 0 | 2 | 2 |
| 0 | 1 | 2 | | | | | | | | | |

## 4.3 Combinatorial method to obtain the GFN codes

Due to the spawning procedure described in the previous section, we have derived the important properties of the GFN codes. However, from the viewpoint of performance, this procedure is slow because each GFN code of $N_n(r)$ is obtained by the series, digit-by-digit, adding and subtracting the coefficient $A_j = r - 1$. The alternative combinatorial method proposed in this section is parallel and therefore faster than the previous one.

In accordance with the above definition, the GFN representation of $N_n(r)$ contains $n$ coefficients $A_j$, which takes any of $2r - 1$ values: 1−r, ..., r−1. We consider all combinations of these values arranged at $n$ positions as a matrix **M** of $(2r-1)^n$ rows by $n$ columns. This matrix is obtained by a full factorial design, where each level is $2r-1$, and the number of levels is $n$. Note that **M** defines the representations of the numbers $\frac{1-(2r-1)^n}{2}, \ldots, \frac{(2r-1)^n-1}{2}$ in the balanced number system with radix $2r - 1$.

For the sake of an exact integer calculation, we multiply both the sides of the GFN formula (4.1.1) by $r^n$. As a result, each GFN weight $r^{-j}$ is replaced by $r^{n-j}$ and we have $n$ powers of the radix $r$, which compose a column-vector **K** of length $n$. Multiplying the matrix **M** by this vector gives a column-vector **F** of length $(2r-1)^n$. To indicate $A_0 = 0$ and $A_0 = 1$ in the GFN codes we introduce a column-vector **B** of the same length. The positive elements of the vector **F** correspond to $A_0 = 0$ and transfer to the vector **B** as zeros, while the negative elements correspond $A_0 = 1$ and transfer to **B** as ones.

We complete the negative elements of the vector **F** to the positive by adding $r^n$ and obtain a column-vector **F'**. This vector contains $(2r-1)^n$ elements, which are the numerators $m = 1, \ldots, r^n - 1$ of all $N_n(r)$. The search for certain $m$ results in several row indexes, while the same rows of **B** and **M** compose the GFN codes for given $N_n(r)$.

For radix $r = 3$ and resolution $n = 2$ the proposed combinatorial method is demonstrated step-by-step in the following, while the full factorial design GFN matrix and the vectors used are shown in Fig. 4.3.1.

1) All combinations of five values -2, -1, 0, 1, 2 in two positions are defined by the full factorial design matrix **M** of $5^2$ rows by 2 columns obtained with the MATLAB command fullfact([5 5])–3.

2) The column-vector **K** = [9; 3; 1], so that the product of **M** and **K** is the column-vector **F** of length $5^2$ comprising the numbers from $1-3^2$ through $3^2-1$.

3) The positive elements of the vector **F** are transferred to the vector **B** as zeros, while the negative elements are transferred to **B** as ones.

4) At the same time, the negative numbers in **F** are completed to positive by adding $3^2$, so that the completed vector **F'** contains the numbers from 0 through $3^2 - 1$.

5) Since **F'** contains $5^2$ elements, any $m = 1, \ldots, 3^2 - 1$ appears in **F'** more than once, and search for certain $m$ results in several row indexes. The same rows of **B** and **M** compose the GFN codes of $N_n(3)$, $n = 1, 2$. Such a gathering is demonstrated in Fig. 4.3.1, where $m = 4$ that corresponds to $N_2(3) = 4/9$.

Comparing the obtained codes with the codes presented in Table 4.2.1 we conclude that the combinatorial method yields the same result as the procedure spawning the GFN codes.

$$\mathbf{B} = \begin{bmatrix} 1 \\ 1 \\ 1 \\ \mathbf{1} \\ 1 \\ \mathbf{1} \\ 1 \\ 1 \\ 1 \\ 1 \\ 1 \\ 1 \\ 0 \\ 0 \\ 0 \\ 0 \\ 0 \\ 0 \\ \mathbf{0} \\ 0 \\ \mathbf{0} \\ 0 \\ 0 \\ 0 \\ 0 \end{bmatrix} \quad \mathbf{M} = \begin{bmatrix} -2 & -2 \\ -2 & -1 \\ -2 & 0 \\ \mathbf{-2} & \mathbf{1} \\ -2 & 2 \\ \mathbf{-1} & \mathbf{-2} \\ -1 & -1 \\ -1 & 0 \\ -1 & 1 \\ -1 & 2 \\ 0 & -2 \\ 0 & -1 \\ 0 & 0 \\ 0 & 1 \\ 0 & 2 \\ 1 & -2 \\ 1 & -1 \\ 1 & 0 \\ \mathbf{1} & \mathbf{1} \\ 1 & 2 \\ \mathbf{2} & \mathbf{-2} \\ 2 & -1 \\ 2 & 0 \\ 2 & 1 \\ 2 & 2 \end{bmatrix} \quad \mathbf{K} = \begin{bmatrix} 3 \\ 1 \end{bmatrix} \quad \mathbf{F} = \begin{bmatrix} -8 \\ -7 \\ -6 \\ -5 \\ -4 \\ -5 \\ -4 \\ -3 \\ -2 \\ -1 \\ -2 \\ -1 \\ 0 \\ 1 \\ 2 \\ 1 \\ 2 \\ 3 \\ 4 \\ 5 \\ 4 \\ 5 \\ 6 \\ 7 \\ 8 \end{bmatrix} \quad \mathbf{F'} = \begin{bmatrix} 1 \\ 2 \\ 3 \\ \mathbf{4} \\ 5 \\ \mathbf{4} \\ 5 \\ 6 \\ 7 \\ 8 \\ 7 \\ 8 \\ 0 \\ 1 \\ 2 \\ 1 \\ 2 \\ 3 \\ \mathbf{4} \\ 5 \\ \mathbf{4} \\ 5 \\ 6 \\ 7 \\ 8 \end{bmatrix}$$

Figure 4.3.1: Full factorial design GFN matrix of size 25×2 and used vectors.

## 4.4 Translating GFN codes to SCC topologies

For certain GFN fraction $N_n(r)$, we consider a step-down SCC system that includes a voltage source $V_{in}$, a set of $n \cdot (r-1)$ flying capacitors and an output capacitor $C_o$ connected in parallel with the load $R_o$. The flying capacitors are divided into $n$ groups of $r-1$ capacitors $C_j$ in each one. These groups and $V_{in}$ are connected in accordance with the GFN codes of $N_n(r)$ in such a way that $C_o$ is continuously charged. In particular, the GFN coefficient $A_0$ is responsible for the connection of $V_{in}$, while each group $j$ of $r-1$ flying capacitors $C_j$ is associated with the GFN coefficient $A_j$. Irrespective of the connection of $V_{in}$, the groups $j$ are always connected in series. Within each group $j$, the type of connection of the flying capacitors $C_j$ is determined by the absolute value $|A_j|$ and its complement $\Delta_j = r - |A_j| - 1$. In order to configure the GFN based SCC topologies we use the following rules.

1) If $A_0 = 1$ then $V_{in}$ is connected.
2) If $A_0 = 0$ then $V_{in}$ is not connected.
3) If $A_j < -1$ then $|A_j|$ capacitors $C_{jx}$ of group $j$ are connected in series with the same polarity and charged. The remaining $\Delta_j$ capacitors are connected in parallel and compose an "equalizing" capacitor $C_{je}$. This capacitor is connected in parallel to each $C_{jx}$ capacitor running thereby over the $|A_j|$ series connection.
4) If $A_j = -1$ then all $r-1$ capacitors $C_{jx}$ of group $j$ are connected in parallel and charged.
5) If $A_j = 0$ then all $r-1$ capacitors $C_{jx}$ of group $j$ are disconnected.
6) If $A_j = 1$ then all $r-1$ capacitors $C_{jx}$ of group $j$ are connected in parallel and discharged.
7) If $A_j > 1$ then $A_j$ capacitors $C_{jx}$ of group $j$ are connected in series with the same polarity and discharged. The remaining $\Delta_j$ capacitors are connected in parallel and compose an "equalizing" capacitor $C_{je}$. This capacitor is connected in parallel with each $C_{jx}$ capacitor running thereby over the $A_j$ series connection.

As an example we translate the GFN codes of $N_1(3) = 3/9$ and $N_2(3) = 4/9$ presented in Table 4.4.1 to the corresponding SCC topologies. Since in the first case of $N_1(3) = 3/9$ the resolution $n = 1$ and the radix $r = 3$, we need a single group of two flying capacitors $C_{1.1}$ and $C_{1.2}$, which is associated with the GFN coefficient $A_1$. Within this group the type of connection of $C_{1.1}$ and $C_{1.2}$ is determined by $|A_1|$ and $\Delta_1 = 2 - |A_1|$. Thus, each GFN code of $N_1(3) = 3/9$ leads to a specific SCC topology as depicted in Figure 4.4.1.

Table 4.4.1.

| $N_1(3) = 3/9$ | | | $N_2(3) = 4/9$ | | |
|---|---|---|---|---|---|
| $A_0$ | $A_1$ | $A_2$ | $A_0$ | $A_1$ | $A_2$ |
| 1 | -2 | 0 | 1 | -1 | -2 |
| 0 | 1 | 0 | 0 | 2 | -2 |
|  |  |  | 1 | -2 | 1 |
|  |  |  | 0 | 1 | 1 |

We assume that within the group $j = 1$ engaged in each SCC topology of Fig. 4.4.1, the flying capacitors $C_{1.1}$ and $C_{1.2}$ keep the same voltage $V_1 = 3^{-1} \cdot V_{in}$. Multiplying $V_{in}$ and $V_1$ by the corresponding coefficients $A_0$ and $A_1$ in the GFN codes of $N_1(3) = 3/9$, we find their algebraic sum, which is equal to the output target voltage $V_o = 3/9 \cdot V_{in}$.

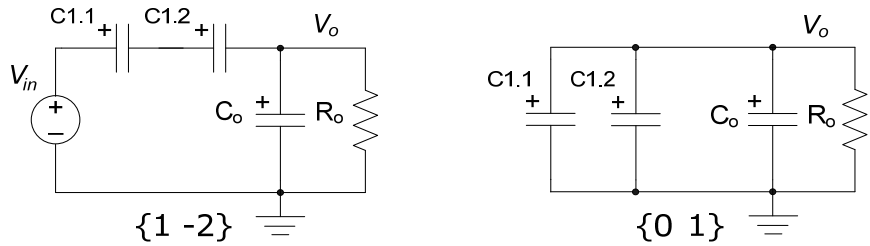

Figure 4.4.1: SCC topologies configured from the GFN codes of $N_1(3) = 3/9$.

Note that the SCC topologies of Fig. 4.4.1 represent the industry-standard SCC with the conversion ratio 1/3, which is actually equal to $N_1(3) = 3/9$. Because the industry-standard SCC uses a single group of the flying capacitors, they can be considered as a sub-class of the GFN based SCC with resolution $n = 1$. The topologies of this SCC sub-class are configured from the GFN codes of $N_1(r)$ where $|A_1|$ is equal to either 1 or $r - 1$. Substituting these values in the formula of the GFN representation yields the conversion ratios of this GFN based SCC sub-class as a pair of complementary fractions $N_1(r) = 1/r$ and $N_1(r) = 1 - 1/r$.

Let us now consider the more complicated case of $N_2(3) = 4/9$, where the resolution $n$ is increased to 2, while the radix $r$ is not changed and equal to 3. In this case, we need two groups of two flying capacitors in each. These groups are numbered $j = 1, 2$ and associated with the coefficients $A_1$ and $A_2$ respectively. Within each group $j$, the flying capacitors are indexed as $C_{j.1}$ and $C_{j.2}$, while the type of their connection is determined by $|A_j|$ and $\Delta_j = 2 - |A_j|$. Thus, in accordance with the above rules, each GFN code of $N_2(3) = 4/9$ leads to a specific SCC topology as depicted in Figure 4.4.2.

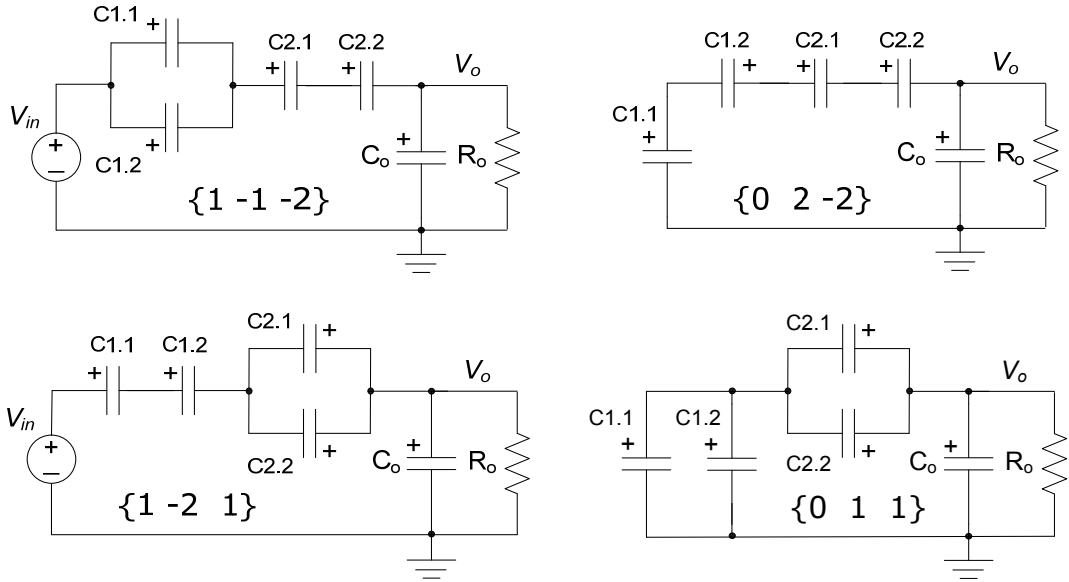

Figure 4.4.2: SCC topologies configured from the GFN codes of $N_2(3) = 4/9$.

To demonstrate that the output voltage $V_o$ in each SCC topology of Fig. 4.4.2 can, in principle, be equal to $4/9 \cdot V_{in}$ we assume that within the groups $j = 1, 2$ the flying capacitors $C_{1.x}$ and $C_{2.x}$ keep the voltages $V_1 = 3^{-1} \cdot V_{in}$ and $V_2 = 3^{-2} \cdot V_{in}$, respectively. As in the previous case, we multiply $V_{in}$, $V_1$, $V_2$ by the corresponding coefficients $A_0$, $A_1$, $A_2$ in the GFN codes of $N_2(3) = 4/9$ and find their algebraic sum, which is equal to the output target voltage $V_o = 4/9 \cdot V_{in}$.

Generally, translating all the GFN codes of certain $N_n(r)$ to the SCC topologies, we ought to obtain the target voltage $V_o = N_n(r) \cdot V_{in}$, under the condition that, within group $j$, each flying capacitor $C_{jx}$ keeps the voltage $V_j = r^{-j} \cdot V_{in}$. In the following we show that all the voltages in the GFN based SCC are self-adjusting to the above specified values and this property is due to Corollaries 1 and 2 of the procedure for spawning the GFN codes.

## 4.5 Self-adjusting voltages in the GFN based SCC

In this section we consider the GFN based SCC under the assumption that the voltages of the capacitors engaged in each SCC topology remain constant, but of unknown value. Applying the KVL to $w$ different SCC topologies we compose a system of $w$ linear equations. If this system has a unique solution, we obtain the target and the radix-$r$-weighted voltages across the output and the flying capacitors respectively.

Consider the GFN based SCC with the conversion ratio $N_1(3) = 1/3$. Its topologies reproduced from the example in the previous section are presented in Fig. 4.5.1.

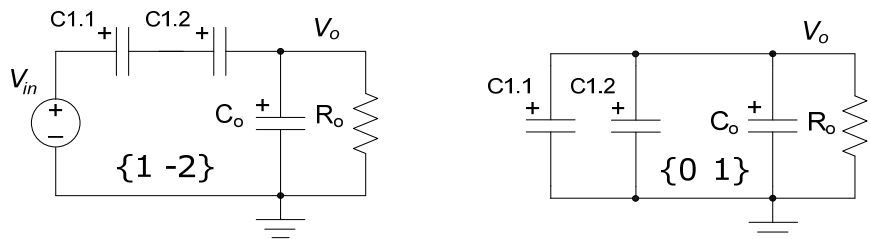

Figure 4.5.1: Topologies of the SCC with conversion ratio 1/3.

Designating the voltages across $C_{1.1}$, $C_{1.2}$ and $C_o$ by $V_1$, $V_{1.2}$ and $V_o$ respectively, we compose the next system of linear equations.

$$\begin{cases} V_{in} - V_{1.1} - V_{1.2} = V_o \\ V_{1.1} = V_{1.2} = V_o \end{cases} \qquad (4.5.1)$$

Because in (5.1.3) both $V_{1.1} = V_o$ and $V_{1.2} = V_o$, we conclude that $C_{1.1}$ and $C_{1.2}$ may be connected to $C_o$ independently, one after another in different time instants. On the other hand, due to the parallel connection of $C_{1.1}$ and $C_{1.2}$ we can introduce $V_1 = V_{1.1} = V_{1.2}$, so that the solution of (4.5.1) will be:

$$V_o = V_{1.1} = V_{1.2} = \frac{1}{3} V_{in} \qquad (4.5.2)$$

For the particular case when each coefficient $A_j$ takes $0, \pm 1, \pm(r-1)$ only, a system of equations for the GFN based SCC may be composed directly from the corresponding GFN codes. As an example, we show that the GFN codes of $N_2(3) = 4/9$ lead not only to the SCC topologies of Fig. 4.5.2, but also to the system of equations (4.5.3).

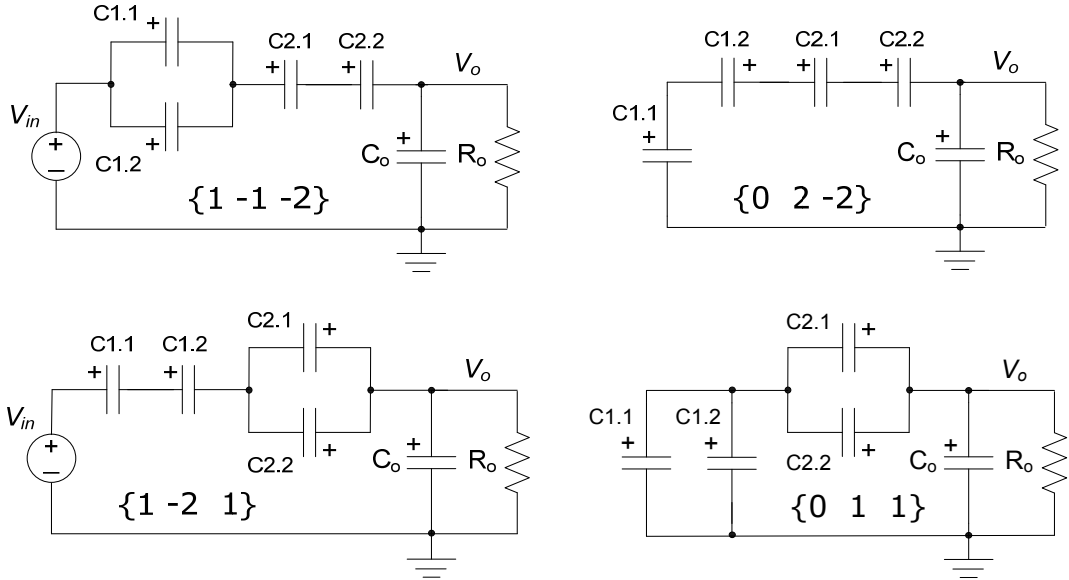

Figure 4.5.2: SCC topologies configured from the GFN codes of $N_2(3) = 4/9$.

$$\begin{cases} V_{in} - V_1 - 2V_2 = V_o \\ 0 + 2V_1 - 2V_2 = V_o \\ V_{in} - 2V_1 + V_2 = V_o \\ 0 + V_1 + V_2 = V_o \end{cases} \quad (4.5.3)$$

The number of equations in (4.5.3) is equal to 4 and equal to the number of all GFN codes of $N_2(3) = 4/9$. Adding one to the resolution $n = 2$ yields the number of unknowns equal to 3. Let us group the unknowns in (4.5.3) on the left hand side:

$$\begin{cases} -V_1 - 2V_2 - V_o = -V_{in} \\ 2V_1 - 2V_2 - V_o = 0 \\ -2V_1 + V_2 - V_o = -V_{in} \\ V_1 + V_2 - V_o = 0 \end{cases} \quad (4.5.4)$$

System (4.5.4) contains two non-zero free terms as the negative value of $V_{in}$. Generally, the connection of $V_{in}$ is provided by Corollary 1 of the procedure spawning the GFN codes as follows. Consider the original code of $N_n(r)$, where the coefficient $A_j$ can take only positive values 0, …, $r - 1$. Due to Corollary 1 the case of $A_j = 0$ is turned to $A_j = 1$, which is used to spawn the coefficient $A_0 = 1$ responsible for the connection of $V_{in}$. So, the GFN based SCC is described by a non-homogeneous system of linear equations.

Returning to (4.5.4), we normalize it to $V_{in}$:

$$\begin{cases} -1 \cdot x_1 - 2 \cdot x_2 - 1 \cdot x_3 = -1 \\ 2 \cdot x_1 - 2 \cdot x_2 - 1 \cdot x_3 = 0 \\ -2 \cdot x_1 + 1 \cdot x_2 - 1 \cdot x_3 = -1 \\ 1 \cdot x_1 + 1 \cdot x_2 - 1 \cdot x_3 = 0 \end{cases} \quad (4.5.5)$$

where

$$\begin{aligned} x_1 &= V_1/V_{in} \\ x_2 &= V_2/V_{in} \\ x_3 &= V_o/V_{in} \end{aligned} \quad (4.5.6)$$

We consider (4.5.5) in the form $\mathbf{Ax} = \mathbf{b}$, where

$$\mathbf{A} = \begin{bmatrix} -1 & -2 & -1 \\ 2 & -2 & -1 \\ -2 & 1 & -1 \\ 1 & 1 & -1 \end{bmatrix} \quad \mathbf{x} = \begin{bmatrix} x_1 \\ x_2 \\ x_3 \end{bmatrix} \quad \mathbf{b} = \begin{bmatrix} -1 \\ 0 \\ -1 \\ 0 \end{bmatrix} \quad (4.5.7)$$

Supplementing **A** with **b** we form the augmented matrix **A1** that allows us to investigate the solvability of (4.5.5).

$$\mathbf{A1} = \begin{bmatrix} -1 & -2 & -1 & -1 \\ 2 & -2 & -1 & 0 \\ -2 & 1 & -1 & -1 \\ 1 & 1 & -1 & 0 \end{bmatrix} \quad (4.5.8)$$

According to the Kronecker-Capelli theorem [31], [36], a nonhomogeneous system has a unique solution if and only if the rank of its coefficient matrix **A** is equal to the rank of its augmented matrix **A1** and equal to the number of unknowns.

For the resolution $n$, the number of columns in the matrix **A** is $n + 1$, while the number of rows is provided by Corollary 1 to be $w \geq n + 1$. According to the above theorem, we conclude that the rank of **A** as well as the rank of **A1** must be equal to $n + 1$, which is exactly the number of unknowns. So, we have to prove rigorously that the procedure spawning the GFN codes provides exactly $n + 1$ linearly independent rows (or columns) in both the matrices **A** and **A1**.

As mentioned above, the GFN based SCC uses $v = n \cdot (r - 1)$ flying capacitors. From a practical point of view, large values of $v$ imply that the SCC will be very expensive for realization and therefore we suppose $1 \leq v \leq 30$ as shown in Table 4.5.1.

Table 4.5.1.

| resolution, $n$ | radix, $r$ | number of caps, $v$ |
|---|---|---|
| 10 | 3 | 20 |
| 8 | 4 | 24 |
| 7 | 5 | 28 |
| 6 | 6 | 30 |

For each pair of $n$ and $r$ we calculated the rank of **A** and the rank of **A1** numerically in MATLAB, and ensured that the procedure for spawning the GFN codes leads to a system of $n+1$ linearly independent equations, while its unique solution is $N_n(r)$. Giving a rigorous theoretical proof of this for arbitrary values of n and r is planned for future work.

In the case of the system of equations (4.5.5), the rank of **A** is 4, which is equal to the rank of **A1** and to the number of unknowns. Thus, system (4.5.5) has a unique solution, which is found in the form $\mathbf{x} = \mathbf{A}^{-1}\mathbf{b}$:

$$\begin{aligned} x_1 &= 1/3 & V_1 &= 1/3 \cdot V_{in} \\ x_2 &= 1/9 & \text{or} \quad V_2 &= 1/9 \cdot V_{in} \\ x_3 &= 4/9 & V_o &= 4/9 \cdot V_{in} \end{aligned} \quad (4.5.9)$$

Since (4.5.9) is unique, the voltages $V_j$ ($j = 1…2$) across each $C_{jx}$ within group $j$ and the output voltage $V_o$ are self-adjusting to the specified values of $3^{-j} \cdot V_{in}$ and $4/9 \cdot V_{in}$.

Generally, for a given fraction $N_n(r)$, there are two ways to provide the self-adjusting voltages $V_j = r^{-j} \cdot V_{in}$ and $V_o = N_n(r) \cdot V_{in}$. One way is to configure only those SCC topologies, which correspond to $n+1$ linearly independent equations. The other way is to introduce redundant SCC topologies in addition to the $n + 1$ mentioned above.

By configuring these $w \geq n + 1$ different SCC topologies periodically we provide a continuous charge transfer through the load. An intuitive explanation for this is that, while the topologies change, each group of flying capacitors goes through a sequence of charging and discharging. This is assured by Corollary 2 of the procedure for spawning the GFN codes, which states that for each $A_j > 0$ in the GFN code there is at least one $A_j < 0$ in the same position. This means that the theory predicts that, in the GFN based SCC, all groups of flying capacitors are going through a sequence of charging and discharging.

According to the number *w* of different SCC topologies, we introduce the same number *w* of time intervals $t_k$, where $k = 1, \ldots, w$. During the interval $t_k$, the corresponding *k*-th topology does not change and repeats as depicted conceptually for $N_2(3) = 4/9$ in Fig. 4.5.3 with a period $T = \sum_{k=1}^{w} t_k$ at the intervals $t_{k+p \cdot w}$, where $p = 0, 1, 2, \ldots$ is the number of period.

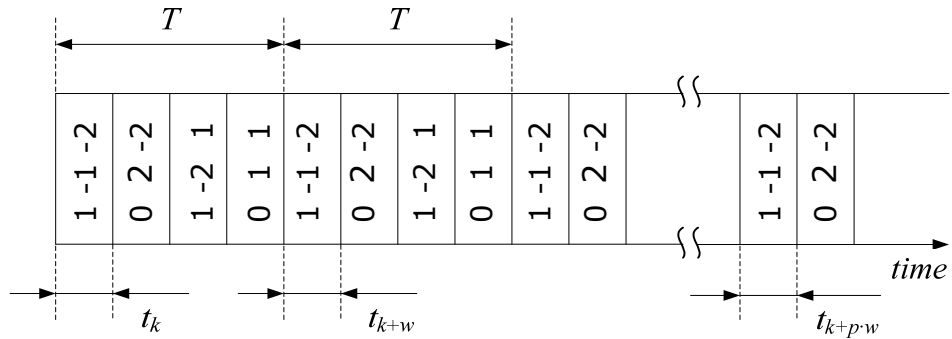

Figure 4.5.3: The perpetual GFN sequence of the SCC with $N_2(3) = 4/9$.

Irrespective of the order in which the topologies repeat, the voltages $V_j$ and $V_o$ eventually reach and stay at the specified values of $r^{-j} \cdot V_{in}$ and $N_n(r) \cdot V_{in}$ even if the SCC starts with zero or arbitrary voltages across the capacitors or when it is subjected to a disturbance. Due to this property, the GFN based SCC can be considered to be hardware for solving a system of linear equations by an iterative method.

## 4.6 The GFN based SCC in step-up mode

As shown above, the step-down EXB based SCC can be simply modified to operate in step-up mode. In this section we demonstrate the same for the step-down GFN based SCC described earlier by system (4.5.3), which is solvable even if the indices of $V_o$ and $V_{in}$ are interchanged. Such a manipulation leads to the system of equations (4.6.1) and from the hardware viewpoint means switching the input and output of the SCC as depicted in Fig 4.6.1.

$$\begin{cases} V_o - V_1 - 2V_2 = V_{in} \\ 0 + 2V_1 - 2V_2 = V_{in} \\ V_o - 2V_1 + V_2 = V_{in} \\ 0 + V_1 + V_2 = V_{in} \end{cases} \quad (4.6.1)$$

The solution of (4.6.1) is:

$$\begin{aligned} V_1 &= 3/4 \cdot V_{in} \\ V_2 &= 1/4 \cdot V_{in} \\ V_o &= 9/4 \cdot V_{in} \end{aligned} \quad (4.6.2)$$

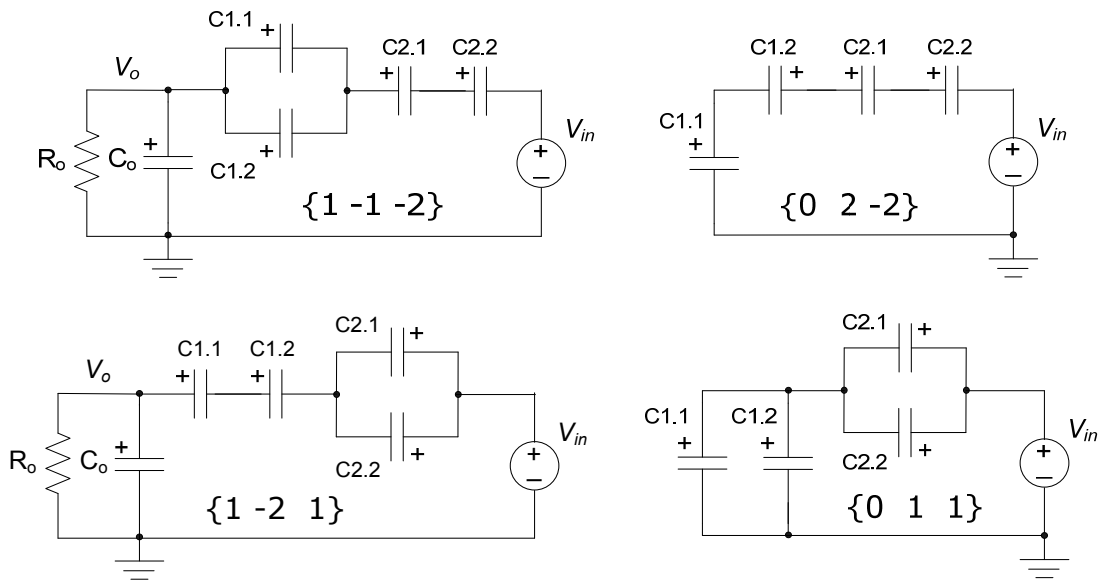

Figure 4.6.1: Topologies for the step-up SCC reciprocal to the case of $N_2(3) = 4/9$.

Generally, the conversion ratios of the step-up GFN based SCC with resolution $n$ are reciprocal to their step-down counterparts $N_n(r)$ and defined by a set of fractions with numerator $r^n$ and denominators $1, 2, \ldots, r^n - 1$. These fractions have no resolution in the common sense and behave as $1/x$. For $n = 1\ldots3$, the step-up conversion ratios $1/N_n(3)$ are depicted in Fig. 4.6.2.

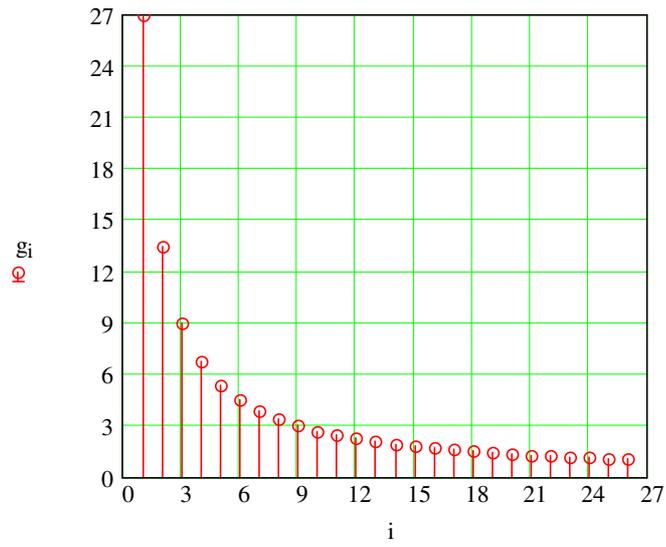

Figure 4.6.2: Step-up conversion ratios $1/N_n(3)$, $n = 1\ldots 3$.

# 5. PROPOSED NUMERICAL ANALYSIS

## 5.1 Investigating the Voltage Convergence Issue

*"Nothing is too wonderful to be true if it be consistent with the laws of nature."*
Michael Faraday

In this section we apply the charge conservation law and Kirchhoff's Voltage Law (KVL) to demonstrate that the SCC output voltage and the voltages across the flying capacitors converge to the target and the radix-*r*-weighted values respectively, under the condition that the charge redistributes immediately. This condition means, in particular, that the SCC load is disconnected, while all the capacitors and switches used are ideal. So, the considered SCC contains no resistive elements, while any SCC topology can be configured instantaneously and, as a consequence, a momentary transition from one SCC topology to another is allowed.

The charge conservation law states that charge can neither be created nor destroyed, only transferred. Applying this law to a SCC, we conclude that when capacitors holding initial voltages are engaged in a certain SCC topology, the existing charge is redistributed proportionally between the capacitances, so that the final voltages are balanced. These final voltages can be used as initial voltages in further iterations when the SCC topology is changed.

The KVL states that the algebraic sum of all voltages around any closed path of a circuit is zero. Since any SCC topology is a closed path circuit, we apply the KVL to the final voltages across the engaged capacitors. Due to the input voltage source $V_{in}$ engaged in some SCC topologies, the KVL algebraic sum can be equal to either $V_{in}$ or zero. Designating the current and next iterations by $i$ and $i+1$, we introduce a total number $m$ of capacitors in the $i+1$ iteration and express the final voltage $V_k^{i+1}$ across each engaged capacitor $C_k$ by initial voltage $V_k^i$ as follows:

$$\begin{cases} V_k^{i+1} = V_k^i \pm \dfrac{Q^{i+1}}{C_k} \\ \sum_{k=1}^{m} V_k^{i+1} = \left| \begin{array}{l} V_{in} \\ 0 \end{array} \right. \end{cases} \quad (5.1.1)$$

Expression (5.1.1) is a system of $m+1$ linear equation with $m+1$ unknowns, so it is solvable in principle.

Solving (5.1.1), we obtain the final voltages $V_k^{i+1}$ for use as the initial conditions in the next iteration. The iterations are done while the SCC topologies change according to either the EXB or the GFN codes. Because the SCC operation is determined by a periodic repetition of SCC topologies, at some iteration the steady-state needs to be reached.

It is important to distinguish between the terms "equilibrium" and "steady-state". Because of the ideal components that have been used, a solution of (5.1.1) assumes the equilibrium or charge (voltage) balance immediately after the SCC topology has been configured. The steady-state for the unloaded SCC with ideal components implies a constant (target) output voltage, while the charge transferred to the output capacitor after some iteration is zero. The number of iterations needed to reach the steady-state defines the adjustment duration.

As an example, we consider the EXB based SCC with the conversion ratio $M_3 = 3/8$. The input voltage source $V_{in}$ is engaged ($A_0 = 1$) in the SCC topology depicted in Fig. 5.1.1(a) and disengaged ($A_0 = 0$) in the SCC topology shown in Fig. 5.1.1(b).

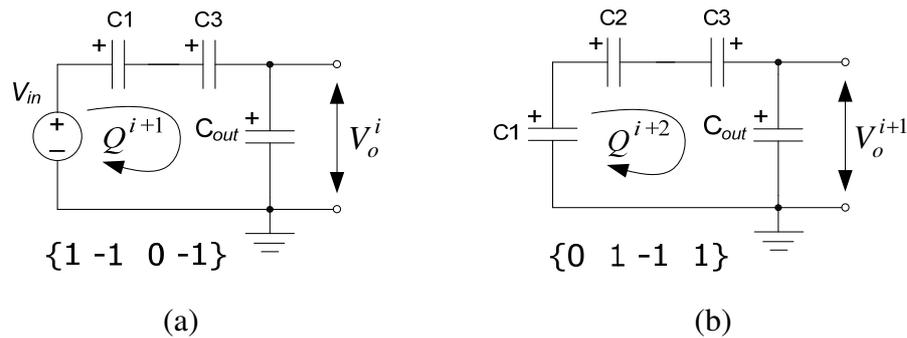

(a)                       (b)

Figure 5.1.1: Topologies of the SCC with conversion ratio $M_3 = 3/8$.

In the SCC topology of Fig. 5.1.1(a), let the initial voltages across the engaged capacitors C1, C2 and $C_{out}$ be $V_1^i$, $V_2^i$ and $V_o^i$, respectively. Using the EXB code {1 -1 0 -1}, we can write the KVL sum:

$$1 \cdot V_{in} - 1 \cdot V_1^{i+1} + 0 \cdot V_2^{i+1} - 1 \cdot V_3^{i+1} = V_o^{i+1} \tag{5.1.2}$$

where $V_1^{i+1}$, $V_2^{i+1}$ and $V_o^{i+1}$ are the final voltages. The total number of capacitors engaged in the SCC topology of Fig. 5.1.1(a) is $m = 3$. Applying the charge conservation law to each one, we have three equations of the kind

$$V_k^{i+1} = V_k^i \pm \frac{Q^{i+1}}{C_k} \tag{5.1.3}$$

The expressions (5.1.2) and (5.1.3) compose the system of $m+1=4$ linear equations with 4 unknowns. For given voltages $V_{in}$ and $V_k^i$ grouped on the right hand side, we obtain:

$$\begin{cases} V_1^{i+1} + 0 + 0 + 0 - \dfrac{Q^{i+1}}{C_1} = V_1^i \\ 0 + 0 + V_3^{i+1} + 0 - \dfrac{Q^{i+1}}{C_3} = V_3^i \\ 0 + 0 + 0 + V_o^{i+1} - \dfrac{Q^{i+1}}{C_o} = V_o^i \\ -V_1^{i+1} + 0 - V_3^{i+1} - V_o^{i+1} = -V_{in} \end{cases} \quad (5.1.4)$$

Because $A_2 = 0$, the capacitor C2 is disengaged and consequently $V_2^{i+1} = V_2^i$. Rewriting the system of equations (5.1.4) in matrix form:

$$\begin{bmatrix} 1 & 0 & 0 & 0 & -1/C_1 \\ 0 & 1 & 0 & 0 & 0 \\ 0 & 0 & 1 & 0 & -1/C_3 \\ 0 & 0 & 0 & 1 & -1/C_o \\ -1 & 0 & -1 & -1 & 0 \end{bmatrix} \times \begin{bmatrix} V_1^{i+1} \\ V_2^{i+1} \\ V_3^{i+1} \\ V_o^{i+1} \\ Q^{i+1} \end{bmatrix} = \begin{bmatrix} V_1^i \\ V_2^i \\ V_3^i \\ V_o^i \\ -V_{in} \end{bmatrix} \quad (5.1.5)$$

It is evident from (5.1.5) that the terms $1/C_j$ are multiplied by the corresponding $A_j$, as well as $V_j^{i+1}$ in the KVL sum. Since the EXB resolution is $n = 3$, the main diagonal of the coefficient matrix in (6.1.5) is formed simply from $n+1$ ones and zero. Solving (5.1.5), we obtain the voltages $V_k^{i+1}$ for use as the initial conditions (given voltages) in the next iteration defined by the EXB code {0 1 -1 1} and the SCC topology of Fig. 6.1.1 (b). The KVL sum for this case is:

$$1 \cdot V_1^{i+1} - 1 \cdot V_2^{i+1} + 1 \cdot V_3^{i+1} = V_o^{i+2} \quad (5.1.6)$$

The total number of capacitors engaged in the SCC topology of Fig. 5.1.1(b) did not changed from $m = 3$. Applying the charge conservation law we have three equations of kind

$$V_k^{i+2} = V_k^{i+1} \pm \dfrac{Q^{i+2}}{C_k} \quad (5.1.7)$$

As in the previous case we group the given voltages $V_k^{i+1}$ of (5.1.7) on the right-hand side and obtain the system of $m+1 = 4$ linear equations with 4 unknowns:

$$\begin{cases} V_1^{i+2} + 0 + 0 + 0 + \dfrac{Q^{i+2}}{C_1} = V_1^{i+1} \\ 0 + V_2^{i+2} + 0 + 0 - \dfrac{Q^{i+2}}{C_2} = V_2^{i+1} \\ 0 + 0 + V_3^{i+2} + 0 + \dfrac{Q^{i+2}}{C_3} = V_3^{i+1} \\ 0 + 0 + 0 + V_o^{i+2} - \dfrac{Q^{i+2}}{C_o} = V_o^{i+1} \\ V_1^{i+2} - V_2^{i+2} + V_3^{i+2} - V_o^{i+2} = 0 \end{cases} \qquad (5.1.8)$$

Rewriting (5.1.8) in matrix form, we have:

$$\begin{bmatrix} 1 & 0 & 0 & 0 & 1/C_1 \\ 0 & 1 & 0 & 0 & -1/C_2 \\ 0 & 0 & 1 & 0 & 1/C_3 \\ 0 & 0 & 0 & 1 & -1/C_o \\ 1 & -1 & 1 & -1 & 0 \end{bmatrix} \times \begin{bmatrix} V_1^{i+2} \\ V_2^{i+2} \\ V_3^{i+2} \\ V_o^{i+2} \\ Q^{i+2} \end{bmatrix} = \begin{bmatrix} V_1^{i+1} \\ V_2^{i+1} \\ V_3^{i+1} \\ V_o^{i+1} \\ 0 \end{bmatrix} \qquad (5.1.9)$$

The coefficient matrix in (5.1.9) is formed in the same manner as in (5.1.5) using the EXB code {0 1 -1 1} and its resolution $n = 3$. Solving (5.1.9), we obtain the voltages $V_k^{i+2}$, which are substituted as the initial conditions in the next iteration.

Using the above technique, we have conducted a convergence analysis of the considered SCC in MATLAB 7.1. The iterations were done for the EXB codes shown in Table 5.1.1.

Table 5.1.1

| \multicolumn{4}{c}{$M_3 = 3/8$} |
|---|---|---|---|
| $A_0$ | $A_1$ | $A_2$ | $A_3$ |
| 1 | -1 | -1 | 1 |
| 0 | 1 | -1 | 1 |
| 1 | -1 | 0 | -1 |
| 0 | 1 | 0 | -1 |
| 0 | 0 | 1 | 1 |

For the sake of convenience, we set the input voltage $V_{in} = 8\,[V]$, the identical flying capacitors $C1 = C2 = C3 = 4.7\,\mu F$ and the output capacitor $C_o = 470\,\mu F$. The initial voltages across all the capacitors are zero. As shown in Fig. 5.1.2, the voltages across the flying capacitors C1, C2 and C3 converge to binary weighted values 4, 2 and 1 [V], respectively. The output voltage converging to the target value 3[V] is depicted in Fig. 5.1.3, while Fig. 5.1.4 depicts the charge decaying to zero.

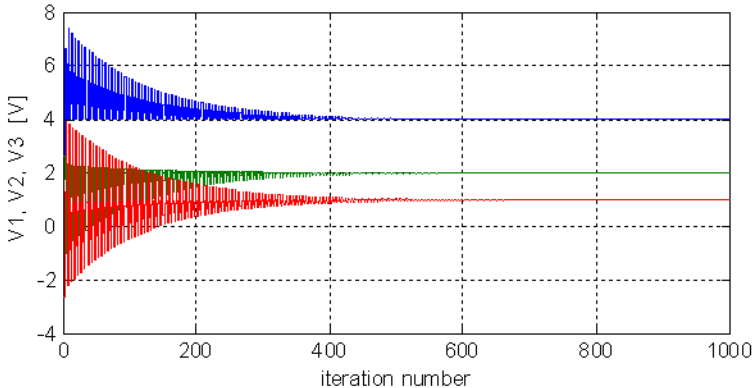

Figure 5.1.2: Convergence of the voltages $V_1$, $V_2$, $V_3$ at zero initial conditions.

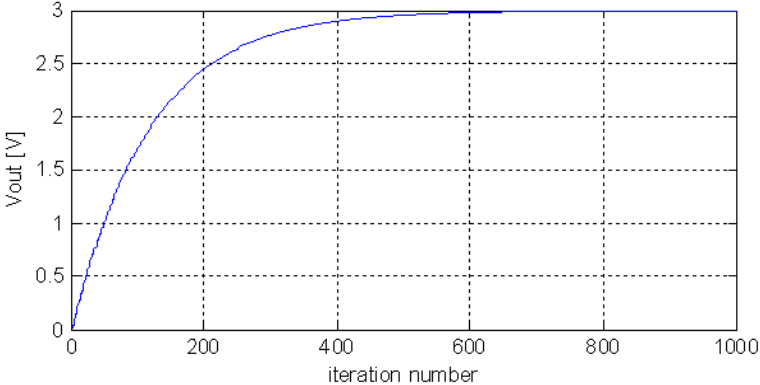

Figure 5.1.3: Convergence of the output voltage at zero initial conditions.

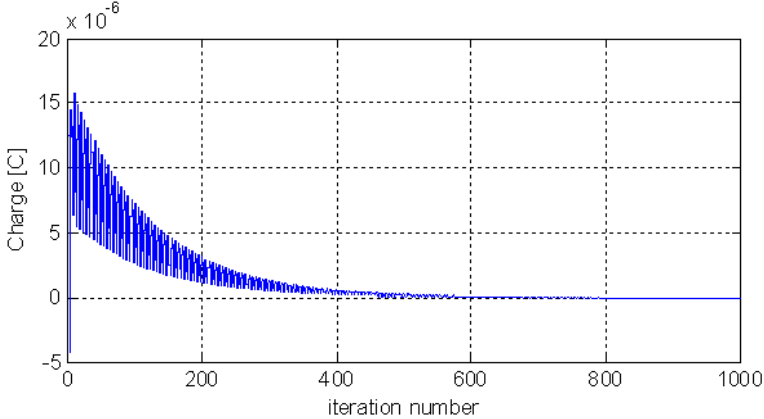

Figure 5.1.4: Decaying to zero charge at zero initial conditions.

Consider a case of non-zero initial conditions: the output capacitor is pre-charged to the initial voltage 1[V], and the flying capacitors C1, C2, C3 are pre-charged to the initial binary weighted voltages 4, 2, 1 [V] respectively. As shown in Fig. 5.1.5 the voltages across C1, C2, C3 change during the charge transfer, but return to their initial values when the steady-state is reached. The initial output voltage 1[V] converges to the target value 3[V] as shown in Fig. 5.1.5, while the decaying to zero charge is depicted in Fig. 5.1.6.

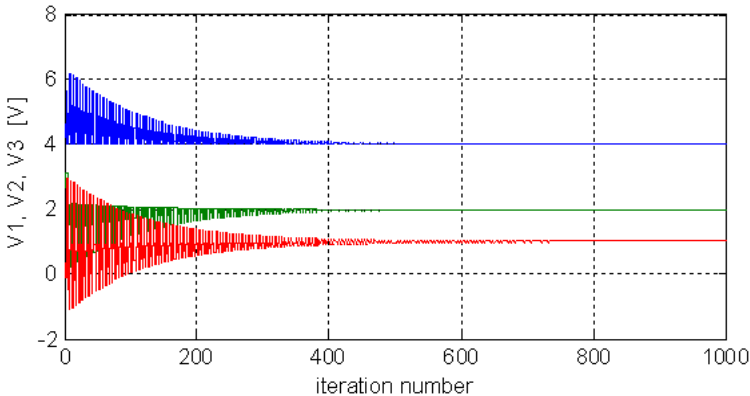

Figure 5.1.5: The voltages $V_1$, $V_2$, $V_3$ returning to binary weighted initial values.

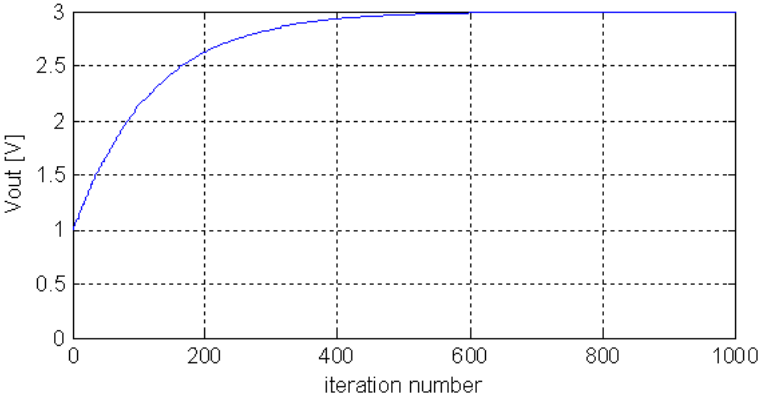

Figure 5.1.6: Convergence of the output voltage at binary weighted initial conditions.

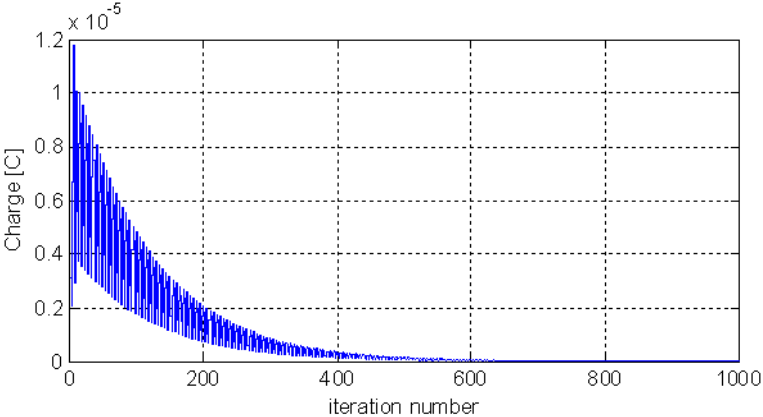

Figure 5.1.7: Decaying to zero charge at binary weighted initial conditions.

An additional way to display the convergence is a locus of the charge decaying to zero, where the number of axes corresponds to the number of SCC topologies, while each axis is the absolute value of the charge. The charge locuses for both cases of zero and binary weighted initial conditions are shown in Fig. 5.1.8. Because the locuses have a winding behavior, the voltages across the capacitors are self-adjusting to the values predicted by the EXB representation.

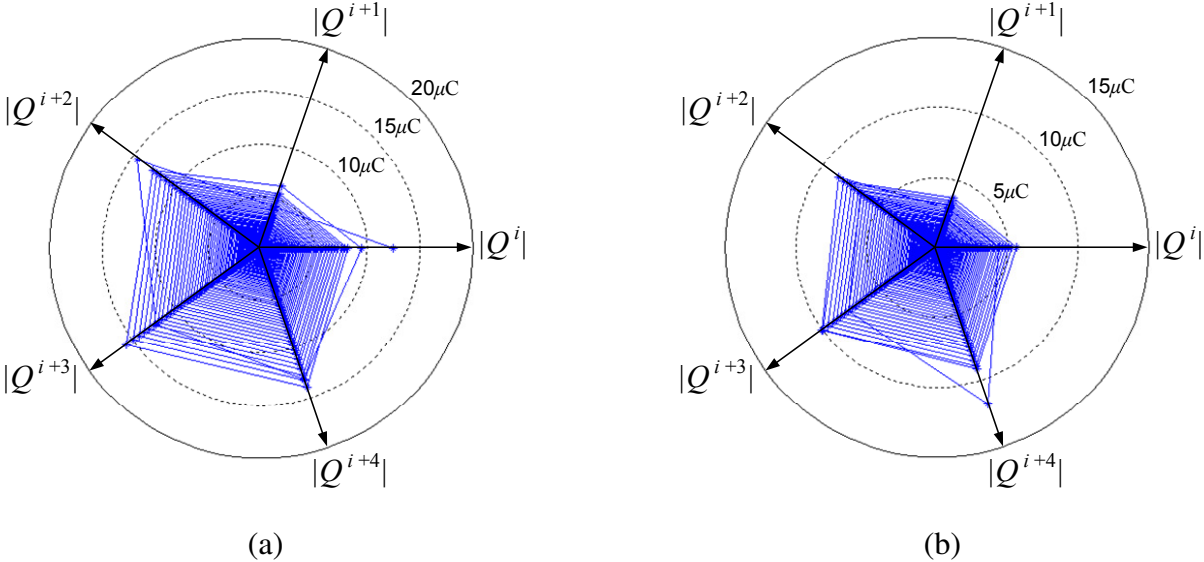

(a)          (b)

Figure 5.1.8: Charge locuses for zero (a) and binary weighted (b) initial conditions.

Due to the elimination procedure described in Section 5.5, we can safely omit the fourth, redundant SCC topology determined by EXB code {0 1 0 -1}. Convergence of the output voltage for this case and zero initial conditions is depicted in Fig. 5.1.9.

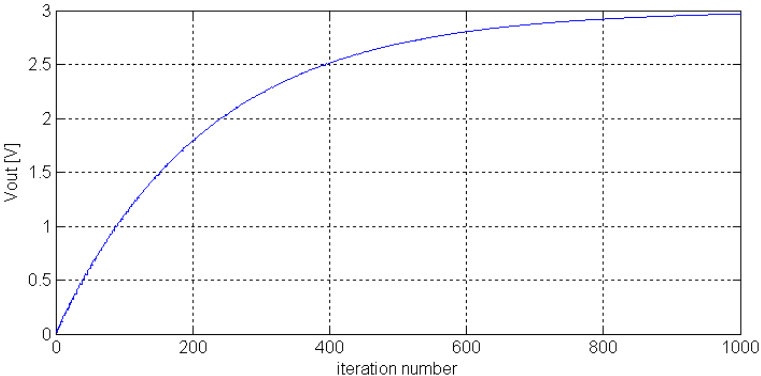

Figure 5.1.9: Convergence of $V_o$ when {0 1 0 -1} SCC topology is omitted.

Let us now consider the case when the elimination procedure described in Section 3.8 is executed for the EXB coefficients $A_j$ ($j$ >0) sorted by the number of zeros in descending order, so that the total capacitor in each SCC topology is increased.

The result of this elimination is presented in Table 5.1.2, while the output voltage converging to the target value at zero initial conditions is depicted in Fig. 5.1.10.

Table 5.1.2

| $M_3 = 3/8$ | | | |
|---|---|---|---|
| $A_0$ | $A_1$ | $A_2$ | $A_3$ |
| 1 | -1 | 0 | -1 |
| 0 | 1 | 0 | -1 |
| 0 | 0 | 1 | 1 |
| 1 | -1 | -1 | 1 |

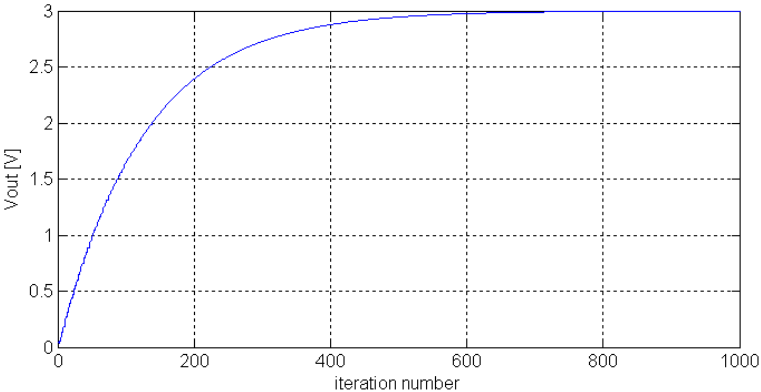

Figure 5.1.10: Convergence of $V_o$ when the total capacitor is increased.

For the EXB codes of Table 5.1.2 and the same flying capacitors, we change the output capacitor from $C_o = 470\mu F$ to $C_o = 220\mu F$. For zero initial conditions, the output voltage converges as shown in Fig. 5.1.11.

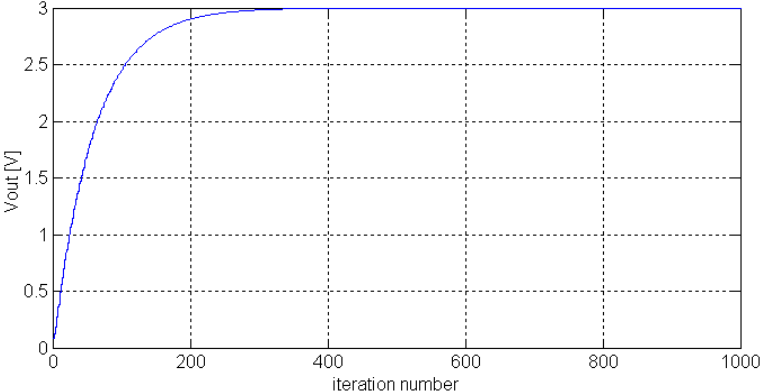

Figure 5.1.11: Convergence of $V_o$ when $C_o$ is changed from 470μF to 220μF.

For the same EXB codes, we set the flying capacitors to be $C1 = 40\mu F$, $C2 = 20\mu F$ and $C3 = 10\mu F$ that is binary weighted, while the output capacitor $C_o = 470\mu F$. Convergence of the output voltage for this case and zero initial conditions is depicted in Fig. 5.1.12.

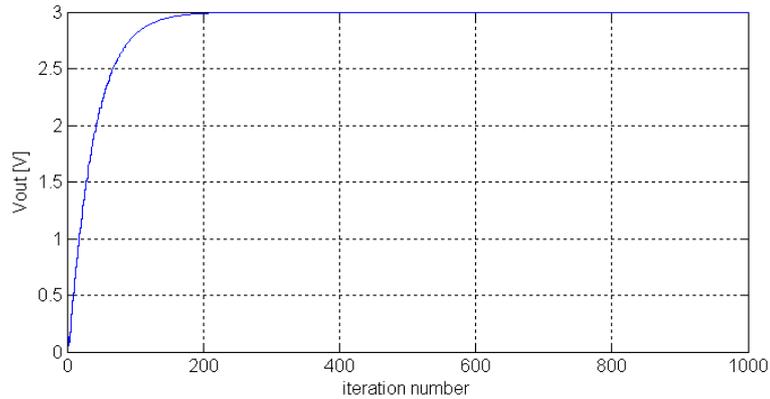

Figure 5.1.12: Convergence of $V_o$ when C1, C2, C3 are changed to be binary weighted.

As follows from Fig. 5.1.8 – Fig. 5.1.12, the adjustment duration depends on the ratios between the flying and the output capacitors and on which redundant SCC topologies were eliminated. From Fig. 5.1.2 and Fig. 5.1.5, we conclude that the voltages across the flying capacitors can change polarity and consequently, in the practical implementation of a SCC, these capacitors must be non-polarized.

## 5.2 Derivation of the Equivalent Resistor Expressions

> *"There exists everywhere a medium in things, determined by equilibrium."*
> Dmitri Mendeleev

The preceding studies [15–25] prove that the power loss in a SCC can be modeled by a single equivalent resistor $R_{eq}$. In this section we derive $R_{eq}$ for the EXB based SCC employing the periodic charge balance condition for each flying capacitor. The derivation is based on the generic and unified average model [26], [54].

Since any SCC in practice comprises parasitic resistances, the charge $Q$ is transferred by an exponentially decaying current $i(t) = I_0 \cdot e^{-t/RC}$ with the initial value $I_0 = \Delta V/R$, where $R$ is the total parasitic resistance. For a time interval $t_1$, we consider the current $i(t)$ as an average current $I_{av} = Q/t_1$ depicted schematically in Fig 5.2.1.

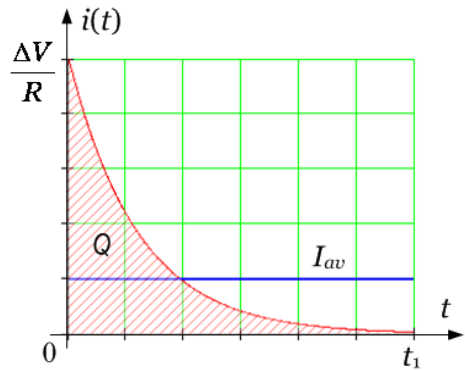

Figure 5.2.1: Exponential and average currents on a time interval.

In the steady state operation of a SCC, the charge received by a flying capacitor must be equal to the delivered charge. Assuming that the SCC topologies are configured for equal time intervals $t_1 = t_2 = \ldots = t_w = t$ and designating the average current in the $k$-th SCC topology by $I_k$, we conclude that the sum of average currents through each flying capacitor is equal to zero:

$$\sum_{k=1}^{w} Q_k = t \cdot \sum_{k=1}^{w} I_k = 0 \qquad (5.2.1)$$

This periodic charge balance condition is verified by simulation for the EXB based SCC with conversion ratio $M_3 = 3/8$ and demonstrated in Fig. 5.2.2 for the flying capacitor $C_1$ in the circuit of Fig. 6.1.2. Applying (5.2.1) to each flying capacitor $C_j$, we express each average current $I_k$ by the average output current $I_o$.

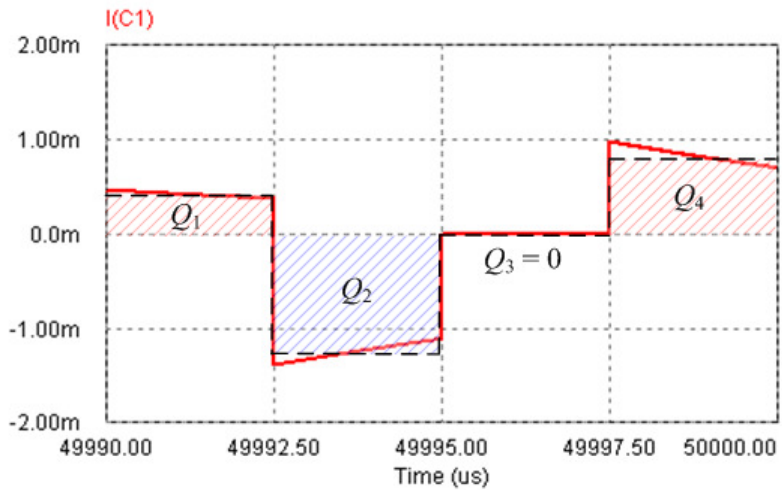

Figure 5.2.2: Charge balance for a single flying capacitor.

The previously spawned EXB codes of $M_3 = 3/8$ are presented in Table 5.2.1, while the SCC topologies are shown in Fig. 5.2.3. The directions of average currents for each flying capacitor $C_j$ are associated with the coefficients of column $A_j$ in Table 5.2.1. Multiplying each column element-by-element by the currents $I_1, \ldots, I_5$ we obtain Table 5.2.2.

Table 5.2.1.

| $M_3 = 3/8$ | | | |
|---|---|---|---|
| $A_0$ | $A_1$ | $A_2$ | $A_3$ |
| 1 | -1 | -1 | 1 |
| 0 | 1 | -1 | 1 |
| 1 | -1 | 0 | -1 |
| 0 | 1 | 0 | -1 |
| 0 | 0 | 1 | 1 |

Table 5.2.2.

| $M_3 = 3/8$ | | |
|---|---|---|
| C1 | C2 | C3 |
| $-I_1$ | $-I_1$ | $I_1$ |
| $I_2$ | $-I_2$ | $I_2$ |
| $-I_3$ | 0 | $-I_3$ |
| $I_4$ | 0 | $-I_4$ |
| 0 | $I_5$ | $I_5$ |

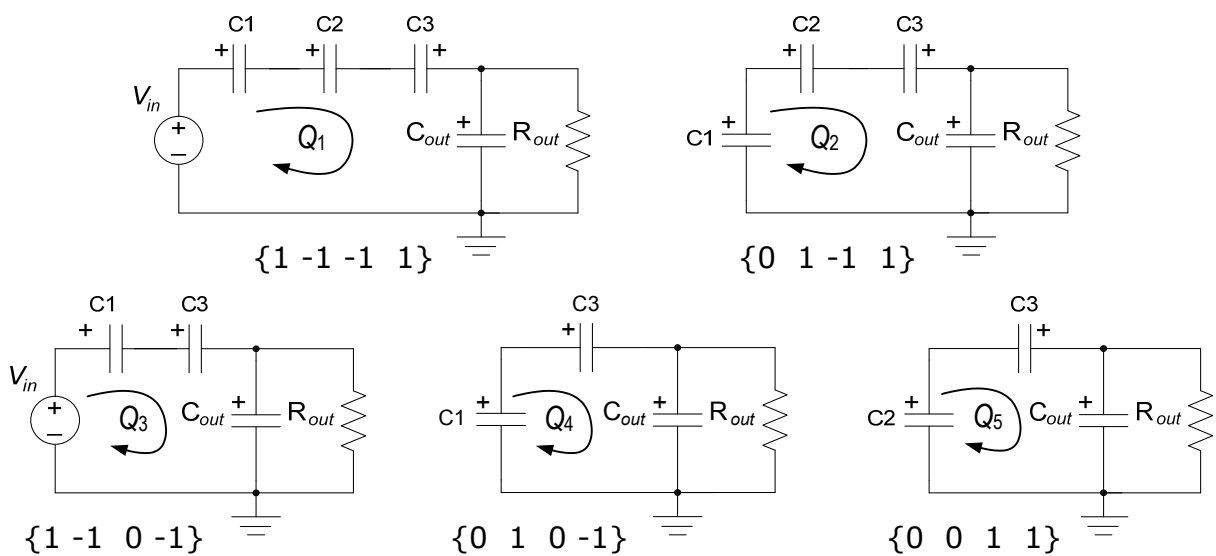

Figure 5.2.3: Topologies of the SCC with conversion ratio $M_3 = 3/8$.

Under the condition (5.2.1), we transpose Table 5.2.2 and obtain three expressions which are equal to zero. Because in each SCC topology of Fig. 5.2.3, the charge $Q_k$ is delivered to the output, an additional expression is $\sum_{k=1}^{w}|I_k| = I_o$. Thus, for the average currents $I_1, \ldots, I_5$ we have a system of linear equations:

$$\begin{cases} -I_1 + I_2 - I_3 + I_4 + 0 = 0 \\ -I_1 - I_2 + 0 + 0 + I_5 = 0 \\ I_1 + I_2 - I_3 - I_4 + I_5 = 0 \\ I_1 + I_2 + I_3 + I_4 + I_5 = I_o \end{cases} \quad (5.2.2)$$

Since Table 5.2.2 comprises redundant rows, the number of equations in (5.2.2) is less than the number of unknowns, and the system (5.2.2) has an infinite number of solutions:

$$\begin{aligned} I_1 - I_4 &= -1/8 \cdot I_o \\ I_2 + I_4 &= 3/8 \cdot I_o \\ I_3 + I_4 &= 1/2 \cdot I_o \\ I_5 &= 1/4 \cdot I_o \end{aligned} \quad (5.2.3)$$

We consider a particular solution of (5.2.2) when $I_4$ is equal to zero. In practice this means that the fourth row in Table 5.2.2 is eliminated and consequently, the fourth SCC topology is not configured at all. For each SCC topology we can find a total capacitor $C_k$ and a total resistor $R_k$, which are substituted to $\beta_k = \dfrac{t}{R_k C_k}$. According to the theory given in Section 2.4, the power dissipated in each topology is:

$$P_k = \frac{I_k^2 T_s}{2 C_k} \coth\left(\frac{\beta_k}{2}\right) \quad (5.2.4)$$

Since each $I_k$ is expressed by $I_o$, each $P_k$ depends on $I_o^2$ as well as the total dissipated power given by the sum of all $P_k$:

$$P = \sum_{k=1}^{w} P_k \quad (5.2.5)$$

We can rearrange (5.2.5) so that $I_o^2$ will be outside the brackets, reducing $I_o^2$ we obtain the equivalent resistor expression dependent of $R_k$ and $C_k$. In all the SCC topologies, the number of switches used is not changed and is equal to 4 as shown in Fig. 5.2.4.

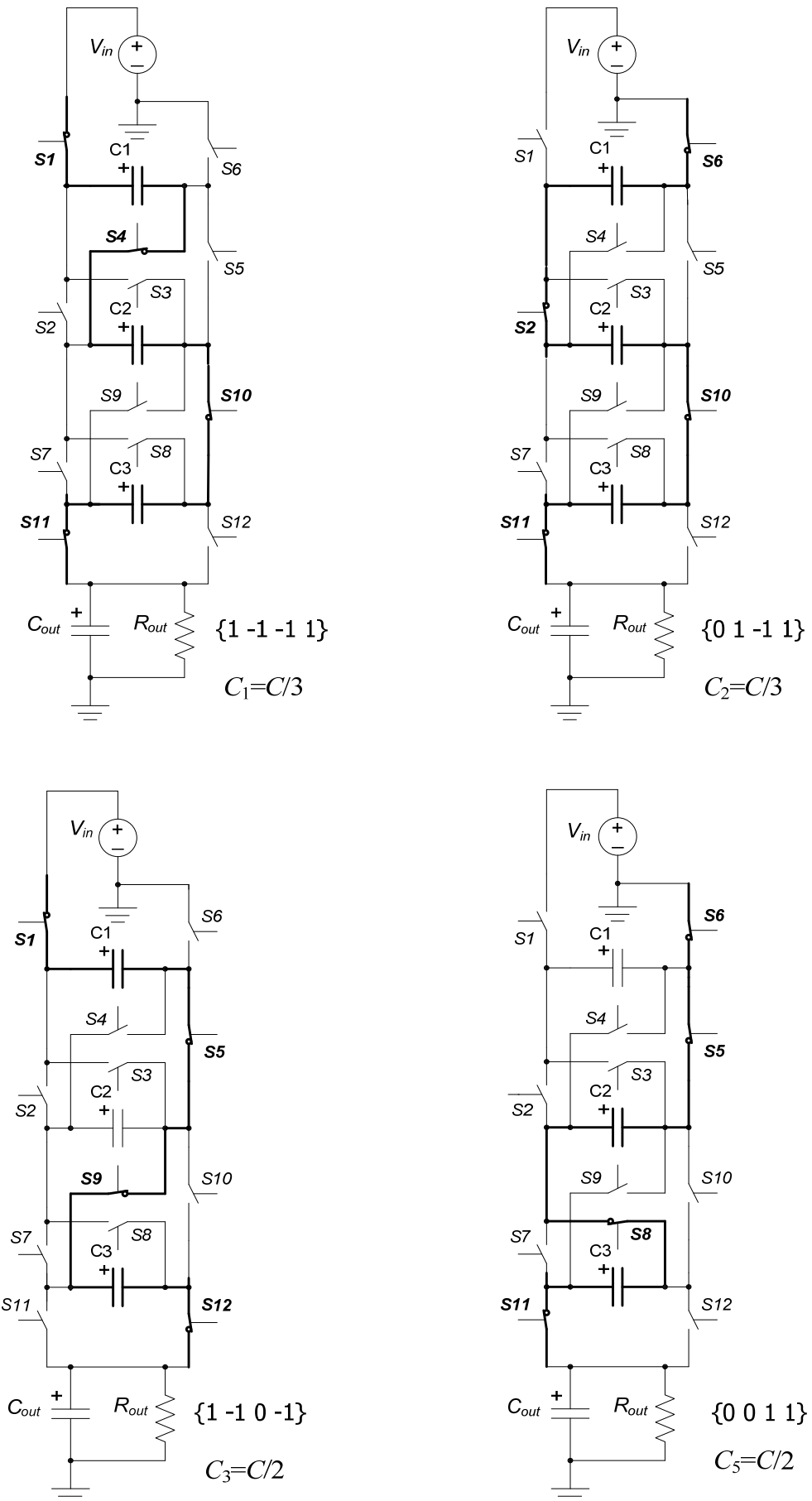

Figure 5.2.4: Switches used in each topology of the EXB based SCC with $M_3 = 3/8$.

Assuming an identical on-resistance *r* for all switches and neglecting other parasitic resistances (e.g. ESR), we define the total on-resistance as $R = 4r$. Since the flying capacitors are always connected in series, the total capacitors $C_k$ can be found using a number of non-zero coefficients $A_j$ ($j > 0$) in Table 5.2.1. In the experimental set-up, identical flying capacitors $C = 4.7\,\mu F$ are used, while $r = 1.2$ Ohm, and the switching frequency $f_s = 100$ kHz.

$$f_s := 100 \cdot 10^3 \quad T_s := \frac{1}{f_s} \quad t := \frac{T_s}{4} \quad R := 4 \cdot 1.2 \quad C := 4.7 \cdot 10^{-6} \quad \beta := \frac{t}{R \cdot C}$$

$$C_1 := \frac{C}{3} \qquad C_2 := \frac{C}{3} \qquad C_3 := \frac{C}{2} \qquad C_5 := \frac{C}{2}$$

$$\beta_1 := 3 \cdot \beta \qquad \beta_2 := 3 \cdot \beta \qquad \beta_3 := 2 \cdot \beta \qquad \beta_5 := 2 \cdot \beta$$

$$I_1 := -\frac{1}{8} \cdot I_o \qquad I_2 := \frac{3}{8} \cdot I_o \qquad I_3 := \frac{1}{2} \cdot I_o \qquad I_5 := \frac{1}{4} \cdot I_o$$

$$P_1 := \frac{I_1^{\,2} \cdot T_s}{2 \cdot C_1} \cdot \coth\!\left(\frac{\beta_1}{2}\right) \quad P_2 := \frac{I_2^{\,2} \cdot T_s}{2 \cdot C_2} \cdot \coth\!\left(\frac{\beta_2}{2}\right) \quad P_3 := \frac{I_3^{\,2} \cdot T_s}{2 \cdot C_3} \cdot \coth\!\left(\frac{\beta_3}{2}\right) \quad P_5 := \frac{I_5^{\,2} \cdot T_s}{2 \cdot C_1} \cdot \coth\!\left(\frac{\beta_5}{2}\right)$$

$$P := \frac{I_o^{\,2} \cdot T_s}{2 \cdot C} \cdot \left[\left(-\frac{1}{8}\right)^2 \cdot 3 \cdot \coth\!\left(\frac{3 \cdot \beta}{2}\right) + \left(\frac{3}{8}\right)^2 \cdot 3 \cdot \coth\!\left(\frac{3 \cdot \beta}{2}\right) + \left(\frac{1}{2}\right)^2 \cdot 2 \cdot \coth(\beta) + \left(\frac{1}{4}\right)^2 \cdot 2 \cdot \coth(\beta)\right]$$

$$P := \frac{I_o^{\,2} \cdot T_s}{2 \cdot C} \cdot \frac{10}{64} \cdot \left(3 \cdot \coth\!\left(\frac{3 \cdot \beta}{2}\right) + 4 \cdot \coth(\beta)\right)$$

So, the equivalent resistor for the SCC conversion ratio $M_3 = 3/8$ is:

$$R_{eq} = \frac{5}{64} \frac{1}{f_s C}\left[3\coth\!\left(\frac{3}{2}\beta\right) + 4\coth(\beta)\right] \tag{5.2.6}$$

Employing $T_s = 4RC\beta$ we can rewrite (5.2.6) in the following form:

$$R_{eq} = \frac{5}{16} R\beta \cdot \left[3\coth\!\left(\frac{3}{2}\beta\right) + 4\coth(\beta)\right] \tag{5.2.7}$$

When beta tends to zero, the expression (5.2.7) is reduced to:

$$\lim_{\beta \to 0} R_{eq} = \frac{15}{8} R = \frac{15}{2} r \tag{5.2.8}$$

The above derivation reveals that the equivalent resistor is inversely proportional to the total capacitor in each SCC topology and, consequently, depends on which redundant topologies were eliminated. The elimination procedure described in Section 3.8 allows increasing the total capacitors when the EXB coefficients $A_j$ ($j > 0$) are sorted by the number of zeros in descending order as shown in Table 5.2.3. For both $M_3 = 3/8$ and $M_3 = 5/8$ we eliminate the fifth row.

Table 5.2.3: Sorted EXB codes of $M_n$, $n = 1\ldots3$.

| $M_3 = 1/8$ | | | | $M_2 = 2/8$ | | | | $M_3 = 3/8$ | | | | $M_1 = 4/8$ | | | |
|---|---|---|---|---|---|---|---|---|---|---|---|---|---|---|---|
| $A_0$ | $A_1$ | $A_2$ | $A_3$ | $A_0$ | $A_1$ | $A_2$ | $A_3$ | $A_0$ | $A_1$ | $A_2$ | $A_3$ | $A_0$ | $A_1$ | $A_2$ | $A_3$ |
| 0 | 0 | 0 | 1 | 0 | 0 | 1 | 0 | 1 | -1 | 0 | -1 | 1 | -1 | 0 | 0 |
| 0 | 0 | 1 | -1 | 1 | -1 | -1 | 0 | 0 | 1 | 0 | -1 | 0 | 1 | 0 | 0 |
| 1 | -1 | -1 | -1 | 0 | 1 | -1 | 0 | 0 | 0 | 1 | 1 | | | | |
| 0 | 1 | -1 | -1 | | | | | 1 | -1 | -1 | 1 | | | | |
| | | | | | | | | *0* | *1* | *-1* | *1* | | | | |

Table 5.2.3: cont'd.

| $M_3 = 5/8$ | | | | $M_2 = 6/8$ | | | | $M_3 = 7/8$ | | | |
|---|---|---|---|---|---|---|---|---|---|---|---|
| $A_0$ | $A_1$ | $A_2$ | $A_3$ | $A_0$ | $A_1$ | $A_2$ | $A_3$ | $A_0$ | $A_1$ | $A_2$ | $A_3$ |
| 1 | 0 | -1 | -1 | 1 | 0 | -1 | 0 | 1 | 0 | 0 | -1 |
| 1 | -1 | 0 | 1 | 1 | -1 | 1 | 0 | 1 | 0 | -1 | 1 |
| 0 | 1 | 0 | 1 | 0 | 1 | 1 | 0 | 1 | -1 | 1 | 1 |
| 1 | -1 | 1 | -1 | | | | | 0 | 1 | 1 | 1 |
| *0* | *1* | *1* | *-1* | | | | | | | | |

The corresponding coefficients for $I_k$ and $C_k$ are summarized in Table 5.2.4.

Table 5.2.4: Coefficients required in $R_{eq}$ derivation.

| k | $M_3 = 1/8$ | | $M_2 = 2/8$ | | $M_3 = 3/8$ | | $M_1 = 4/8$ | |
|---|---|---|---|---|---|---|---|---|
| | $I_k/I_o$ | $C_k/C$ | $I_k/I_o$ | $C_k/C$ | $I_k/I_o$ | $C_k/C$ | $I_k/I_o$ | $C_k/C$ |
| 1 | 1/2 | 1 | 1/2 | 1 | 1/8 | 1/2 | 1/2 | 1 |
| 2 | 1/4 | 1/2 | 1/4 | 1/2 | 3/8 | 1/2 | 1/2 | 1 |
| 3 | 1/8 | 1/3 | 1/4 | 1/2 | 1/4 | 1/2 | | |
| 4 | 1/8 | 1/3 | | | 1/4 | 1/3 | | |

Table 5.2.4: cont'd.

| k | $M_3 = 5/8$ | | $M_2 = 6/8$ | | $M_3 = 7/8$ | |
|---|---|---|---|---|---|---|
| | $I_k/I_o$ | $C_k/C$ | $I_k/I_o$ | $C_k/C$ | $I_k/I_o$ | $C_k/C$ |
| 1 | 1/4 | 1/2 | 1/2 | 1 | 1/2 | 1 |
| 2 | 1/8 | 1/2 | 1/4 | 1/2 | 1/4 | 1/2 |
| 3 | 3/8 | 1/2 | 1/4 | 1/2 | 1/8 | 1/3 |
| 4 | 1/4 | 1/3 | | | 1/8 | 1/3 |

Assuming an identical total resistor *R* in the SCC topologies and using Table 5.2.3 and Table 5.2.4, the equivalent resistor expressions were derived for all the ratios $M_n$, $n = 1, \ldots, 3$ as presented in Table 5.2.5. An important issue in this derivation is that the same equivalent resistor is obtained for a pair of complementary conversion ratios $M_n$ and $1 - M_n$.

Table 5.2.5: Equivalent resistors for all the ratios $M_n$, $n = 1, \ldots, 3$.

| $M_n$ | Equivalent resistor expression | $\lim_{\beta \to 0} R_{eq}$ | $R_{eq}$, Ohm |
|---|---|---|---|
| $M_3 = 1/8$, $M_3 = 7/8$ | $R_{eq} = \dfrac{T_s}{64C}\left[8\coth\left(\dfrac{\beta}{2}\right) + 4\coth(\beta) + 3\coth\left(\dfrac{3}{2}\beta\right)\right]$ | $\dfrac{11}{8}R$ | 6.615 |
| $M_2 = 2/8$, $M_2 = 6/8$ | $R_{eq} = \dfrac{T_s}{8C}\left[\coth\left(\dfrac{\beta}{2}\right) + \coth(\beta)\right]$ | $\dfrac{9}{8}R$ | 5.42 |
| $M_3 = 3/8$, $M_3 = 5/8$ | $R_{eq} = \dfrac{T_s}{32C}\left[7\coth(\beta) + 3\coth\left(\dfrac{3}{2}\beta\right)\right]$ | $\dfrac{9}{8}R$ | 5.428 |
| $M_1 = 4/8$ | $R_{eq} = \dfrac{T_s}{4C}\coth\left(\dfrac{\beta}{2}\right)$ | $R$ | 4.82 |

# 6. Simulation Results

*"Science is what we understand well enough to explain to a computer.*
*Art is everything else we do."*
Donald E. Knuth

## 6.1. Verification of the equivalent resistor values

In order to configure the topologies of the EXB based SCC, each flying capacitor needs to have three types of connections {-1, 0, 1}. Since the pair {-1, 1} is responsible for the connection polarity, we consider a bridge switched circuit. Compose two capacitor bridges and connect one to another by rotation, the obtained double bridge cascade is presented in Fig. 6.1.1.

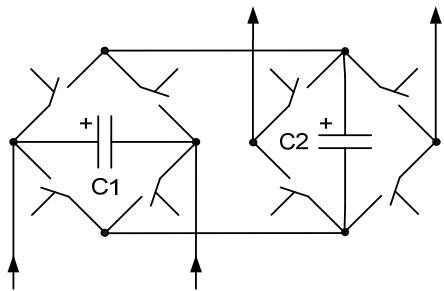

Figure 6.1.1: Double-bridge cascade.

The adjacent switches in the cascade allow either capacitor to be disabled, while the opposite switches can parallel a capacitor with the neighbor. So, this cascade provides not only three types of connections {-1, 0, 1}, but can also be used for building the GFN based SCC.

Using the double-bridge cascade, four additional switches and one additional capacitor, we compose the simulation circuit of the EXB based SCC with resolution $n = 3$ as depicted in Fig. 6.1.2. The simulations are done in the PSIM 7 simulator especially designed for power electronics, motor control, and dynamic system. Some of its advantages include fast simulation and lack of convergence problems. The PSIM bidirectional switches are ideal, i.e. have zero internal resistance, whereas the calculations done in the previous section assume an identical on-resistance $r = 1.2\ \Omega$ of all the switches. So, we need to introduce sub-circuits that comprise the ideal switch and a series resistor of 1.2 Ω as shown in Fig. 6.1.3.

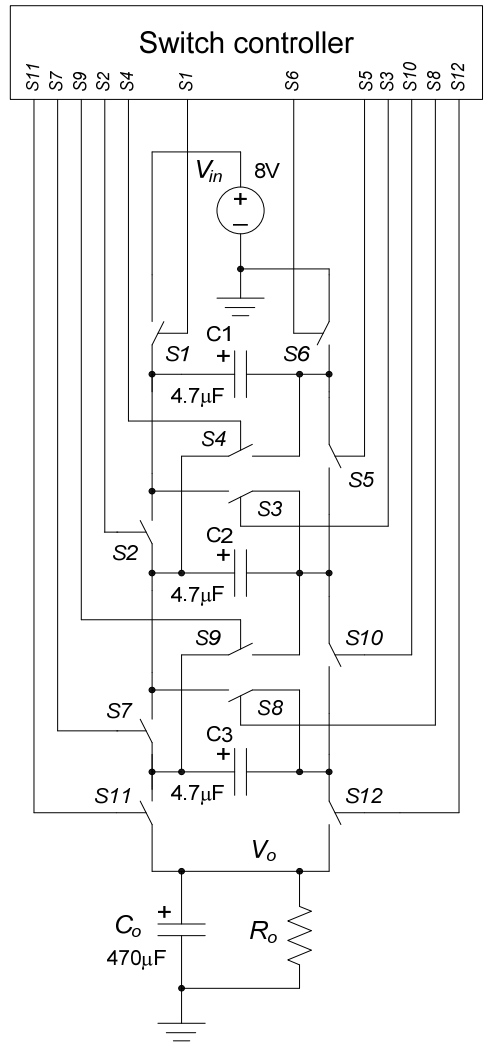

Figure 6.1.2: Simulation circuit for the EXB based SCC.

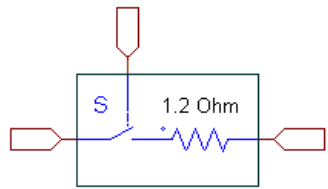

Figure 6.1.3: Switch subcircuit.

The switch controller in the circuit of Fig. 6.1.2 is built using the PSIM embedded C Script Block, which is clocked by an internal oscillator. The duration of each clock is the switching period $T_s$ = 10μs divided by $n + 1$, where $n \leq 3$ is the resolution of the EXB based SCC.

To verify the equivalent resistor values given in Table 5.2.5 we run the simulations for the EXB codes given in Table 5.2.3 and measure the output voltage $V_o$ for different load resistances $R_o$ = 100, 200, …, 500Ω when the SCC has reached the steady-state.

Fig. 6.1.4 depicts the measured voltages $V_o$ for the case of $M_1 = 4/8$. The measurements for the complementary ratios $M_n$ and $1 - M_n$ are presented side by side in Fig. 6.1.5 - 6.1.7.

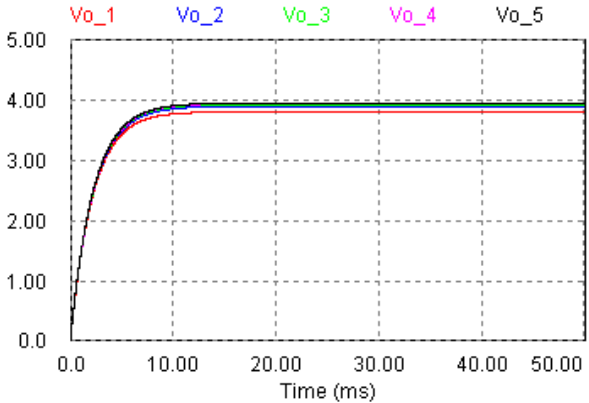

Figure 6.1.4: Simulation result for $M_1 = 4/8$.

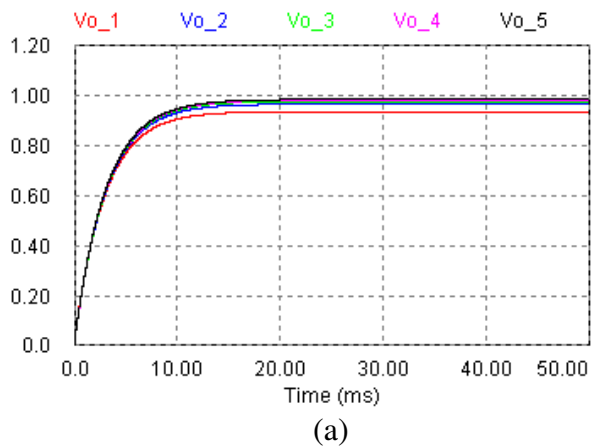 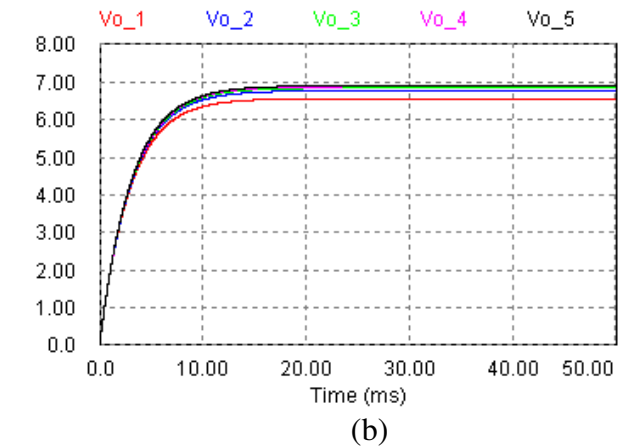

(a) (b)

Figure 6.1.5: Simulation result for $M_3 = 1/8$ (a) and $M_3 = 7/8$ (b).

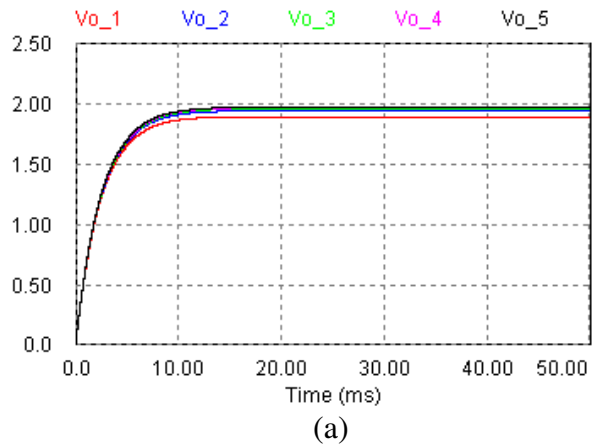 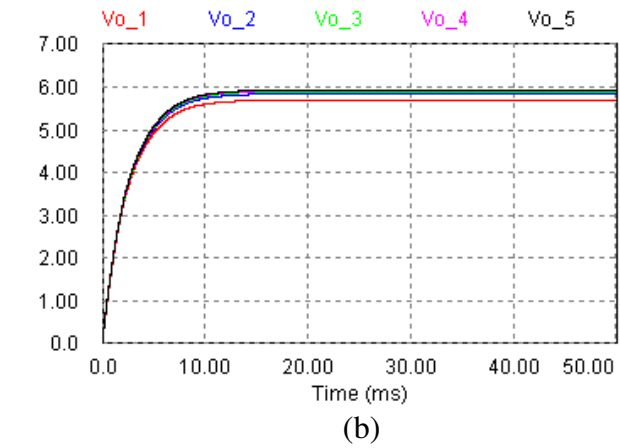

(a) (b)

Figure 6.1.6: Simulation result for $M_2 = 2/8$ (a) and $M_2 = 6/8$ (b).

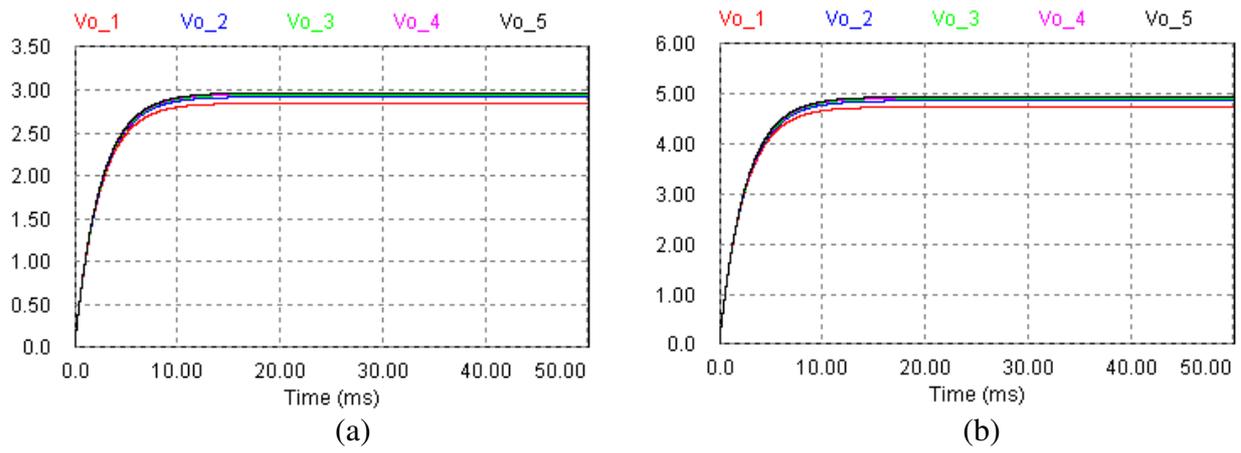

(a)                                                                       (b)

Figure 6.1.7: Simulation result for $M_3 = 3/8$ (a) and $M_3 = 5/8$ (b).

Considering the SCC equivalent circuit (Fig. 6.1.8), we can write the output voltage $V_o$ using the formula of voltage divider (6.1.1), where the target voltage $V_{TRG} = M_n \cdot V_{in}$.

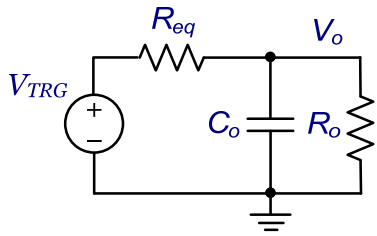

Figure 6.1.8: The SCC equivalent circuit.

$$V_o = \frac{R_o}{R_{eq} + R_o} \cdot V_{TRG} \qquad (6.1.1)$$

Rearranging (6.1.1) we obtain the equivalent resistor as:

$$R_{eq} = \left(\frac{V_{TRG}}{V_o} - 1\right) \cdot R_o \qquad (6.1.2)$$

For each ratio from $M_3 = 1/8$ through $M_3 = 7/8$, the values of $V_o$ measured at $t = 50$ms and the values of $R_{eq}$ calculated according to (6.1.2) are summarized in Table 6.1.1.

Table 6.1.1: Parameter sweep data for all the ratios $M_n$, $n = 1, \ldots, 3$.

| $R_o$, Ohm | $M_3 = 1/8$ | | $M_2 = 2/8$ | | $M_3 = 3/8$ | | $M_1 = 4/8$ | |
|---|---|---|---|---|---|---|---|---|
| | $V_o$,(V) | $R_{eq}$,($\Omega$) | $V_o$,(V) | $R_{eq}$,($\Omega$) | $V_o$,(V) | $R_{eq}$,($\Omega$) | $V_o$,(V) | $R_{eq}$,($\Omega$) |
| 100 | 0.938 | 6.61 | 1.897 | 5.43 | 2.846 | 5.411 | 3.816 | 4.822 |
| 200 | 0.968 | 6.612 | 1.947 | 5.444 | 2.921 | 5.409 | 3.906 | 4.813 |
| 300 | 0.978 | 6.748 | 1.964 | 5.499 | 2.947 | 5.395 | 3.937 | 4.801 |
| 400 | 0.984 | 6.504 | 1.973 | 5.474 | 2.959 | 5.542 | 3.952 | 4.858 |
| 500 | 0.987 | 6.586 | 1.978 | 5.561 | 2.968 | 5.391 | 3.962 | 4.796 |

Table 6.1.1: cont'd.

| $R_o$, Ohm | $M_3 = 5/8$ | | $M_2 = 6/8$ | | $M_3 = 7/8$ | |
|---|---|---|---|---|---|---|
| | $V_o$,(V) | $R_{eq}$,($\Omega$) | $V_o$,(V) | $R_{eq}$,($\Omega$) | $V_o$,(V) | $R_{eq}$,($\Omega$) |
| 100 | 4.743 | 5.419 | 5.692 | 5.411 | 6.565 | 6.626 |
| 200 | 4.868 | 5.423 | 5.841 | 5.444 | 6.776 | 6.612 |
| 300 | 4.911 | 5.437 | 5.894 | 5.395 | 6.849 | 6.614 |
| 400 | 4.933 | 5.433 | 5.919 | 5.474 | 6.886 | 6.622 |
| 500 | 4.946 | 5.459 | 5.936 | 5.391 | 6.908 | 6.659 |

The averaged values of $R_{eq}$ are presented in Table 6.1.2. These values are found to be in excellent agreement with their theoretical counterparts given in Table 5.2.5.

Table 6.1.2: The values of $R_{eq}$ the obtained by simulations and theoretically.

| $M_n / R_{eq}$ | $M_3 = 1/8$ | $M_2 = 2/8$ | $M_3 = 3/8$ | $M_1 = 4/8$ | $M_3 = 5/8$ | $M_2 = 6/8$ | $M_3 = 7/8$ |
|---|---|---|---|---|---|---|---|
| Averaged | 6.612 | 5.482 | 5.43 | 4.818 | 5.434 | 5.423 | 6.627 |
| Theoretical | 6.615 | 5.42 | 5.428 | 4.82 | 5.428 | 5.42 | 6.615 |

## 6.2 Utilizing the EXB based SCC in step-up mode

According to the concept demonstrated in Section 3.7.1, the step-down EXB based SCC can be utilized for a step-up conversion when its input and output are interchanged. To test the validity of this concept, we switch the input and output in the SCC circuit of Fig. 6.1.2 as shown in Fig 6.2.1.

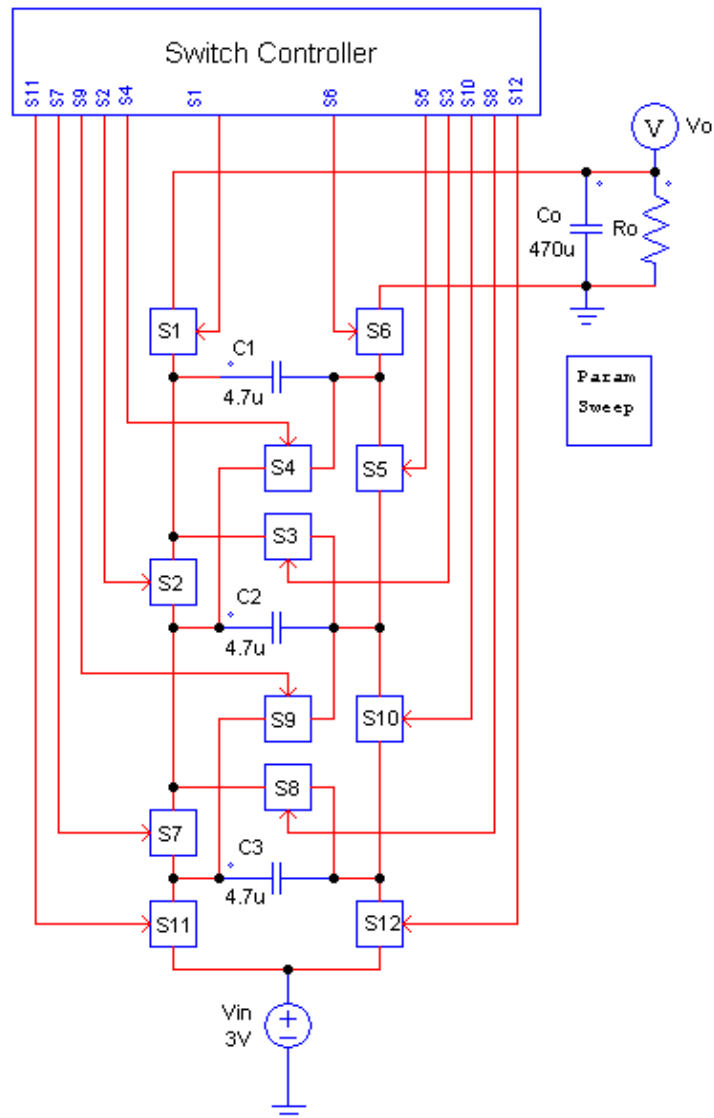

Figure 6.2.1: Simulation circuit for the step-up case.

As in the previous section we run the simulations for the EXB codes given in Table 5.2.3 and measure the output voltage $V_o$ for different load resistances $R_o$ = 1kΩ, 1.5kΩ, …, 3kΩ when the SCC has reached the steady-state.

For the EXB codes of $M_3 = 3/8$ and $M_3 = 5/8$, the measured voltages $V_o$ are depicted in Fig. 6.2.2. As expected, these voltages are the reciprocals of their step-down counterparts.

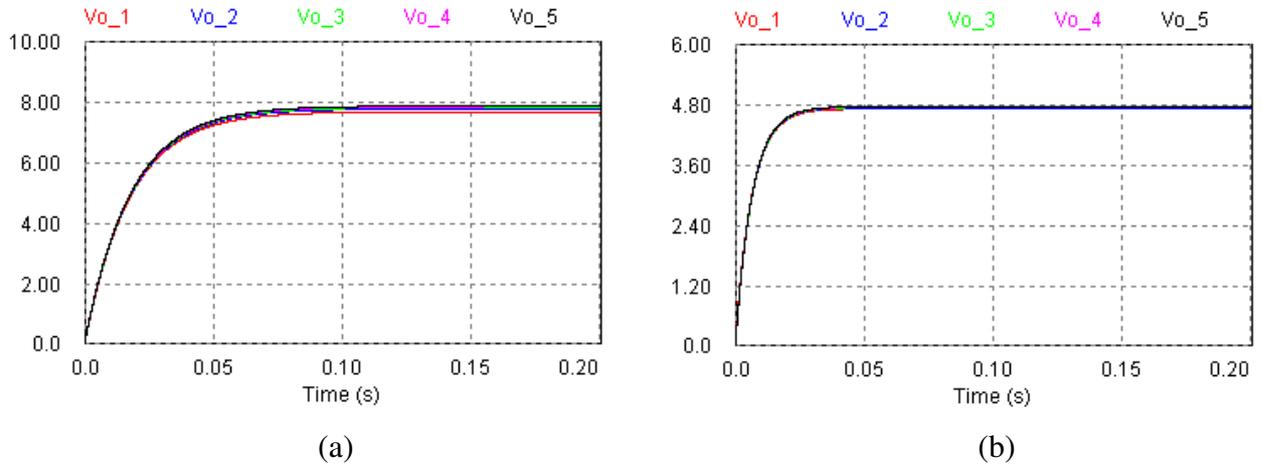

(a)  (b)

Figure 6.2.2: Simulation result for $1/M_3 = 8/3$ (a) and $1/M_3 = 8/5$ (b).

The values of $V_o$ measured at $t = 0.2$s and the values of $R_{eq}$ calculated according to (6.1.2) for both conversion ratios $1/M_3 = 8/3$ and $1/M_3 = 8/5$ are summarized in Table 6.2.1.

Table 6.2.1

| $R_o$, Ohm | $1/M_3 = 8/3$ | | $1/M_3 = 8/5$ | |
|---|---|---|---|---|
| | $V_o$,(V) | $R_{eq}$,($\Omega$) | $V_o$,(V) | $R_{eq}$,($\Omega$) |
| 1000 | 7.703 | 38.556 | 4.734 | 13.942 |
| 1500 | 7.799 | 38.659 | 4.756 | 13.877 |
| 2000 | 7.849 | 38.476 | 4.767 | 13.845 |
| 2500 | 7.878 | 38.715 | 4.773 | 14.142 |
| 3000 | 7.898 | 38.744 | 4.778 | 13.813 |

As seen from Table 6.2.1, the value of $R_{eq}$ for $1/M_3 = 8/3$ is greater then for $1/M_3 = 8/5$. This can be explained by considering the coefficients $A_0$ in the EXB sequences of $M_3 = 3/8$ and $M_3 = 5/8$ (Table 5.2.3). When $A_0 = 1$, the load is connected and therefore the EXB sequence where the number of $A_0 = 1$ is large provides low $R_{eq}$. The EXB sequence of $M_3 = 5/8$ comprises three $A_0 = 1$, while the sequence of $M_3 = 3/8$ comprises two, so the value of $R_{eq}$ for $1/M_3 = 8/3$ is expected to be greater.

## 6.3 Test for unipolar voltages across the switches

The switches used in the above simulations are bidirectional, however from a practical point of view, these switches are expensive. In this section we examine the feasibility of replacing some bidirectional switches with unidirectional devices, which comprise a diode between the terminals (for example a MOSFET). When this diode is forward biased the regular operation of the SCC can be disturbed. So, we need to check the polarity of voltages across all the switches. In the case of unipolar voltages, the corresponding switches can be unidirectional. The connection of voltage probes to the circuit of Fig. 6.1.2 is depicted in Fig. 6.3.1.

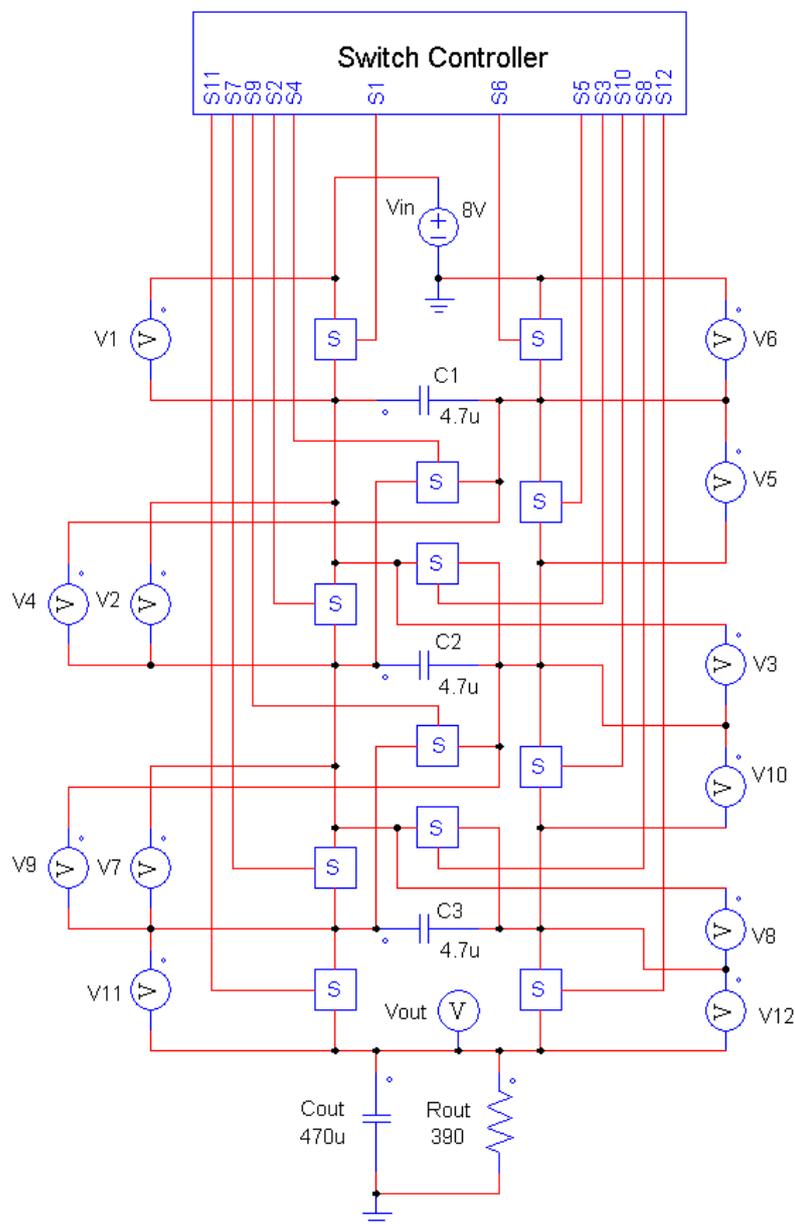

Figure 6.3.1: Measuring the voltages across the switches.

Fig. 6.3.2 depicts the measured voltages for the case of $M_1 = 4/8$. The measurements for the complementary ratios $M_n$ and $1 - M_n$, $n = 1\ldots3$, are shown side by side in Fig. 6.3.3 - 6.3.5.

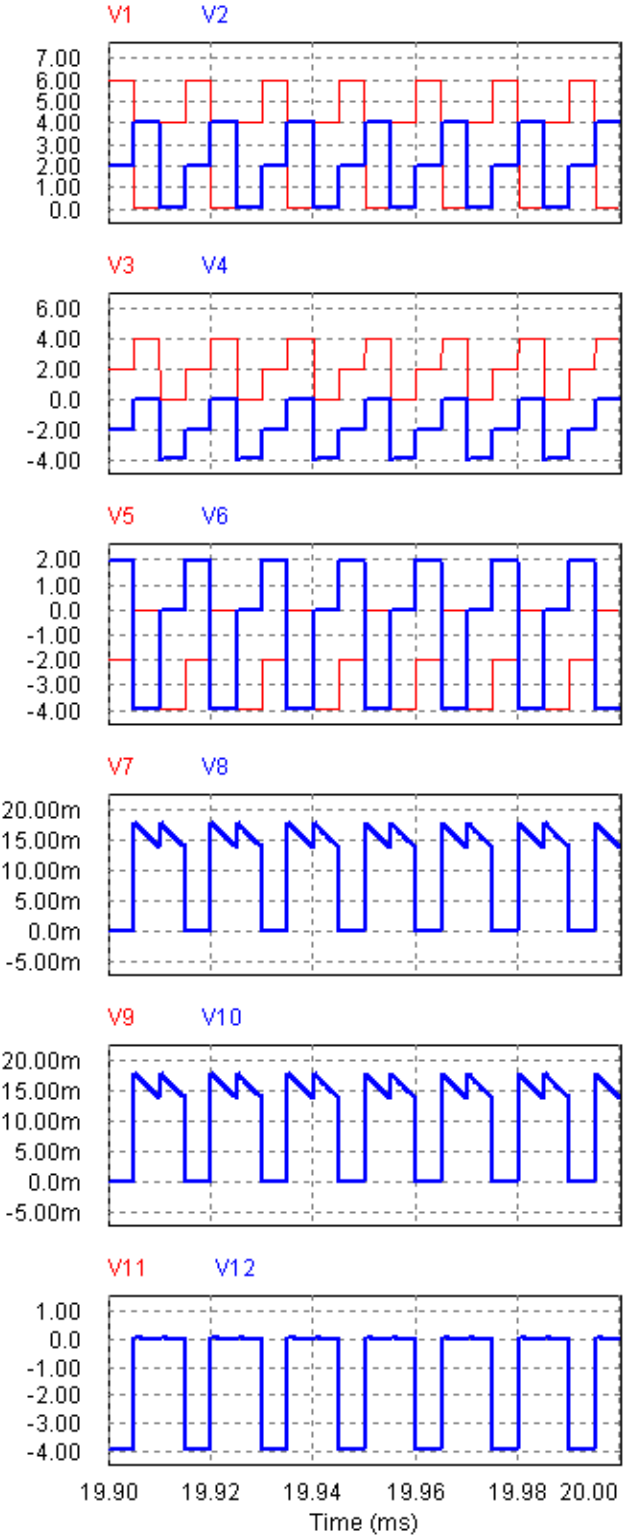

Figure 6.3.2: Measured voltages for $M_1 = 4/8$.

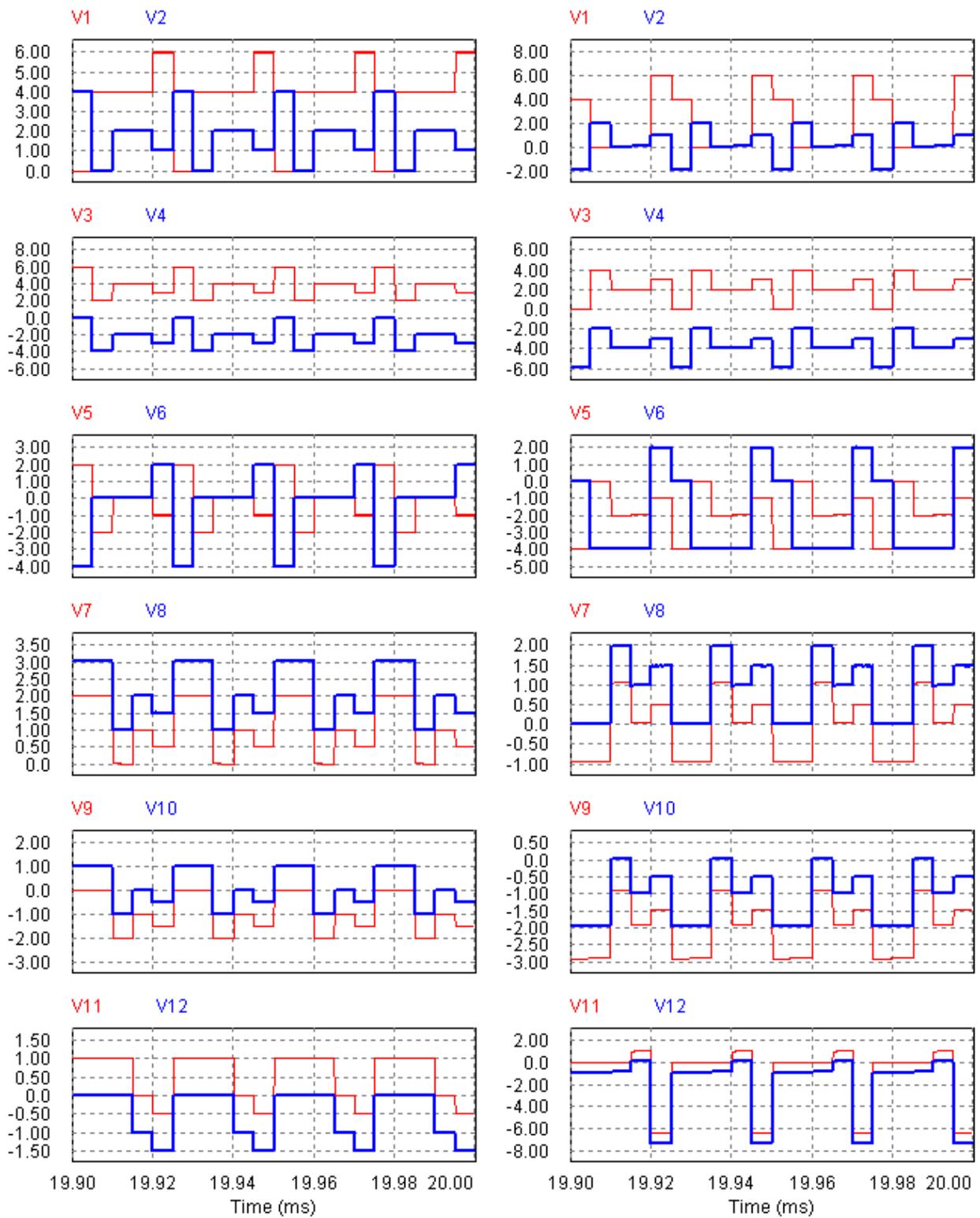

Figure 6.3.3: Measured voltages for $M_3 = 1/8$ (a) and $M_3 = 7/8$ (b).

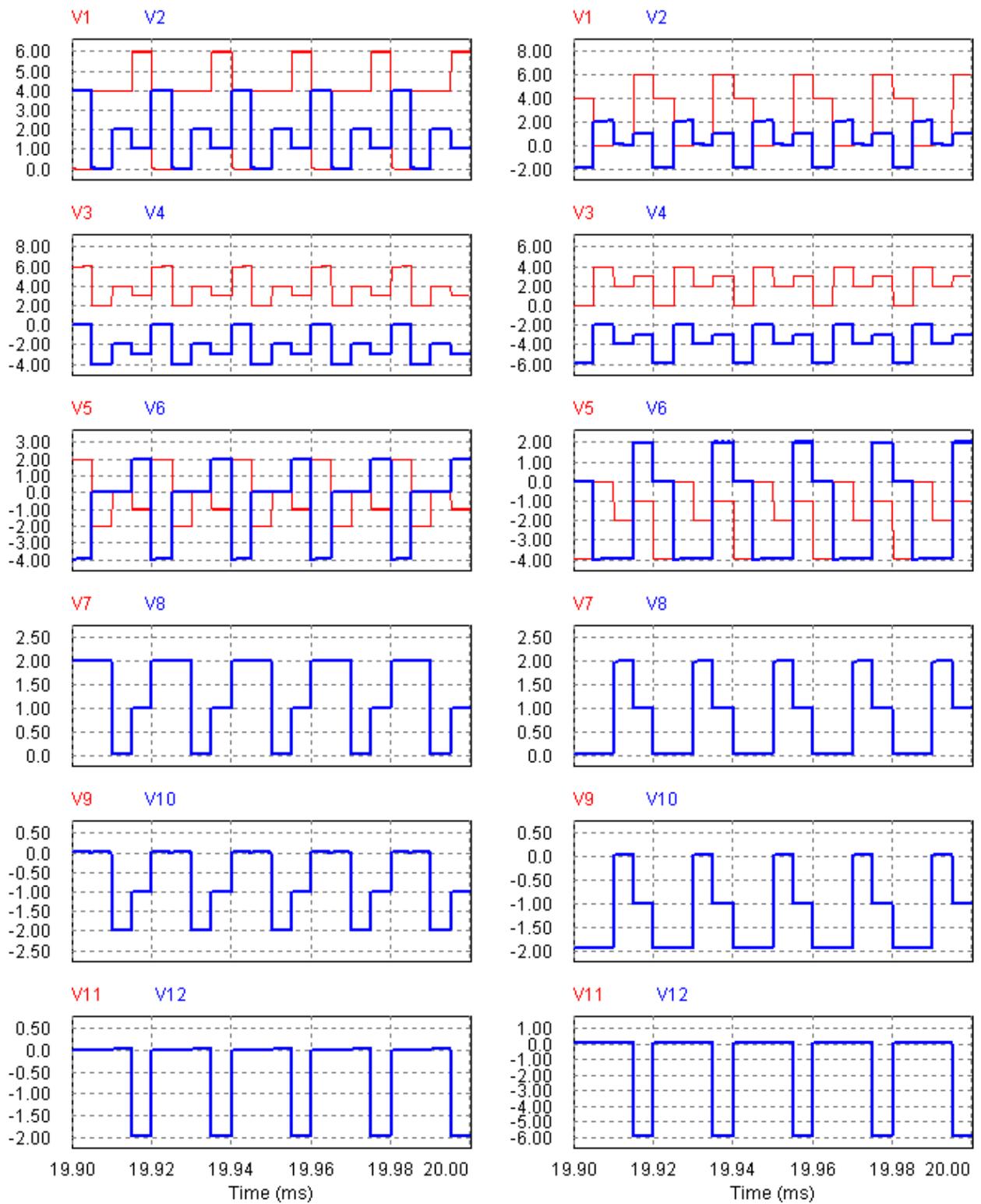

Figure 6.3.4: Measured voltages for $M_2 = 2/8$ (a) and $M_2 = 6/8$ (b).

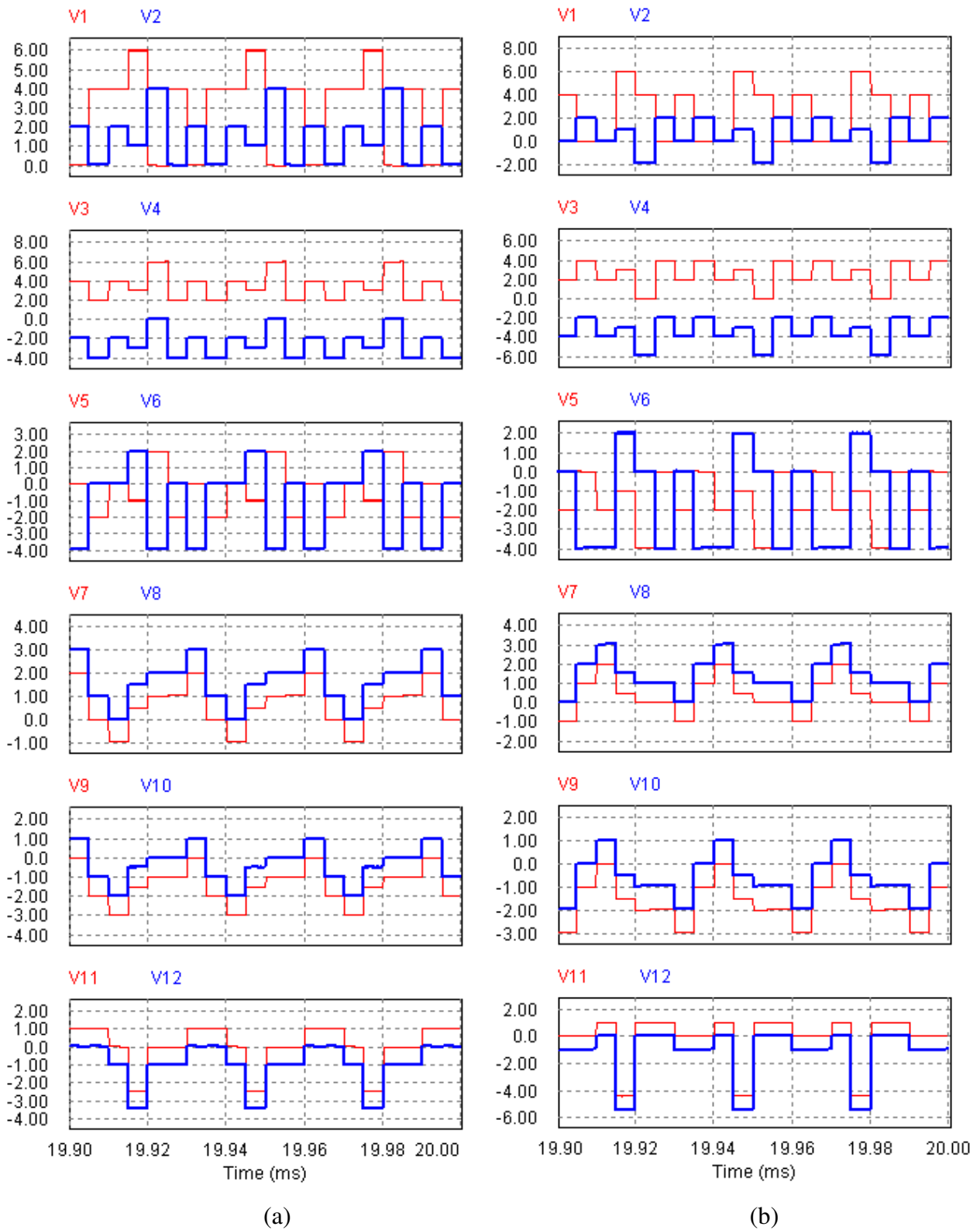

Figure 6.3.5: Measured voltages for $M_3 = 3/8$ (a) and $M_3 = 5/8$ (b).

Inspecting Fig. 6.3.2 - Fig. 6.3.5 for unipolar voltages, we map the switches used in the circuit of Fig. 6.1.2 as presented in Table 6.3.1, where "u" and "b" designate unidirectional and bidirectional switch respectively.

Table 6.3.1: Mapping the switches used.

| $M_n/S$ | 1/8 | 2/8 | 3/8 | 4/8 | 5/8 | 6/8 | 7/8 |
|---|---|---|---|---|---|---|---|
| S1 | u | u | u | u | u | u | u |
| S2 | u | u | u | u | b | b | b |
| S3 | u | u | u | u | u | u | u |
| S4 | u | u | u | u | u | u | u |
| S5 | b | b | b | u | u | u | u |
| S6 | b | b | b | b | b | b | b |
| S7 | u | u | b | u | b | u | b |
| S8 | u | u | u | u | u | u | u |
| S9 | u | u | u | u | u | u | u |
| S10 | b | u | b | u | b | u | u |
| S11 | b | u | b | u | b | u | b |
| S12 | u | u | u | u | u | u | u |

# 7. EXPERIMENTAL RESULTS

*"We judge ourselves by what we feel capable of doing, while others judge us by what we have already done."*

Henry W. Longfellow

## 7.1 Response to a step in input voltage

The measurements depicted in Fig. 7.1.1 and Fig. 7.1.2 are done for the conversion ratios $M_3 = 3/8$ and $M_3 = 5/8$, respectively, while $V_{in} = 8V$ and the load resistance $R_o = 3.6k\Omega$.

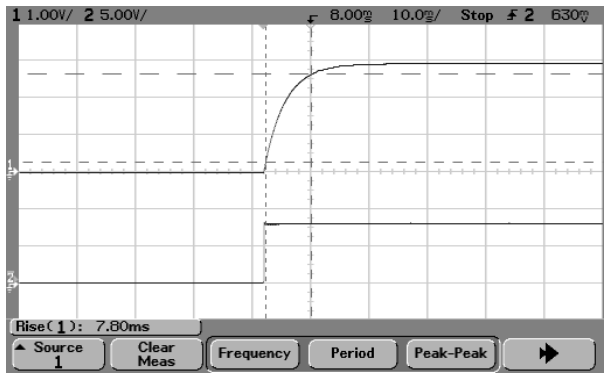
(a)

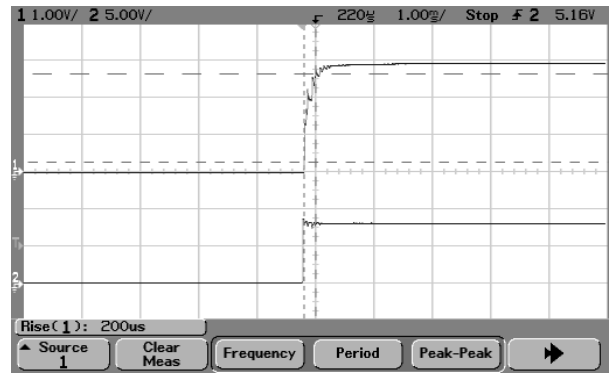
(b)

Figure 7.1.1: The SCC cold start, $M_3 = 3/8$, $C_o = 470\mu F$ (a) and $C_o = 22\mu F$ (b).

Vertical scale: 1V/div; Horizontal scale: 10ms/div.

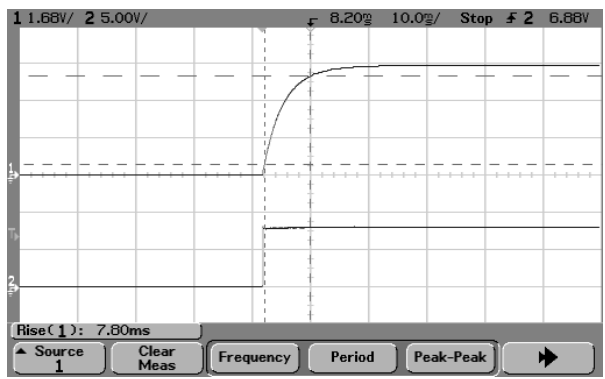
a)

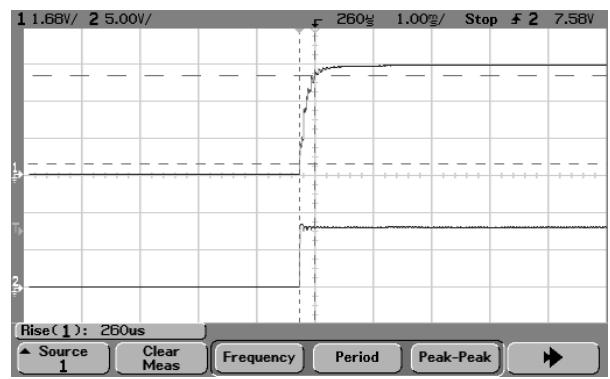
(b)

Figure 7.1.2: The SCC cold start, $M_3 = 5/8$, $C_o = 470\mu F$ (a) and $C_o = 22\mu F$ (b).

Vertical scale: 1.68V/div; Horizontal scale: 10ms/div.

## 7.2 Response to a step in load resistance

The load is switched from the nominal resistor 3.6kΩ to the indicated value. The top trace is the output voltage, while the bottom trace is the control signal.

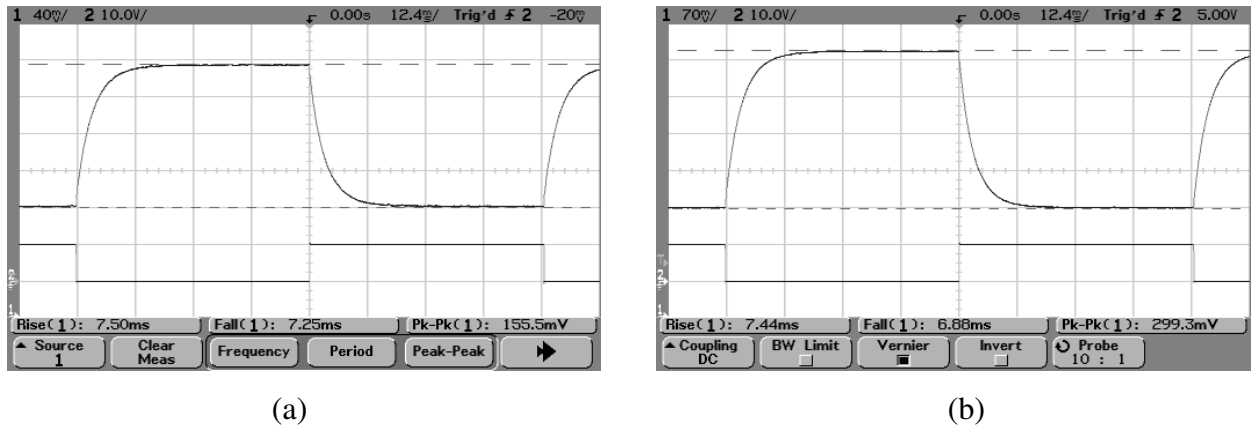

(a)                 (b)

Figure 7.2.1: The SCC response, $M_3 = 3/8$, $C_o = 470\mu F$, $R_o = 128\Omega$ (a) and $R_o = 62\Omega$ (b).

Vertical scale: 40mV/div (a) and 70mV/div (b); Horizontal scale: 12.4ms/div.

Peak to peak output voltage is 155.5mV (a) and 299.3mV (b).

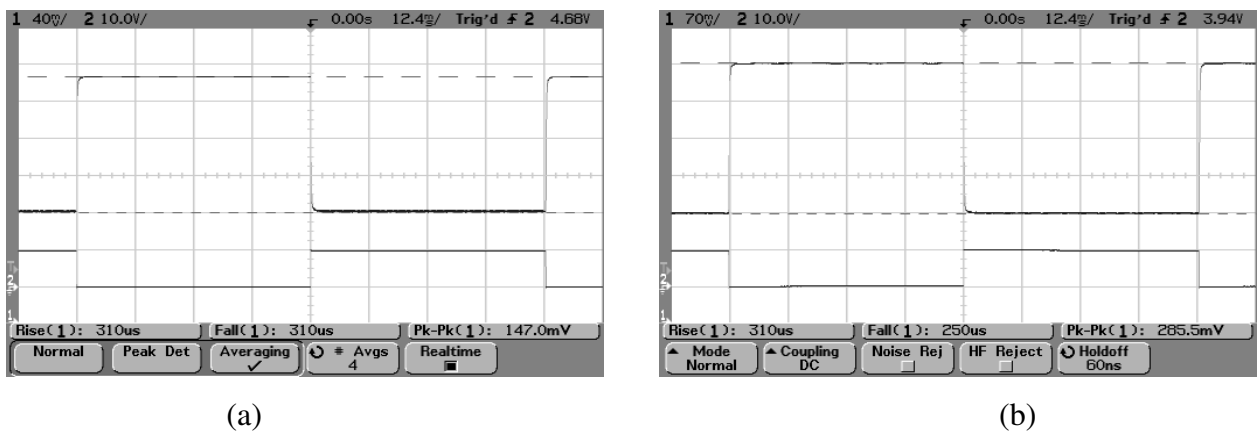

(a)                 (b)

Figure 7.2.2: The SCC response, $M_3 = 3/8$, $C_o = 22\mu F$, $R_o = 128\Omega$ (a) and $R_o = 62\Omega$ (b).

Vertical scale: 40mV/div (a) and 70mV/div (b); Horizontal scale: 12.4ms/div.

Peak to peak output voltage is 155.5mV (a) and 299.3mV (b).

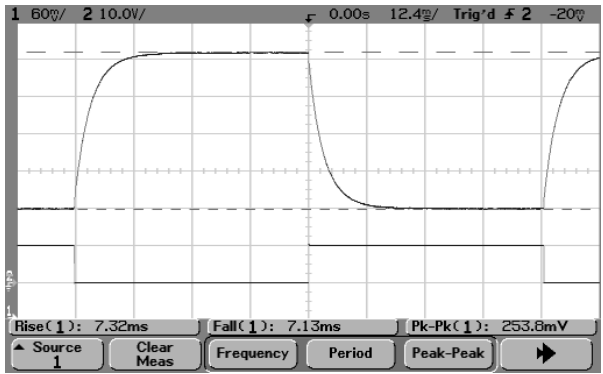 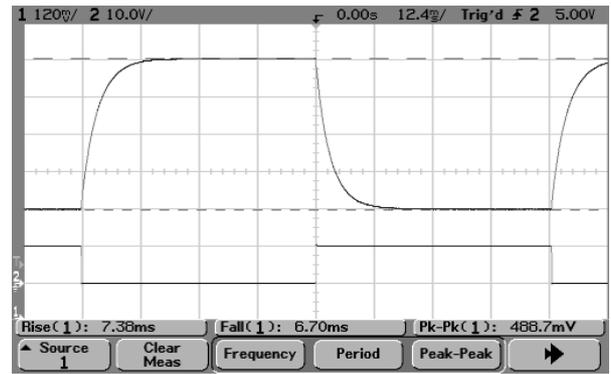

(a)                                                    (b)

Figure 7.2.3: The SCC response, $M_3 = 5/8$, $C_o = 470\mu F$, $R_o = 128\Omega$ (a) and $R_o = 62\Omega$ (b).

Vertical scale: 60mV/div (a) and 120mV/div (b); Horizontal scale: 12.4ms/div.

Peak to peak output voltage is 253.8mV (a) and 488.7mV (b).

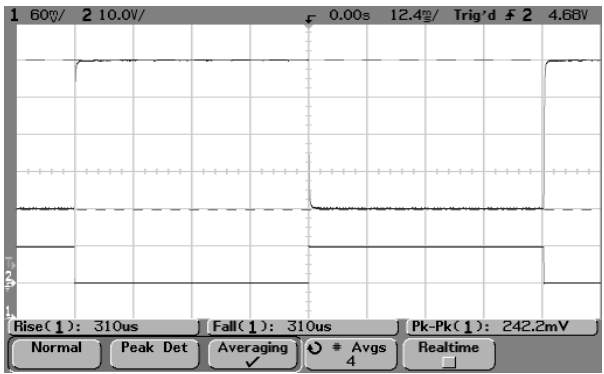 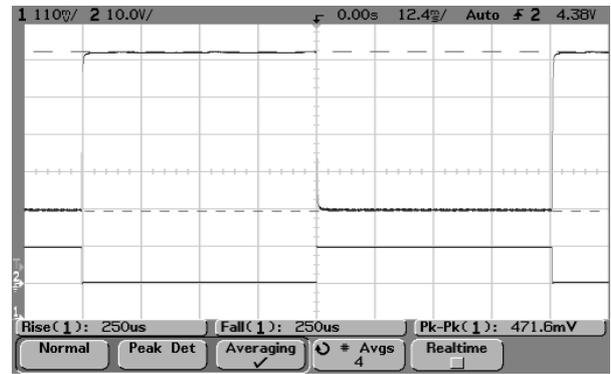

(a)                                                    (b)

Figure 7.2.4: The SCC response, $M_3 = 5/8$, $C_o = 22\mu F$, $R_o = 128\Omega$ (a) and $R_o = 62\Omega$ (b).

Vertical scale: 60mV/div (a) and 110mV/div (b); Horizontal scale: 12.4ms/div.

Peak to peak output voltage is 242.2mV (a) and 471.6mV (b).

## 7.3 Efficiency versus load resistance

The efficiencies measured for the step-down EXB based SCC are presented in Fig. 7.3.1 and Fig. 7.3.2, while Fig. 7.3.1 depicts the efficiency for the step-up case as given in [55].

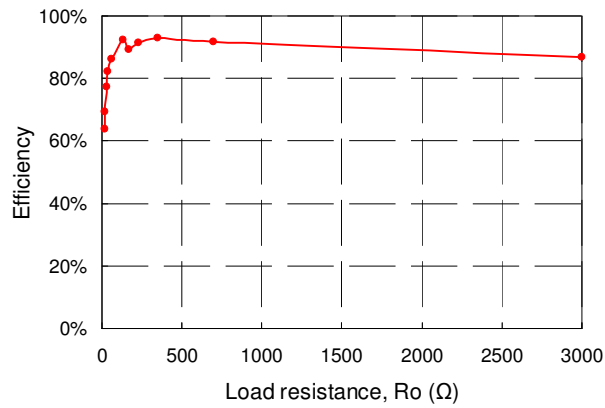

Figure 7.3.1: Efficiency of the step-down SCC with $M_3 = 3/8$.

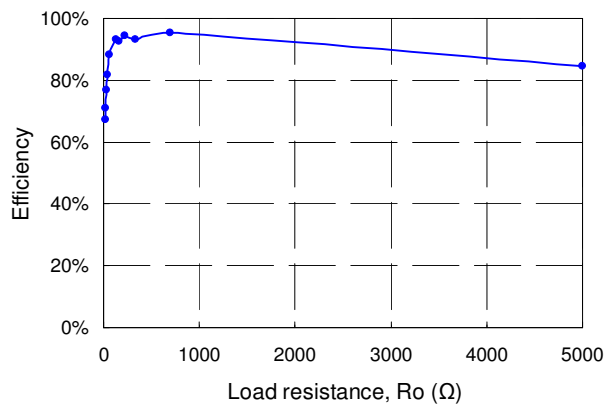

Figure 7.3.2: Efficiency of the step-down SCC with $M_3 = 5/8$.

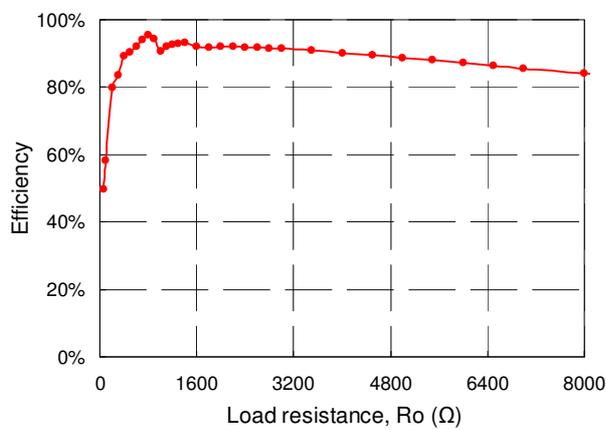

Figure 7.3.3: Efficiency of the step-up SCC with $1/M_3 = 8/3$.

## 7.4 Load characteristics and effect of $R_{eq}$

To test the validity of the equivalent circuit (Fig. 7.4.1) over a wide operational range, we consider the formula for the voltage divider (7.4.1), where the target voltage $V_{TRG} = M_n \cdot V_{in}$.

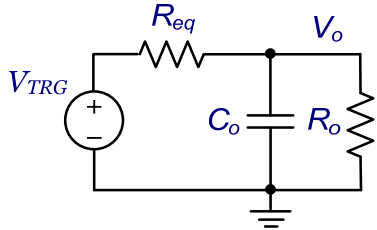

Figure 7.4.1: The SCC equivalent circuit.

$$I_o = \frac{V_{TRG}}{R_{eq} + R_o} = \frac{V_o}{R_o} \qquad (7.4.1)$$

Defining $y$ as:
$$y = \frac{R_o}{V_o} = \frac{R_{eq}}{V_{TRG}} + \frac{R_o}{V_{TRG}} \qquad (7.4.2)$$

One can rewrite (7.4.2) as:
$$y(x) = ax + b \qquad (7.4.3)$$

where $x = R_o$, $a = 1/V_{TRG}$ and $b = R_{eq}/V_{TRG}$. Hence, when plotting $R_o/V_o$ versus $R_o$, one should get a straight line with a slope of $1/V_{TRG}$ (Fig. 7.4.2) and intersection of the X axis ($y = 0$) at $R_{eq}/V_{TRG}$, from which $R_{eq}$ can be calculated.

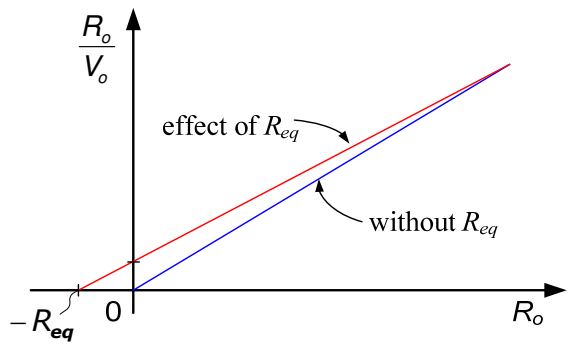

Figure 7.4.2: Intersection yielding the equivalent resistor.

The output resistance was tested for $V_{in} = 8V$ and $f_s = 100\text{kHz}$. Each SCC topology is configured for the time interval $t = T_s/(n + 1)$, where $T_s = 10\mu s$, and $n < 3$ is the resolution.

For *n* = 3 and the conversion ratios $M_3 = 3/8$ and $M_3 = 5/8$ all five SCC topologies are configured (one topology is redundant) and $t = T_s/5$.

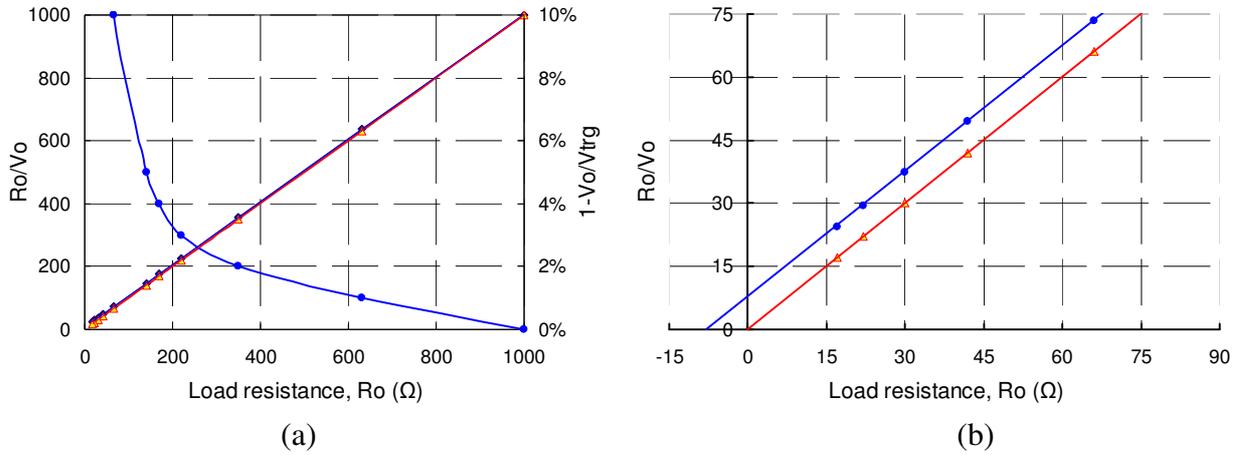

(a)                  (b)

Figure 7.4.3: Experimental result for $M_3 = 1/8$ (a) and close-up view (b).

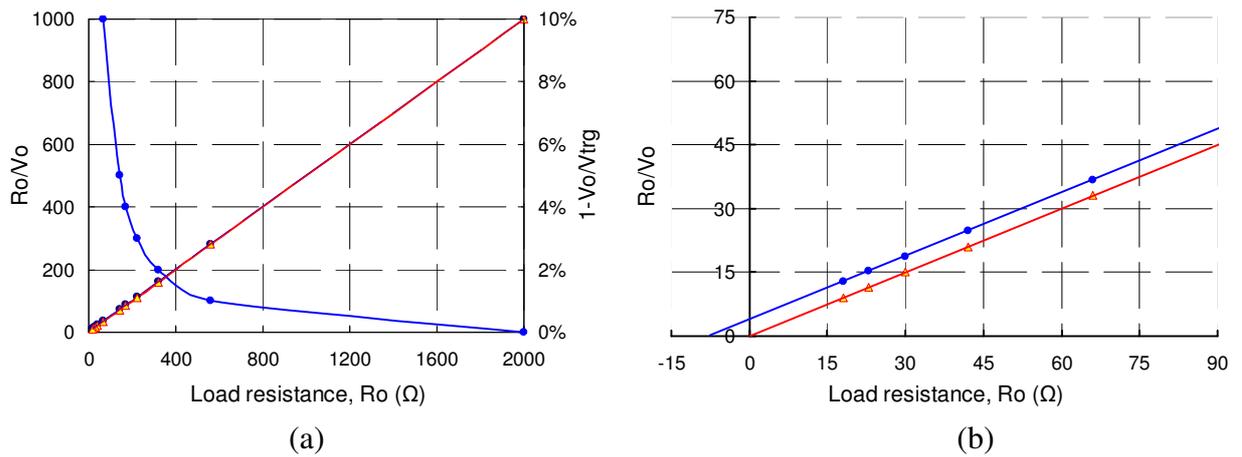

(a)                  (b)

Figure 7.4.4: Experimental result for $M_2 = 2/8$ (a) and close-up view (b).

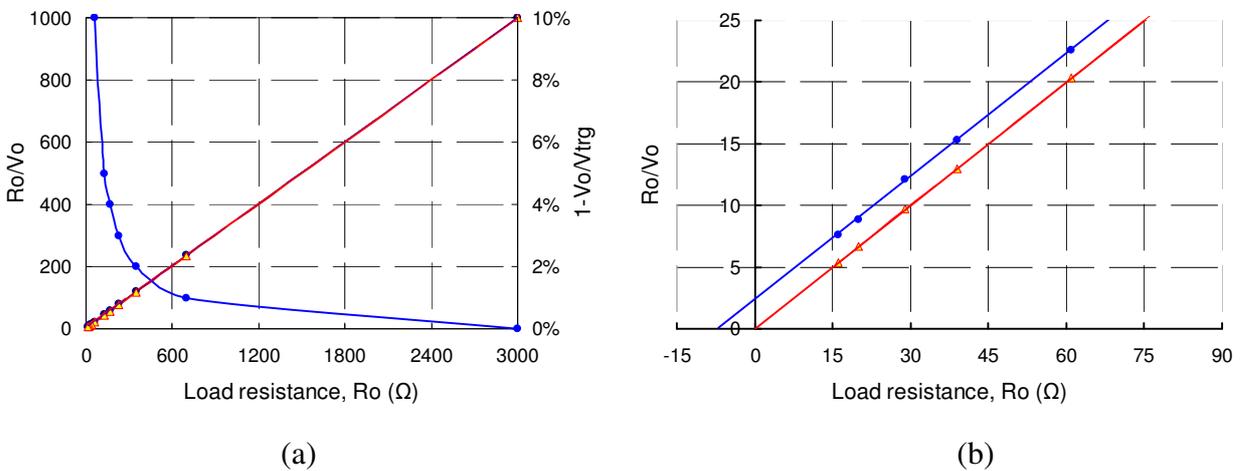

(a)                  (b)

Figure 7.4.5: Experimental result for $M_3 = 3/8$ (a) and close-up view (b).

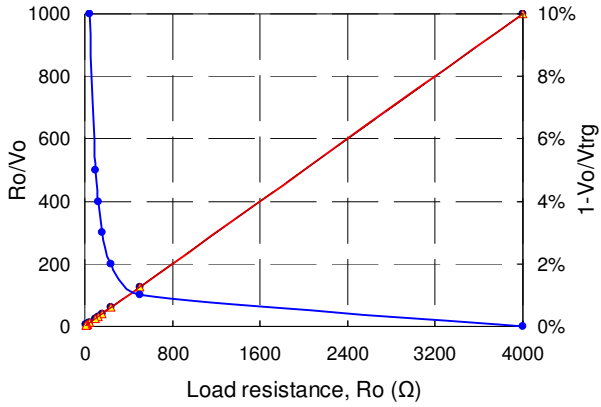
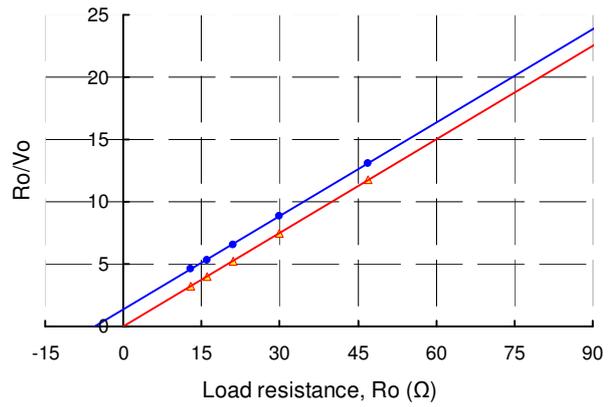

(a)                              (b)

Figure 7.4.6: Experimental result for $M_1 = 4/8$ (a) and close-up view (b).

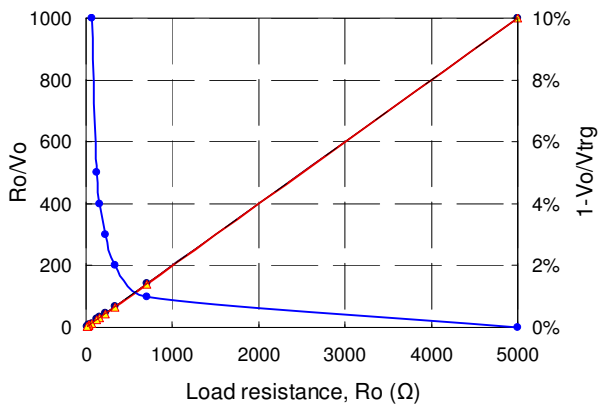
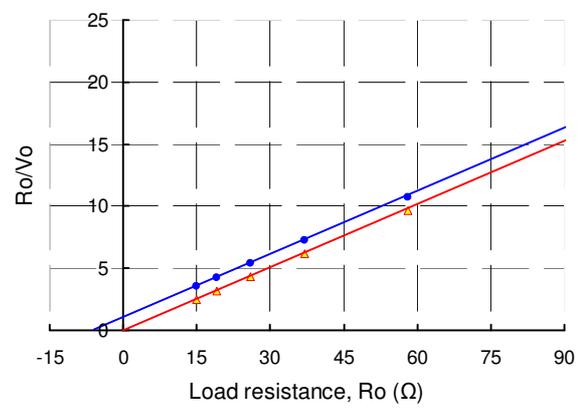

(a)                              (b)

Figure 7.4.7: Experimental result for $M_3 = 5/8$ (a) and close-up view (b).

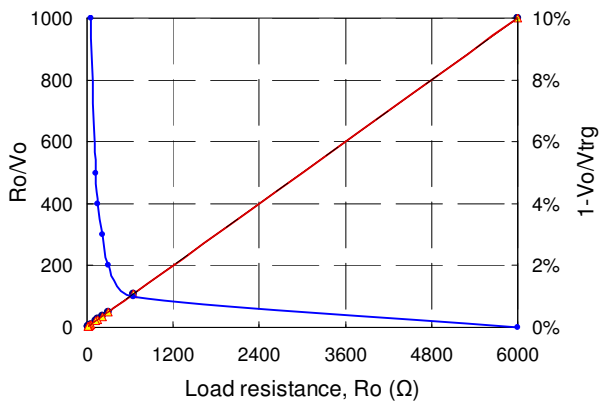
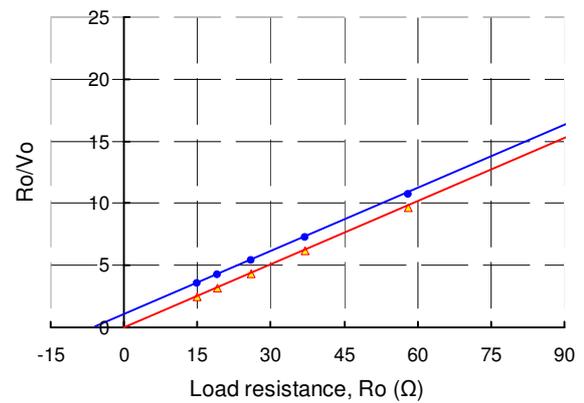

(a)                              (b)

Figure 7.4.8: Experimental result for $M_2 = 6/8$ (a) and close-up view (b).

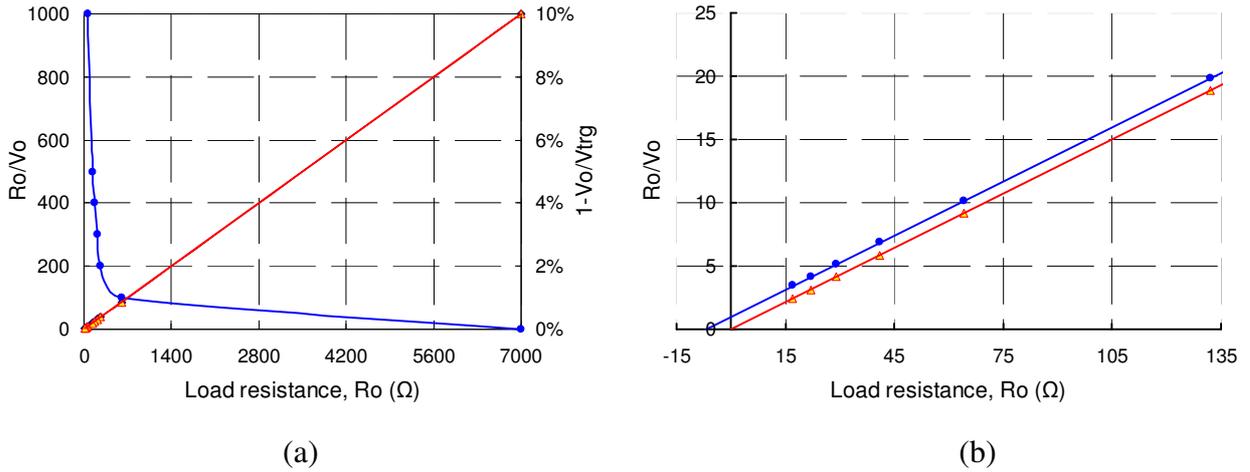

Figure 7.4.9: Experimental result for $M_1 = 7/8$ (a) and close-up view (b).

The experimentally obtained values of $R_{eq}$ are presented in Table 7.4.1. Although the theoretical and simulated values (Table 5.2.5 and Table 6.1.2) close, the experimental values are somewhat apart. This is probably due to the effects of parasitic elements that have not been considered in the theoretical model.

Table 7.4.1: The values of $R_{eq}$ obtained experimentally, by simulations and by theory.

| $M_n / R_{eq}(\Omega)$ | $M_3 = 1/8$ | $M_2 = 2/8$ | $M_3 = 3/8$ | $M_1 = 4/8$ | $M_3 = 5/8$ | $M_2 = 6/8$ | $M_3 = 7/8$ |
|---|---|---|---|---|---|---|---|
| Experimental | 7.65 | 6.82 | 7.35 | 5.3 | 7.16 | 6.53 | 6.94 |
| Simulated | 6.612 | 5.482 | 5.43 | 4.818 | 5.434 | 5.423 | 6.627 |
| Theoretical | 6.615 | 5.42 | 5.428 | 4.82 | 5.428 | 5.42 | 6.615 |

## 7.5 Output voltage regulation

As mentioned above, output voltage regulation of a SCC can be achieved by changing the equivalent resistor $R_{eq}$. We propose here two alternative approaches that are compatible with the structures of the EXB and GFN based SCC. In these SCC, it would be desirable to keep the voltages across the capacitors at their nominal values and not change them by partial charges or discharges. One approach to accomplish this is to use dithering that is, switching from one transfer ratio to another. In the regular one ratio mode, the SCC will scan over all the codes that correspond to the desired conversion ratio $M_n$. For conversion ratios which are in between the discrete $M_n$ values one can dither between two neighboring ratios as depicted in Fig. 7.5.1.

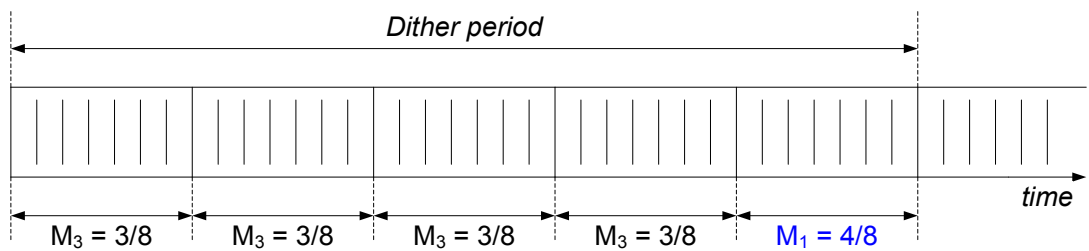

Figure 7.5.1: Dithering between $M_3 = 3/8$ and $M_1 = 4/8$ (in 4:1 ratio).

In the case depicted in Fig. 7.5.1 the dither duration is 5 sequences, 4 of 3/8 and one for 4/8, consequently the average ratio will be:

$$\frac{V_o}{V_{in}} = \frac{4}{5} \cdot \frac{3}{8} + \frac{1}{5} \cdot \frac{4}{8} = \frac{2}{5} = 0.4 \qquad (7.5.1)$$

For $V_{in} = 8V$ and $R_o = 437\Omega$ the output voltage is depicted in Fig 7.5.2, where the vertical scale is 10mV/div, while the horizontal scale is 100μS/div.

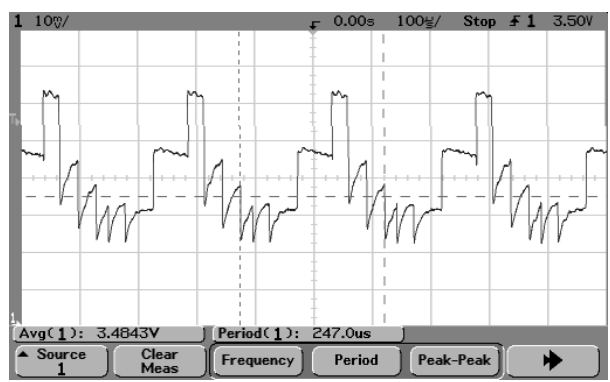

Figure 7.5.2: Output ripple. Dithering between $M_3 = 3/8$ and $M_1 = 1/2$.

Another method proposed here for output voltage control is to introduce a linear, low dropout (LDO) voltage regulator at the output (Fig. 7.5.3). In this case, the LDO will provide the regulation for the LSB while the SCC maintains a low voltage across the LDO.

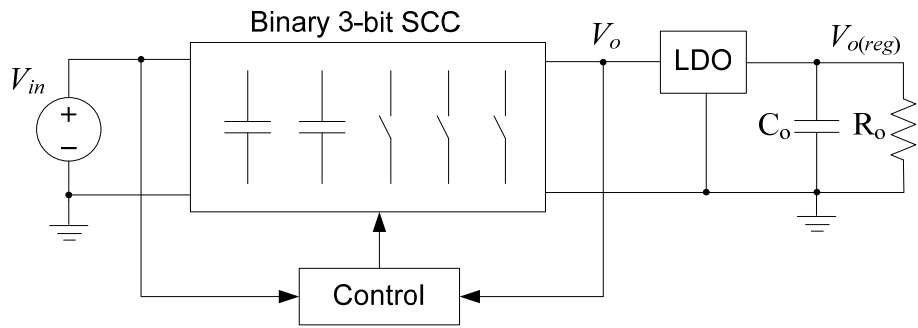

Figure 7.5.3: Block diagram of output voltage regulation by LDO at the output.

The concept of regulation with a LDO at the output was tested [55] using a LT1121-3.3 (Linear Technology) with a fixed output voltage of 3.3V. As shown in Fig. 7.5.3, the output and input voltages were sampled and the control was programmed to select the minimal conversion ratio which provides an output voltage greater than 3.6V. The LT1121 has a minimum dropout voltage of 0.3V, so at least 3.6V is required at the input. This implies that the upper limit of the efficiency is 3.3/3.6 = 0.92 and that is without taking into account the losses of the SCC.

In this preliminary, proof of concept, closed loop experiment, the EXB based SCC was configured to operate in both step-up and step-down modes by introducing extra switches that could flip the input and output terminals. The measured efficiency [55] for the input voltage range of 1.8V to 10V is depicted in Fig. 7.5.4.

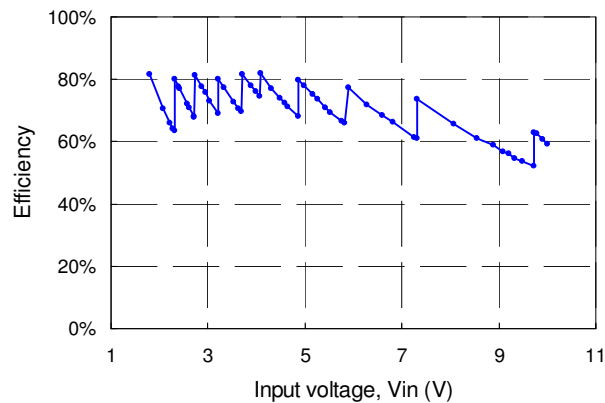

Figure 7.5.4: Efficiency of EXB based SCC operating with an LDO.

# 8. Discussion and Conclusions

*"Never express yourself more clearly than you are able to think."*
Niels Bohr

An EXB fractional representation is proposed and extended to the general radix case defined as GFN. In the case of the EXB, the radix is 2, while the general GFN can be defined for any radix $r$. Hence, the particular GFN case of $r = 2$ and the corresponding fractions $N_n(2)$ is in fact the EXB case. Based on these new fractional representations, a novel procedure is proposed for the design of high resolution multi-target SCC that emulate the EXB and GFN codes. It is shown that these SCC can be considered as hardware computational systems that solve a set of equations determined by the EXB, or in the general case, by the GFN representations. It is further shown that, for a given number of capacitors, one can generate many target voltages by configuring the flying capacitors interconnections according to different EXB or GFN codes.

These codes are used to derive the equivalent resistor, which defines both the output voltage drop and the power loss due to conversion inefficiencies. The experimentally obtained values of the equivalent resistor were found to be in good agreement with both the theoretically predicted values of $R_{eq}$ and those obtained in simulations. The new theoretically supported concepts were verified by simulation and experiment for static and dynamic responses. The experiments were conducted on the step-down and step-up EXB based SCC. Several control schemes were tested, including linear and dithering approaches to provide continuous regulation of the output voltage. Both of the proposed control approaches were found to function properly, however the dithering scheme gave rise to a higher output ripple. This could be explained by the fact that in this control method the SCC has to be reconfigured dynamically between two $M_n$ values.

Notwithstanding the good results, a number of theoretical issues are still open and require further investigation. Although some realizations of the proposed SCC have been described by way of illustration, they can be put into practice with many modifications that are within the scope of application engineers.

*Prototypes of binary SCC*:

# APPENDIX A. Circuit diagrams

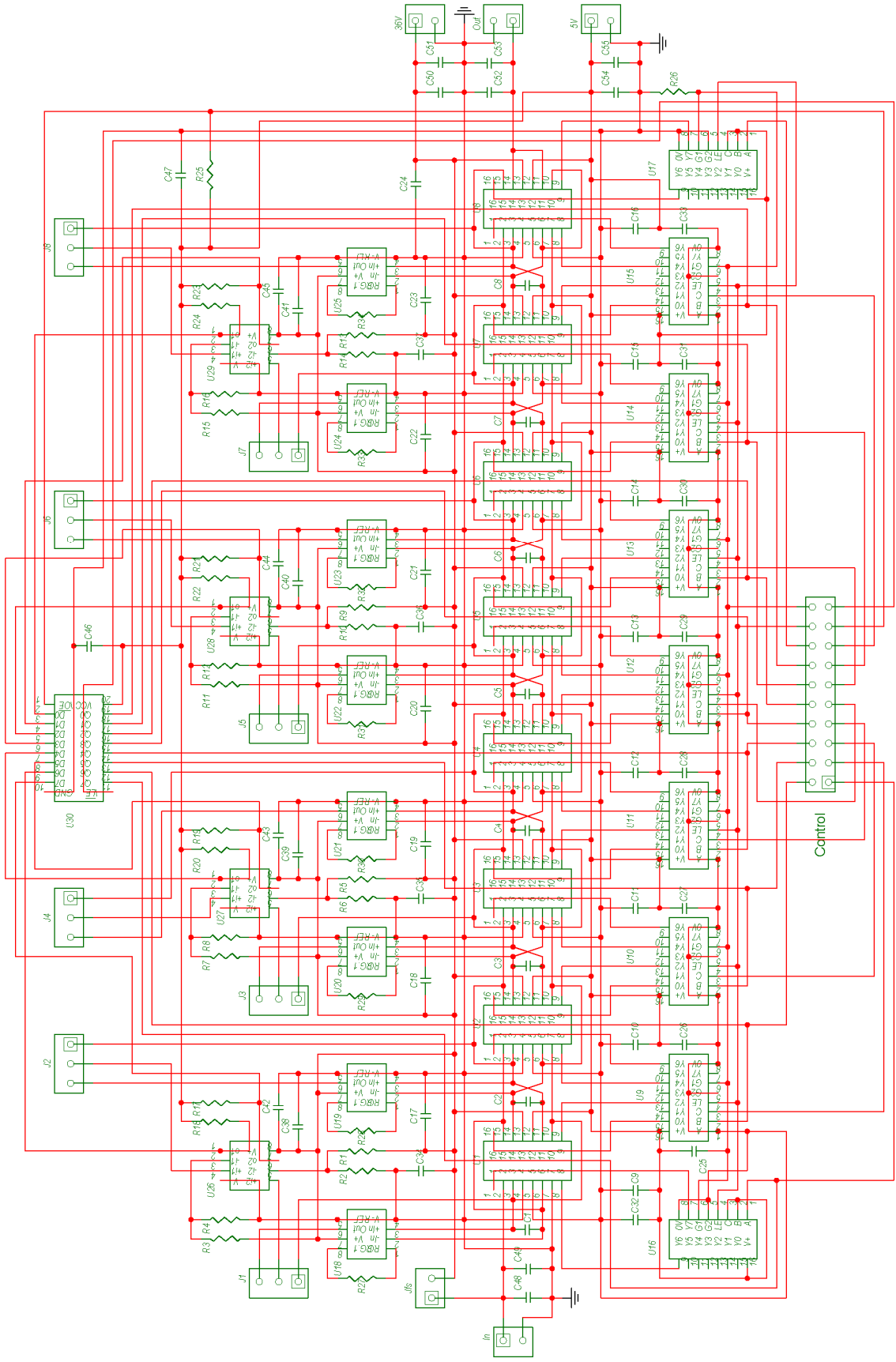

Circuit diagram of the EXB based SCC. Part I: Power stage 8-bit with measuring elements.

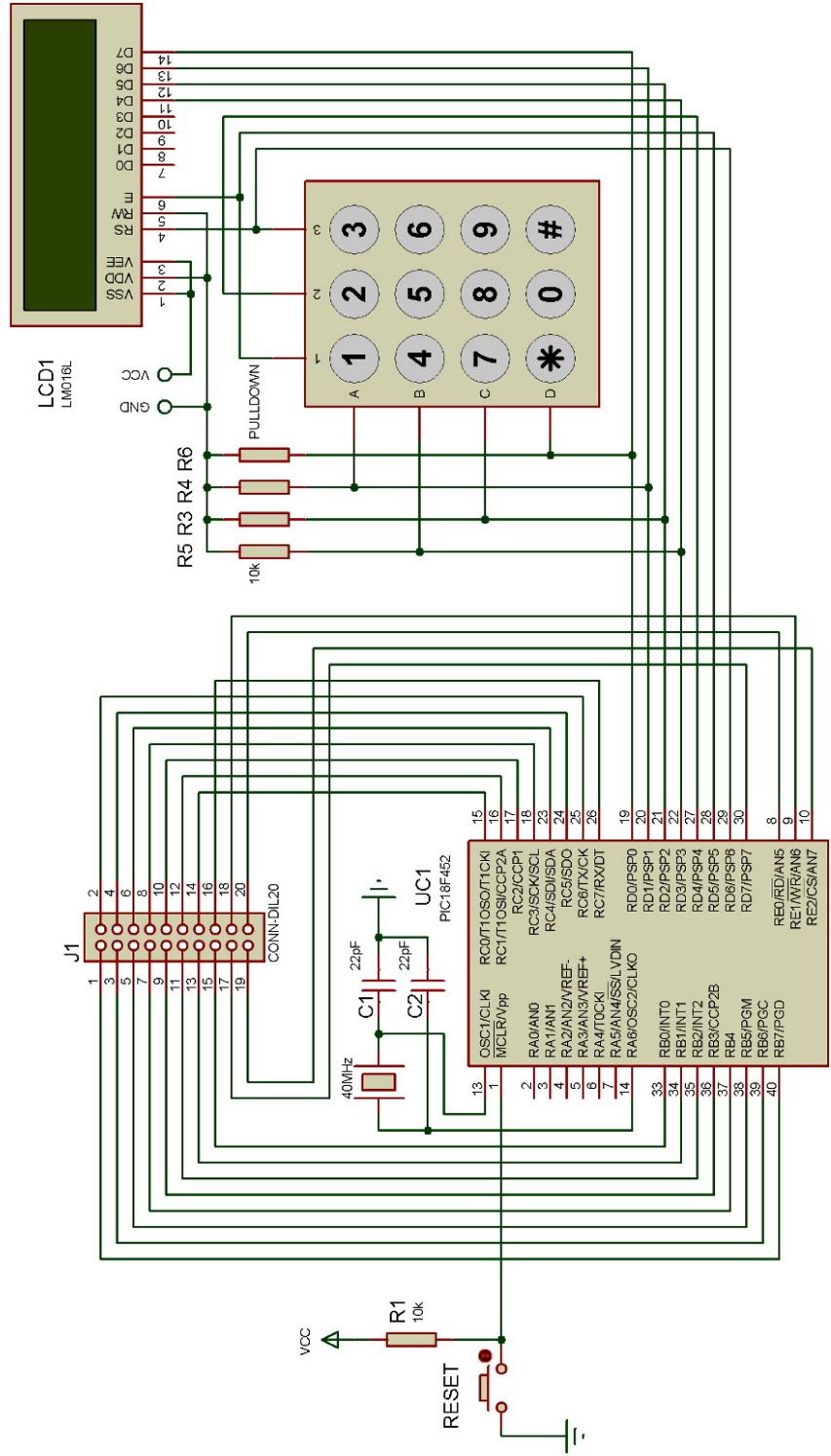

Circuit diagram of the EXB based SCC. Part II: Microcontroller connection to 8-bit power stage.

Microcontroller PIC18F452 providing $f_s$ = 100kHz.
MAX4678 quad SPST normally-off CMOS switch
74HC237 decoder 3-to-8 lines with address latches
4.7μF×50V MLCC marked as C340C475M5U5TA
470μF×50V electrolytic capacitor
0.1μF ceramic bypass capacitors

Unused components:

74HC573 octal 3-state transparent latch
INA118 instrumental amplifier
LM393 comparators
Precision resistors

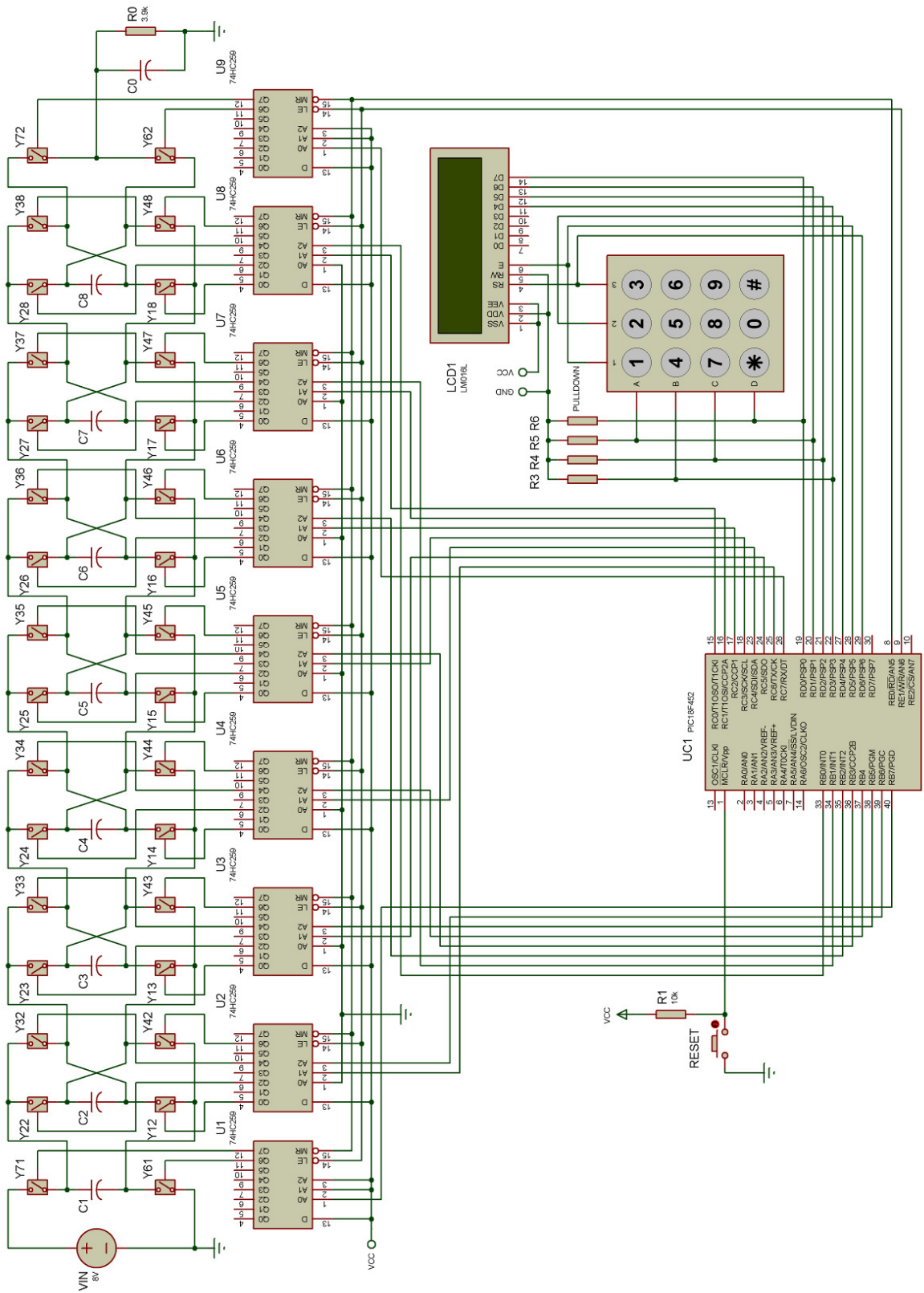

Alternative schematic of the EXB based SCC for simulation in the ISIS Proteus 7 software

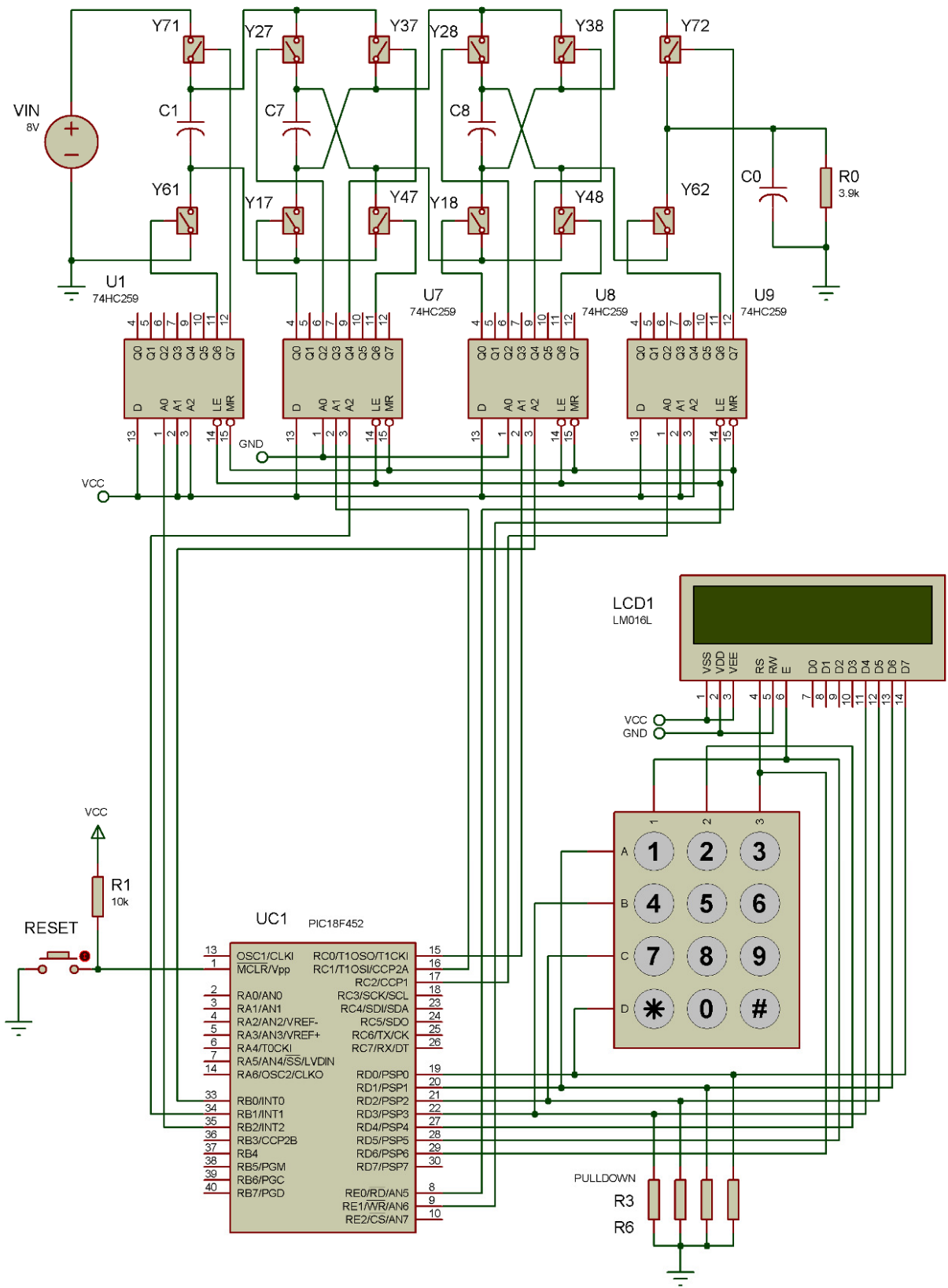

Example how to reduce the circuit to 3-bit by eliminating the intermediate stages

# APPENDIX B. Program listings

```
%*********************************************************
%              Spawning the EXB codes
%*********************************************************
clc
clear all
M=fullfact([3 3 3])-2;
K=[4 2 1];
F=M*K';
F1=F;
h=numel(F);
B=zeros(h,1);
i=find(F<0);
msb=sum(K)+1;
F1(i,:)=F1(i,:)+msb;
B(i,:)=1;
for g=1:msb-1
i=find(F1==g);
exb=[B(i) M(i,:)]
end

%*********************************************************
%              Spawning the GFN codes
%*********************************************************
clc
clear all
M=fullfact([5 5])-3;
K=[3 1];
F=M*K';
F1=F;
h=numel(F);
B=zeros(h,1);
i=find(F<0);
msb=2*sum(K)+1;
F1(i,:)=F1(i,:)+msb;
B(i,:)=1;
for g=1:msb-1
i=find(F1==g);
exb=[B(i) M(i,:)]
end

%*********************************************************
%       Investigating the voltage convergence issue
%*********************************************************
clc;
clear all;
close all;
% input EXB codes of M3=3/8
exb=[
    1    -1    -1     1
    0     1    -1     1
    1    -1     0    -1
    0     1     0    -1
    0     0     1     1];
```

```matlab
% forming the main diagonal
% due to the resolution of M3
n=numel(exb(1,:))+1;
E=eye(n);
% grouping A0 at the right
A1=exb(:,2:4);
A2=-exb(:,1);
A=[A1 A2];
% values of the capacitors
uF=10^(-6);
Co=470*uF;
flycap=4.7*uF.*ones(3,1);
caps=[flycap; Co];
Vin=8; Vo=0;
% initial voltages
pres_Vc=[0; 0; 0; Vo; 0];
% convergence double loop
old=0;
Q=[0, 0];
Vall=pres_Vc';
for j=1:200
for i=1:5
% adding the EXB row
E(n,:)=[A(i,:),0];
% coefficient -1 for Vo
E(n,n-1)=-1;
% adding the EXB collumn
A(i,1:n-1);
E(:,n)=E(n,:);
% dividing by the capacitances
E(1:n-1,n)=E(1:n-1,n)./caps;
% put Vin to the initial voltages
pres_Vc(n)=Vin.*A(i,n-1);
% solving
next_Vc=E\pres_Vc;
% very small charge value
% next_Vc(n)
% changing the initial voltages
pres_Vc=next_Vc;
% extracting the final values
cnt=old+i;
all=pres_Vc(1:n-1);
Vall(cnt,:)=[cnt all'];
Q(cnt,:)=[cnt, pres_Vc(n)];
end
old=j*i;
end
% displaying the flying voltages
subplot(2,1,1)
stairs(Vall(:,1), Vall(:,2:n-1),'LineWidth',1 );
axis normal
grid
xlabel('iteration number');
ylabel('V1, V2, V3  [V]')
% displaying the output voltage
subplot(2,1,2)
stairs(Vall(:,1), Vall(:,n),'LineWidth',1 );
axis normal
```

```matlab
grid
xlabel('iteration number');
ylabel('Vout [V]')
% displaying the decaying charge
close all;
subplot(2,1,1)
stairs(Q(:,1), Q(:,2),'LineWidth',1 );
axis normal
grid
xlabel('iteration number');
ylabel('Charge [C]')
% displaying the charge locus
nt=5; x=0;
angle=2*pi/nt;
fin=numel(Q(:,2));
for shift=1:nt
   x(shift:nt:fin)=(shift-1).*angle;
end
close all;
theta=x';
rho=abs(Q(:,2));
polar(theta,rho,'-*b')
```

*C code for the PSIM script block:*

```c
//  ***   ARRAYS FOR 1/8 Vdd
// the EXB code {0 0 0 1}
const unsigned char SW18_1[] =
{//  S1 S2 S3 S4 S5 S6 S7 S8 S9 S10 S11 S12
    0, 0, 0, 0, 1, 1, 0, 0, 0, 1,  1,  0};
// the EXB code {0 0 1 -1}
const unsigned char SW18_2[] =
{//  S1 S2 S3 S4 S5 S6 S7 S8 S9 S10 S11 S12
    0, 0, 0, 0, 1, 1, 1, 0, 0, 0,  0,  1};
// the EXB code {1 -1 -1 -1}
const unsigned char SW18_3[] =
{//  S1 S2 S3 S4 S5 S6 S7 S8 S9 S10 S11 S12
    1, 0, 0, 1, 0, 0, 0, 0, 1, 0,  0,  1};
// the EXB code {0 1 -1 -1}
const unsigned char SW18_4[] =
{//  S1 S2 S3 S4 S5 S6 S7 S8 S9 S10 S11 S12
    0, 1, 0, 0, 0, 1, 0, 0, 1, 0,  0,  1};

//  ***   ARRAYS FOR 2/8 Vdd
// the EXB code {0 0 1 0}
const unsigned char SW28_1[] =
{//  S1 S2 S3 S4 S5 S6 S7 S8 S9 S10 S11 S12
    0, 0, 0, 0, 1, 1, 1, 0, 0, 0,  1,  0};
// the EXB code {1 -1 -1 0}
const unsigned char SW28_2[] =
{//  S1 S2 S3 S4 S5 S6 S7 S8 S9 S10 S11 S12
    1, 0, 0, 1, 0, 0, 0, 0, 0, 1,  0,  1};
```

```
//  the EXB code {0 1 -1 0}
const unsigned char SW28_3[] =
{//  S1  S2  S3  S4  S5  S6  S7  S8  S9  S10 S11 S12
     0,  1,  0,  0,  0,  1,  0,  0,  0,  1,  0,  1};

//  ***    ARRAYS FOR 3/8 Vdd
//  the EXB code {1 -1 0 -1}
const unsigned char SW38_1[] =
{//  S1  S2  S3  S4  S5  S6  S7  S8  S9  S10 S11 S12
     1,  0,  0,  0,  1,  0,  0,  0,  1,  0,  0,  1};
//  the EXB code {0 1 0 -1}
const unsigned char SW38_2[] =
{//  S1  S2  S3  S4  S5  S6  S7  S8  S9  S10 S11 S12
     0,  1,  0,  0,  0,  1,  1,  0,  0,  0,  0,  1};
//  the EXB code {0 0 1 1}
const unsigned char SW38_3[] =
{//  S1  S2  S3  S4  S5  S6  S7  S8  S9  S10 S11 S12
     0,  0,  0,  0,  1,  1,  0,  1,  0,  0,  1,  0};
//  the EXB code {1 -1 -1 1}
const unsigned char SW38_4[] =
{//  S1  S2  S3  S4  S5  S6  S7  S8  S9  S10 S11 S12
     1,  0,  0,  1,  0,  0,  0,  0,  0,  1,  1,  0};

//  ***    ARRAYS FOR 4/8 Vdd
//  the EXB code {1 -1 0 0}
const unsigned char SW48_1[] =
{//  S1  S2  S3  S4  S5  S6  S7  S8  S9  S10 S11 S12
     1,  0,  0,  0,  1,  0,  0,  0,  0,  1,  0,  1};
//  the EXB code {0 1 0 0}
const unsigned char SW48_2[] =
{//  S1  S2  S3  S4  S5  S6  S7  S8  S9  S10 S11 S12
     0,  1,  0,  0,  0,  1,  1,  0,  0,  0,  1,  0};

//  ***    ARRAYS FOR 5/8 Vdd
//  the EXB code {1 0 -1 -1}
const unsigned char SW58_1[] =
{//  S1  S2  S3  S4  S5  S6  S7  S8  S9  S10 S11 S12
     1,  1,  0,  0,  0,  0,  0,  0,  1,  0,  0,  1};
//  the EXB code {1 -1 0 1}
const unsigned char SW58_2[] =
{//  S1  S2  S3  S4  S5  S6  S7  S8  S9  S10 S11 S12
     1,  0,  0,  0,  1,  0,  0,  0,  0,  1,  1,  0};
//  the EXB code {0 1 0 1}
const unsigned char SW58_3[] =
{//  S1  S2  S3  S4  S5  S6  S7  S8  S9  S10 S11 S12
     0,  0,  1,  0,  0,  1,  0,  0,  0,  1,  1,  0};
//  the EXB code {1 -1 1 -1}
const unsigned char SW58_4[] =
{//  S1  S2  S3  S4  S5  S6  S7  S8  S9  S10 S11 S12
     1,  0,  0,  0,  1,  0,  1,  0,  0,  0,  0,  1};
```

```
// ***   ARRAYS FOR 6/8 Vdd
// the EXB code {1 0 -1 0}
const unsigned char SW68_1[] =
{//  S1  S2  S3  S4  S5  S6  S7  S8  S9  S10  S11  S12
     1,  1,  0,  0,  0,  0,  0,  0,  0,  1,   0,   1};
// the EXB code {1 -1 1 0}
const unsigned char SW68_2[] =
{//  S1  S2  S3  S4  S5  S6  S7  S8  S9  S10  S11  S12
     1,  0,  0,  0,  1,  0,  1,  0,  0,  0,   1,   0};
// the EXB code {0 1 1 0}
const unsigned char SW68_3[] =
{//  S1  S2  S3  S4  S5  S6  S7  S8  S9  S10  S11  S12
     0,  0,  1,  0,  0,  1,  1,  0,  0,  0,   1,   0};

// ***   ARRAYS FOR 7/8 Vdd
// the EXB code {1 0 0 -1}
const unsigned char SW78_1[] =
{//  S1  S2  S3  S4  S5  S6  S7  S8  S9  S10  S11  S12
     1,  1,  0,  0,  0,  0,  1,  0,  0,  0,   0,   1};
// the EXB code {1 0 -1 1}
const unsigned char SW78_2[] =
{//  S1  S2  S3  S4  S5  S6  S7  S8  S9  S10  S11  S12
     1,  1,  0,  0,  0,  0,  0,  0,  0,  1,   1,   0};
// the EXB code {1 -1 1 1}
const unsigned char SW78_3[] =
{//  S1  S2  S3  S4  S5  S6  S7  S8  S9  S10  S11  S12
     1,  0,  0,  0,  1,  0,  0,  1,  0,  0,   1,   0};
// the EXB code {0 1 1 1}
const unsigned char SW78_4[] =
{//  S1  S2  S3  S4  S5  S6  S7  S8  S9  S10  S11  S12
     0,  0,  1,  0,  0,  1,  0,  1,  0,  0,   1,   0};

// *** MAIN FOR Vo=3/8Vdd ***

unsigned char a=0;
unsigned char b=0;
unsigned char num=0;
static unsigned char cnt=0;

b = in[0]; // clock pulse
if ( b !=a ) // each transition
{ a = b;   cnt++;
if ( cnt >4 ) { cnt = 1; }
switch (cnt)
{ case 1:
   for (num = 0; num < 12; num++)
        {out[num] = SW38_1[num];}
   break;
   case 2:
   for (num = 0; num < 12; num++)
        {out[num] = SW38_2[num]; }
```

```
       break;
    case 3:
    for (num = 0; num < 12; num++)
          {out[num] = SW38_3[num];}
    break;
    case 4:
    for (num = 0; num < 12; num++)
          {out[num] = SW38_4[num];}
    break; }
 }
// end for Vo=3/8Vdd
```